\def\real{{\tt I\kern-.2em{R}}}   
\def\nat{{\tt I\kern-.2em{N}}}    
\def\snat{{\rm I\kern-.2em{N}}}    
\def\eps{\epsilon}                 
\def\realp#1{{\tt I\kern-.2em{R}}^#1}
\def\natp#1{{\tt I\kern-.2em{N}}^#1}
\def\hyper#1{\,^*\kern-.2em{#1}}
\def\hy#1{\,^*\kern-.2em{#1}}

\def\St#1{{\tt st}#1}
\def\st#1{{\tt st}(#1)}
\def\hyperreal{{^*{\real}}}
\def\hyperrealp#1{{\tt ^*{I\kern-.2em{R}}}^#1} 
\def\hypernat{{^*{\nat }}}
\def\hypernatp#1{{{^*{{\tt I\kern-.2em{N}}}}}^#1} 
\def\eskip{\hskip.25em\relax}

\def\Hyper#1{\hyper {\eskip #1}}
\def\Hy#1{\hyper {\eskip #1}}
\def\leaderfill{\leaders\hbox to 1em{\hss.\hss}\hfill}
\def\srealp#1{{\rm I\kern-.2em{R}}^#1}

\def\power#1{{{\cal P}(#1)}}
\def\iff{\leftrightarrow}
\def\qed{{\vrule height6pt width3pt depth2pt}\ }
\def\pars{\par\smallskip}
\def\parm{\par\medskip}
\def\r#1{{\rm #1}}

\def\ref#1{$^{#1}$}

\def\sig{{^\sigma}}
\def\m@th{\mathsurround=0pt}
\def\rightarrowfill{$\m@th \mathord- \mkern-6mu \cleaders\hbox{$\mkern-2mu 
\mathord- \mkern-2mu$}\hfil \mkern-6mu \mathord\rightarrow$}
\def\leftarrowfill{$\mathord\leftarrow
\mkern -6mu \m@th \mathord- \mkern-6mu \cleaders\hbox{$\mkern-2mu 
\mathord- \mkern-2mu$}\hfil $}
\def\noarrowfill{$\m@th \mathord- \mkern-6mu \cleaders\hbox{$\mkern-2mu 
\mathord- \mkern-2mu$}\hfil$}
\def\orgate{$\bigcirc \kern-.80em \lor$}
\def\andgate{$\bigcirc \kern-.80em \land$}
\def\inverter{$\bigcirc \kern-.80em \neg$}
\magnification=1000         
\tolerance 10000                
\baselineskip=14pt
\font\eightsl=cmsl8
\hoffset=.25in 
\hsize 6.00 true in 
\vsize 8.75 true in 

\def\id{\par\hangindent2\parindent\textindent}
\def\textindent#1{\indent\llap{#1}}
\pageno=1
\headline={\ifnum\pageno= 1\hfil {\quad}
\hfil\else\ifodd\pageno\rightheadline \else\leftheadline\fi\fi} 
\def\rightheadline{\eightsl \hfil Nonstandard Analysis Simplified  \hfil} 
\def\leftheadline{{} \eightsl \hfil Nonstandard Analysis Simplified \hfil} 
\voffset=1\baselineskip
 
\def\id{\par\hangindent2\parindent\textindent}
\def\textindent#1{\indent\llap{#1}}
\font\large=cmdunh10 at 18truept
\font\next=cmr12 at 14truept

{\quad}
\vskip 1.50in
\centerline{\large NONSTANDARD ANALYSIS}\bigskip
\centerline{\next - A SIMPLIFIED APPROACH - }
\vskip 1.5in
\centerline{\next ROBERT A. HERRMANN}\vfil\eject
{\quad}\vskip 2.5in 
\centerline{Copyright \copyright\  2003 by Robert A. Herrmann}
\vfil\eject
\centerline{}
\centerline{\bf CONTENTS}
\bigskip
\indent {\bf Chapter 1}\parm
\indent {\bf Filters} {\eightsl Ultrafilters, Cofinite Filter, Principle Ultrafilters, \par Free Ultrafilters} \leaderfill {\bf 6}.\par 
\bigskip
\indent {\bf Chapter 2}\parm
{\bf A Simple Nonstandard Model for Analysis} {\eightsl Equivalence\par Classes of Sequences, Totally Ordered Field of Equivalence Classes,\par The Hyper-extension of Sets and Relations, The Standard Object\par Generator}\leaderfill {\bf 9}.\par
\bigskip                                                      
\indent {\bf Chapter 3}\parm
{\bf Hyper-set Algebra\hfil}\par 
{\bf Infinite and Infinitesimal Numbers} {\eightsl The Behavior of the Hyper\par and Standard Object Generators; Infinitesimals $\mu (0)$, The Infinite and\par Finite Numbers, *-Transform Process, Maximum Ideal $\mu (0)$}\leaderfill {\bf 15.}
\bigskip
\indent {\bf Chapter 4}\parm
{\bf Basic Sequential Convergence} {\eightsl Bounded Sequences, Convergent\par Sequences, Accumulation Point, Subsequences, Cauchy Criterion,\par Nonstandard Characteristics and Examples} \leaderfill {\bf 23}\bigskip

\indent {\bf Chapter 5}\parm
{\bf Advanced Sequential Convergence} {\eightsl Double Sequences, Iterated \par Limits, Upper and Lower Limits, Nonstandard Characteristics and\par Examples} \leaderfill {\bf 28.}
\bigskip
\indent {\bf Chapter 6}\parm
{\bf Basic Infinite Series} {\eightsl Hyperfinite Summation, Standard Results, \par Nonstandard Characteristics and Examples}\leaderfill {\bf 33.}
\bigskip
\indent {\bf Chapter 7}\parm
{\bf An Advance Infinite Series Concept} {\eightsl Multiplying of Infinite\par  Series, Nonstandard Characteristics and Examples} \leaderfill {\bf 37.}
\bigskip
\indent {\bf Chapter 8}\parm
{\bf Additional Real Number Properties} {\eightsl Interior, Closure,\par Cluster, Accumulation, Isolated Points, Boundedness, Compactness,\par Nonstandard Characteristics} \leaderfill {\bf 41.}
\bigskip \vfil\eject   
\indent {\bf Chapter 9}\parm
{\bf Basic Continuous Function Concepts} {\eightsl All Notions\hfil\par Generalized to Cluster Points, One-sided Limits and Continuity, Sum,\par Product and Composition of Continuous Functions, Extreme and\par Intermediate Value Theorems, Nonstandard Characteristics} \leaderfill {\bf 46.}
\bigskip
\indent {\bf Chapter 10} \parm
{\bf Slightly Advanced Continuous Function Concepts} {\eightsl Non-\hfil \par standard Analysis and Bolzano's Product Theorem, Inverse Images of \par Open Sets, Additive Functions, Uniform Continuity, Extensions of\par Continuous Functions} \leaderfill {\bf 51.}
\bigskip 
\indent {\bf Chapter 11} \parm
{\bf Basic Derivative Concepts} {\eightsl Nonstandard Characteristics\hfil \par for Finite and Infinite Derivatives at Cluster Points, The Infinitesimal\par Differential, The Fundamental Theorem of Differentials, Order Ideals,\par Basic Theorems, Generalized Mean Value Theorem, L'Hospital's\hfil\par Rule} \leaderfill {\bf 55.}
\bigskip
\indent {\bf Chapter 12}\parm
{\bf Some Advanced Derivative Concepts} {\eightsl Nonstandard Analysis\par and the nth-Order Increments, nth-Order Ideals, Continuous \par Differentiability, Uniformly Differentiable, the Darboux Property, and \par Inverse Function Theorems} \leaderfill {\bf 61.}
\bigskip
\indent {\bf Chapter 13}\parm
{\bf Riemann Integration} {\eightsl The Simple Partition, Fine Partitions,\par Upper and Lower Sums, Upper and Lower Hyperfinite Sums,\par The Simple Integral, The Equivalence of The Simple Integral \par and The Riemann Integral, The Basic Integral Theorems and\par How The Generalized and Lebesgue Integral Relate to Fine \par Partitions.} \leaderfill {\bf 68 .}
\bigskip
\indent {\bf Chapter 14}\parm
{\bf What Does the Integral Measure?} {\eightsl Additive Functions and\par The Rectangular Property}\leaderfill {\bf 75.}
\bigskip
\indent {\bf Chapter 15}\parm
{\bf Generalizations} {\eightsl Metric and Normed Linear Spaces} \leaderfill {\bf 78.}
\bigskip
{\bf Appendix} {\eightsl Existence of Free Ultrafilters, Proof of *-Transform\par  Process} \leaderfill {\bf 80.}\bigskip
\line{\indent {\bf References\leaderfill 82.}}\bigskip
\vfil\eject
{\quad}\vskip 2.0in
{\leftskip 0.5in \rightskip 0.5in \noindent Although the material in this book is copyrighted by the author, it may be reproduced in whole or in part by any method, without the payment of any fees, as long as proper credit and identification is given to the author of the material reproduced.\par
\centerline{Disclaimer}
\noindent The material in this monograph has not been independently edited nor independently verified for correct content. It may contain various typographical or content errors. Corrections will be made only if such errors are significant. A few results, due to the simplicity of this approach, may not appear to be convincingly established. However, a change in our language or in our structure would remove any doubts that the results can be established rigorously. \par}\par  
\vfil\eject
\centerline{\bf 1. FILTERS}
\parm
For over three hundred years, a basic question about the calculus remained unanswered. Do the infinitesimals, as conceptional understood by Leibniz and Newton, exist as formal mathematical objects? This question was answer affirmatively by Robinson (1961) and the subject termed ``Nonstandard Analysis'' (Robinson, 1966) was introduced to the scientific world. As part of this book, the mathematical existence of the infinitesimals is established and their properties investigated and applied to basic real analysis notions.\pars
 I intend to write this ``book'' informally and it's certainly about time that technical books be presented in a more ``friendly'' style. What I'm not going to do is to present an introduction filled with various historical facts and self-serving statements; statements that indicate what an enormous advancement in mathematics has been achieved by the use of Nonstandard Analysis. Rather, let's proceed directly to this simplified approach, an approach that's correct but an approach that cannot be used to analysis certain areas of mathematics that are not classified as elementary in character. These areas can be analysis but it requires one to consider additional specialized mathematical objects. Such specialized objects need only be considered after an individual becomes accustomed to the basic methods used within this simple approach. There are numerous  exciting and thrilling new concepts and results that cannot be presented using the simple approach discussed. The ``internal'' objects, objects that ``bound'' sets that represent ``concurrent'' relations and saturated models are for your future consideration. The main goal is to present some of the basic nonstandard results that can be obtained without investigating such specialized objects. \pars
I'll present a complete ``Proof'' for a stated result. However, one only needs to have confidence that the stated ``Theorem'' has been acceptably established. Indeed, if you simply are interested in how these results parallel the original notions of the ``infinitesimal'' and the like, you need not bother to read the proofs at all. \pars
There's an immediate need for a few set-theoretic notions. We let $\nat$ be the set of all natural numbers, which includes the zero as the first one. I assume that you understand some basic set-theoretic notation. {\bf Further, throughout this first chapter, $X$ will always denote a nonempty set.} Recall that for a given set $X$ the set of all subsets of $X$ exists and is called the {{\bf power set}}. It's usually denoted by the symbol $\power {X}.$ For example, let $X = \{0,1,2\}$. The power set of $X$ contains 8 sets. In particular, $\power {X} = \{\emptyset, X, \{0\},\{1\},\{2\},\{0,1\},\{0,2\}, \{1,2\}\},$ where $\emptyset$ denotes the {{\bf empty set}}, which can be thought of as a set which contains ``no members.'' \parm
{\bf Definition 1.1.} ({{\bf The Filter.}}) (The symbol $\subset$ means ``subset'' and includes the possible equality of sets.) A (nonempty) $\emptyset \not= {\cal F}\subset \power {X}$ is called a (proper) {{\bf filter on $X$}} if and only if \pars
\indent\indent (i) for each $A,\ B \in {\cal F}$, $A \cap B \in {\cal F}$;\pars
\indent\indent (ii) if $A \subset B \subset X$ and $A \in {\cal F},$ then  $B \in {\cal F}$; \pars
\indent\indent (iii) $\emptyset \notin \cal F$.\parm  
{\bf Example 1.2.} (i) Let $\emptyset \not= A \subset X$. Then $[A]\uparrow $ is the set of all subsets of $X$ that contain $A,$ or, more formally, $[A]\uparrow = \{ x \mid x \subset X\ {\rm and}\ A \subset x\},$ is a filter on $X$ called the {{\bf principal filter}}.\parm
How to properly define what one means by a ``{finite}'' set has a long history. But as {Suppes} states ``The common sense notion is that a set is finite just when it has ``$m$'' members for some non-negative integer $m.$ [You can use our set $\nat$ to get such an ``$m$.''] This common sense idea is technically sound . . . .'' (1960, p. 98). The notion of finite can also be related to {\bf constants} that name members of a set and the formal expression that characterizes such a set in terms of the symbols ``='' and ``$\lor$'' (i.e.``or''). Notice that the empty set is a ``finite set.'' This association to the natural numbers is denoted by subscripts that actually represent the range values of a function. I also assume certain elementary properties of the finite sets and those sets that are not finite. For $X$, let $\emptyset \not= {\cal B} \subset \power {X}$ have the {{\bf finite intersection property}}. This means that each $B_i \in {\cal B}$ is nonempty and that the ``intersection'' of all of the members of any other finite subset of $\cal B$ is not the empty set. Given such a ${\cal B} \not= \power {X}$, then we can generate the ``smallest'' filter that contains $\cal B.$ This is done by first letting ${\cal B}'$ be the set of all subsets of $\power {X}$ formed by taking the intersection of each nonempty finite subset of ${\cal B},$ where the intersection of the members of a set that contains but one member is that one member. Let's consider some basic mathematical abbreviations. The formal symbol {$\land$} means ``and'' and the formal ``quantifier'' {$\exists$} means ``there exists such and such.'' Now take a member $B$  of ${\cal B}'$ and build a set composed of all subsets of $X$ that contain $B.$ Now do this for all members of ${\cal B}'$ and gather them together in a set to get the set $\langle {\cal B}\rangle.$ Hence, $\langle {\cal B}\rangle$ is the set of all subsets $x$ of $X$ such that there exists some set $B$ such that $B$ is a member of ${\cal B}'$ and $B$ is a subset of $x$. Formally,  $\langle {\cal B}\rangle = \{x \mid (x \subset X)\land (\exists B((B \in {\cal B}') \land (B \subset x))\}.$ I have  ``forced'' $\langle {\cal B}\rangle$ to contain $\cal B$ and to have the necessary properties that makes it a filter.\parm
There exists a very significant ``filter'' $\cal C$ on infinite $X$ defined by the notion of not being finite. This object, once we have shown that it is a filter on $X$, is called the {{\bf cofinite}} filter. [Cofinite means that the relative complement is a finite set].  \parm
{\bf Definition 1.3.} ({{\bf QED and iff.}}) Let ${\cal C} = \{x \mid (x \subset X) \land (X - x)\ {\rm is\ finite})\}.$ Also I will use the symbol \qed  for the statement ``QED,'' which indicates the end of the proof. Then ``iff'' is an abbreviation for the phrase ``if and only if.''\parm
{\bf Theorem 1.4.} {\it The set $\cal C$ is a filter on infinite $X$ and, the intersection of all members of ${\cal C},$ $\bigcap\{F\mid F \in {\cal C}\} = \emptyset$.}\pars
Proof. Since $X \not= \emptyset,$ then there is some $a \in X.$ Further, $X - (X -\{a\}) = \{a\}$ implies that ${\cal C} \not= \emptyset.$ So, assume that $A,B \in \cal C.$ Then since $X -(A\cap B) = (X - A)\cup (X -B)$, $X - (A\cap B)$ is a finite subset of $X.$
Thus, since $X -(X -(A\cap B)) = A\cap B,$ then $A \cap B \in \cal C.$ Now suppose  $A \subset C \subset X.$ Then $X -C \subset X -A$. Thus, because $X -A$ is finite, then $X -C$ is finite. Hence, $C \in \cal C.$ 
Also, since $X$ is infinite, then $X -\emptyset = X$ implies that $\emptyset \notin \cal C.$ Consequently, $\cal C$ is a filter on $X.$ \pars
Now observe that $K = X - \bigcap\{F\mid F \in {\cal C}\} = \bigcup \{X - F\mid F \in {\cal C}\} = X,$ for if $a \in X,$ then $X - \{a\} \in \cal C$ implies that $X -(X -\{a\}) = \{a\} \subset K.$ Hence, we must have that $\cap\{F\mid F \in {\cal C}\}= \emptyset$. \qed\parm
I mention that $\cal C$ is also called the {Fr\'echet} filter. It turns out that we are mostly interested in a maximum filter that contains $\cal C.$ \parm
{\bf Definition 1.5.} ({{\bf Ultrafilter.}}) A filter $\cal U$ on $X$ is called an {{\bf ultrafilter}} iff whenever there's a filter $\cal F$ on X such that 
${\cal U} \subset {\cal F}$, then ${\cal U} = {\cal F}.$ \parm
Prior to showing that ultrafilters exist, let's see if they have any additional useful properties.\parm
{\bf Theorem 1.6.} {\it Suppose that $\cal U$ is an ultrafilter on $X$. If $A  \cup B \in \cal U,$ then $A \in \cal U$ or $B \in \cal U.$}\pars
Proof. Let $A \notin \cal U$ and $B \notin \cal U$ but $A \cup B \in \cal U$. Let ${\cal G} = \{x \mid (x \subset X)\land (A\cup x \in {\cal U}\}.$ We show that $\cal G$ is a filter on $X$. [Note: We now begin to use variables such as $x,y,z$ etc. as mathematical variables representing members of sets. These symbols are used in two context, however. The other context is as a variable in our formal logical expressions.] Let ${x,y, \in \cal G}.$ Then $A \cup x,\ A \cup y \in \cal U.$ Hence, $(A \cup x)\cap (A \cup y) = A \cup (x \cap y) \in \cal U$ implies that $x\cup y \in {\cal G}.$ Now suppose that $x\in \cal G$ and $x\subset y \subset X.$ Then $A \cup x \subset A \cup y $ implies that $A \cup y \in \cal U$. Hence, $y \in \cal G.$ Also $A \cup \emptyset \notin \cal U$ implies that $\emptyset \notin \cal G.$ Thus $\cal G$ is a filter on $X$. \pars
Let $C \in \cal U$. Then $C \subset A \cup C \in \cal U$ implies that $C \in \cal G.$ Therefore, ${\cal U} \subset \cal G$. But, $B\in \cal G$ implies that ${\cal U} \not= {\cal G.}$ This contradicts the maximum aspect for $\cal U$. \qed\parm 

{\bf Theorem 1.7.} {\it Let $\cal F$ be a filter on $X.$ Then $\cal F$ is an ultrafilter iff for each $A \subset X$, either $A \in \cal U$ or $X-A \in \cal U$, not both.}\pars
Proof. Assume $\cal F$ is a ultrafilter. Then $X = X \cup (X - A)$ implies that either $A \in \cal U$ or $X - A \in \cal U.$ Both $A$ and $X-A$ cannot be members of $\cal U$, for if they were then $A \cap (X -A) =\emptyset \in \cal U$; a contradiction.\pars
Conversely, suppose that for each $A \subset X$, either $A \in \cal F$ or $X-A \in \cal F.$ Let $\cal G$ be filter on $X$ such that ${\cal F}\subset {\cal G} .$ Let $A \in \cal G$. Then $X-A \notin \cal G$ since $\cal G$ is a filter. Thus, $X-A \notin \cal F.$ Hence, $A \in \cal F$. Thus ${\cal G} \subset \cal F$ implies that ${\cal G} = {\cal F}$. \qed\pars  
Given any filter $\cal F$ on $X$ a major question is whether there exists an ultrafilter $\cal U$ on $X$ such that ${\cal F}\subset {\cal U}$ The answer to this question can take on, at least, two forms. The next result states that such ultrafilters always exist. The proof in the appendix uses a result, Zorn's Lemma, that is equivalent to the Axiom of Choice. The Axiom of Choice, although it's consistent with the other axioms of set theory, may not be ``liked'' by some.   
There's an axiom that is also consistent with the other usual axioms of set theory that is weaker than the Axiom of Choice. What it states is that such an ultrafilter always exits. So, you can take your pick. \parm
{\bf Theorem 1.8.} {\it Let $\cal F$ be a filter on $X$. Then there exists an ultrafilter $\cal U$ on $X$ such that ${\cal F}\subset {\cal U}$.}\pars
Proof. See the appendix.\qed\pars
A natural study is to see if we can partition the set of all ultrafilters defined on $X$ into different categories. And, why don't we use the symbol ${\cal F}_X$ [resp. ${\cal U}_X$] to always denote a filter [resp. ultrafilter] on $X$. It turns out there are two basic types of ${\cal U}_X,$ the principal ones and those that contain ${\cal C}_X.$ \pars
{\bf Theorem 1.9.} {\it Let $p \in X.$ Then $[p]\uparrow$ is an ${\cal U}_X.$}\pars
Proof. Let nonempty $A \subset X.$ Then either $p \in  A$ or $p \in (X-A)$ and  not both. Thus $A \in [p]\uparrow$ or $(X - A) \in [p]\uparrow.$  Hence, by Theorem 1.3, $[p]\uparrow$ is an ${\cal U}_X.$\parm

{\bf Theorem 1.10.} {\it Assume that ${\cal U}_X$ is not a principal ultrafilter. Then ${\cal C}_X \subset {\cal U}_X.$} \pars
Proof. Let arbitrary nonempty finite $\{p_0,\ldots,p_k\} \subset X.$ Since ${\cal U}_X$ is non-principal, then ${\cal U}_X \not= [p_i]\uparrow,\ i = 0,\ldots,k.$ Hence, for each $i= 0,\ldots, k$ there exists some $A_i \subset X$ such that $A_i \in {\cal U}_X$ and $p_i \notin A_i$. For, otherwise, if $p_i \in A_i$ for any $A_i \in {\cal U}_X$, then $[p_i]\uparrow \subset {\cal U}_X$ (they are =.)  Consequently, $\{p_0,\ldots,p_k\} \cap (A_0\cap \cdots \cap A_k) = \emptyset.$ However, $(A_0\cap \cdots \cap A_k)\in {\cal U}_X.$ Therefore,
$\{p_0,\ldots,p_k\} \notin {\cal U}_X.$ Theorem 1.3 implies that $X- \{p_0,\ldots, p_k\} \in {\cal U}_X.$ Thus, ${\cal C}_X \subset {\cal U}_X$. \qed \pars  
Non-principal ultrafilters are also called {{\bf free}} ultrafilters. This comes from Theorems 1.4 and 1.10 which imply that ${\cal U}_X$ is  free iff $\bigcap\{F\mid F \in {\cal U}_X\} = \emptyset.$ Also, another characterization is that ${\cal U}_X$ is free iff there does not exist a nonempty finite $F \subset X$ such that $F \in {\cal U}_X.$ If we let $X = \nat,$ then there are a lot of free ${\cal U}_X$. Unless otherwise stated, the free ultrafilter that's used will not affect any of the stated results.  
\vfil\eject
\centerline{\bf 2. A SIMPLE NONSTANDARD MODEL FOR ANALYSIS}
\parm
 We let {$\real$} denote the {{\bf real numbers}}. The set $\real$ uses various operators and relations to obtain results within analysis. For this simplified approach, most of what we need is defined from the basic addition $+$, multiplication $\cdot$, total order $\leq$ properties and few other ones accorded $\real.$  For convenience, I denote this fact by the {{\bf structure}} notation $\langle \real,+,\cdot,\leq , \Phi_i\rangle,$ where the $\Phi_i$ are any other relations one might consider for $\real$ whether definable from the basic relations or not. {\it Further, it is always understood that each structure includes the $=$ relation, which is but set-theoretic equality or identity for members of $\real.$} I'll have a little more to write about how these ``numbers'' should be viewed later. But, first to our construction. Let $\real^{\snat}$ represent the {{\bf set of all sequences with domain $\nat$ and range values (images) in $\real.$}} Of course, sequences are functions, (maps, mappings, etc.) that are often displayed as a type of ``ordered'' set in the form $\{s_0,s_1,s_2,\ldots \}.$ You can define binary operators $+$ and $\cdot,$ among others, for sequences by simply taking any two $f,g \in \real^\snat$ and defining $f + g = h$ to be the sequence 
$h$ where the values of $h$ are $h(n) = f(n) + g(n)$ and $f\cdot g = fg = k$ 
to be the sequence $k$ where the values of $k$ are $k(n) = f(n)g(n)$ for each $n \in \nat.$ This forms, at the very least, what is called a {\bf ring} with unity. What I'll do later is to show that there's a subset of $\real^\snat$ that ``behaviors''  like the real numbers, with respect to the defined relations, and we'll us this subset as if it is the real numbers. In all the follows, ${\cal U} = {\cal U}_\snat$ will always be a free ultrafilter and the symbol $U$ is used to represent members of ${\cal U}$. Now to make things symbolically simple {{\bf capital letters from the beginning of the alphabet}} $A, B, C, \ldots$ will always denote members of $\real^\snat.$ Also, we usually use the subscript notation for the images. Now let us begin our construction of a nonstandard model for real analysis. \parm
{\bf Definition 2.1.} ({\bf Equality in} ${\cal U}$) Let $A,B \in \real^\snat$. Define $A =_{\cal U} B$ iff $\{n\mid A_n = B_n \}=U \in \cal U.$ (The set of all $\nat$ such that the values of the sequences $A$ and $B$ are equal.)\parm
It has been said that the most important binary relation within mathematics is the {{\bf equivalence}} relation. This relation, in general, behaves like $=$ except that you may not be allowed to ``substitute'' one equivalent object for another. Recall that for a set $X$ a binary relation $R$ is an equivalence relation on $X$ iff it has the following properties. For each $x,y,z \in X$,  (i) $xRx$ {(reflexive property)}; (ii) if $xRy$, then $yRx$ {(symmetric property)}; [Note that if this holds, then $xRy$ iff $yRx.$] (iii) if $xRy$ and $yRz$, then $xRz$ {(transitive property)}. Hence, it is almost an ``equality.''  \parm 
{\bf Theorem 2.2.} {\it The relation $=_{\cal U}$ is an equivalence relation on $\real^\snat.$}\pars
Proof. Of course, properties of the $=$ for members of $\real$ are used. First, notice that $\{n\mid A_n = A_n\} = \nat \in \cal U$ for any $A \in \real^\snat.$ Thus, the relation is reflexive.\pars
Clearly, for any $A,B \in \real,$ if $\{n\mid A_n = B_n\} \in \cal U$, then $\{n\mid B_n = A_n\} \in \cal U$.\pars
Finally, suppose that $A,B,C \in \real^\snat$ and $A =_{\cal U} B$ and $B =_{\cal U} C.$ Hence, $\{n\mid A_n = B_n\} \in \cal U$ and  $\{n\mid B_n = C_n\} \in \cal U.$ The word ``and'' implies, since $\cal U$ is a filter, that $\{n\mid A_n = B_n\} \cap \{n\mid B_n = C_n\} \in \cal U.$ Of course, this ``intersection'' need not give all the values of $\nat$ that these three sequences have in common, but that does not matter since the ``superset'' property for a filter implies from the result  
$$ \{n\mid A_n = B_n\} \cap \{n\mid B_n = C_n\}\subset \{n\mid A_n = C_n\}$$
that $\{n\mid A_n = C_n\} \in \cal U.$ \qed \parm
[Note: In the above ``proof,'' the two step process of getting the common members by the ``intersection'' and using the superset property is a major proof method.]\parm
{\bf Definition 2.3.} ({\bf Equivalence classes.}) We now use the relation $=_{\cal U}$ to define actual subsets of $\real^\snat$. For each $A \in \real^\snat$, let the set $[A] = \{x\mid (x \in \real^\snat) \land(x =_{\cal U} A)\}.$\parm
It is easy to show that for each $A, B \in \real^\snat$, either $[A] = [B]$ or $[A] \cap [B] = \emptyset.$ (The ``='' here is the set-theoretic equality.) Further, $\real^\snat = \bigcup \{[x]\mid x \in \real^\snat\}.$ That is the set $\real^\snat$ is completely {\bf partitioned} (separated into, broken up into) these non-overlapping nonempty sets. Because of these properties, we can use any member of the set $[A]$ to generate the set. That is if $B,C \in [A],$ then $[A] =[B] =[C].$ As to notation, when I'm not particular interested in a sequence that generates the equivalence class, I'll denote them by {{\bf lower case}} letters $a,b,c,\ldots$. \pars
Denote the set of all of these equivalence classes by $\hyperreal$ and call this set the set of all {{\bf hyperreal numbers}}. (The $\hyper {}$ is often translated as ``{{\bf hyper}}.'') Consequently, $\hyperreal = \{[A]\mid A \in \real^\snat\}.$ After various relations are defined on $\hyperreal,$ the resulting ``structure'' is generally called an {{\bf ultrapower}}. Indeed, it's this ultrapower that will act as our nonstandard model for portions of real analysis. There's still a lot of work to do to turn $\hyperreal$ into a such a model, but to motive this work I'll simply mention that if you take a sequence $s$ that converges in the normal calculus sense to 0, then $[s]$ is one of our {infinitesimals}. What will be done, after the ultrapower model is constructed, is to ``embed'' $\langle \real, +, \cdot, \leq, \Phi_i \rangle$ into the ultrapower so that comparisons can be easily made between the ``standard'' objects that represent the properties of the actual real numbers and other objects in the ultrapower. The notation $\langle \real, +, \cdot, \leq, \Phi_i \rangle$ identifies the {{\bf carrier}}, $\real,$ as well as certain specialized relations defined for (on) the carrier.\pars
There are two approaches to analyze this ultrapower, a direct and tedious method, and a method that uses notions from Mathematical Logic. Once  everything is constructed and the embedding is secured, then the embedded objects become our {{\bf standard}} objects. The set of {{\bf nonstandard}} objects is the remainder of the ultrapower.  \parm
{\bf Definition 2.4.} ({\bf Addition and multiplication for the $\hyperreal$.}) 
Consider any $a,b,c\in \hyperreal.$ Define $a \hyper {+}\ b = c$ iff $\{n\mid A_n + B_n = C_n \} \in \cal U.$ [Note: such definitions assume that you have selected some sequences $A_n \in a,\ B_n \in b, \ C_n \in c.$ Now define
$a \hyper {\cdot}\ b = c$ iff $\{n \mid (A_n)\cdot(B_n)= C_n\} \in \cal U.$\parm
Whenever such definitions are made by taking members of a set that contains more than one member it is always necessary to show that they are {{\bf well-defined}} in that the result is not dependent upon the member one chooses. The next result shows how this is done and gives insight as to how it will be done later in completely generality.\parm
{\bf Theorem 2.5.} {\it The operations defined in definition 2.4 are well-defined.}\pars
Proof. Let $[A], [D] \in a, \ [B], [F] \in b$. Notice that $\{n\mid A_n = D_n \}\in \cal U\}$ and $\{n\mid B_n = F_n \}\in \cal U$ implies that  $\{n\mid A_n = D_n \}\cap \{n\mid B_n = F_n \}\in \cal U$ and $\{n\mid A_n = D_n \}\cap \{n\mid B_n = F_n \}\subset \{n\mid A_n + B_n = D_n + F_n\}$ implies by the superset property that $\{n\mid A_n + B_n = D_n + F_n\} \in \cal U.$ Thus the $\hyper {+}$ is well-defined. (Note: Processes of this type that use filter properties that imply something is a member of a filter will be abbreviated.) In like manner, for the $\hyper {\cdot}\ .$ \qed \parm  
Thus far, the fact the $\cal U$ is an ultrafilter has not been used. But, for the structure $\langle\hyperreal, \hyper +, \hyper \cdot \rangle$ to have all the necessary mathematical ``field'' properties, this ultrafilter property is significant. That is so that the $\hy +,\ \hy \cdot$ arithmetic behaves for $\hyperreal,$ like $+,\ \cdot$ behave for real number arithmetic. \parm\vfil\eject
{\bf Theorem 2.6.} {\it For the structure $\langle \hyperreal, \hyper +,\hyper \cdot \rangle$\pars
\indent\indent {\rm (i)} $[0]$ is the additive identity;\par
\indent\indent {\rm (ii)} for each $a = [A]\in \hyperreal,$ $-a = [-A]$ is the additive inverse;\par
\indent\indent {\rm (iii)} $[1]$ is the multiplicative identity; \par
\indent\indent {\rm (iv)} If $a \not= [0],$ then there exists $b =[B] \in \hyperreal $ such that $a\hyper \cdot\, b = [1].$\pars
\indent\indent {\rm (v)} For each $n\in \nat$ if $D_n = A_n + B_n$ and $E_n = A_nB_n,$ then $[A] \hy + [B] = [D],\ [A]\hy \cdot [B] = [E].$ That is our definitions for addition and multiplication of sequences and the hyper-operators $\hy +,\ \Hy \cdot$ are compatible.}\pars
Proof.  (i) Let $[A] \hyper +\ [0] = [C].$ Considering that $\{n\mid A_n + 0_n = C_n\} \in \cal U$ and $\{n\mid A_n + 0_n = C_n\}\subset \{n\mid A_n = C_n\}\in \cal U$, then $[A]= [C].$\pars
(ii) Let $[-A] = [B].$ Then once again $\{n\mid A_n + (- A_n)= 0=0_n\}= \nat \in \cal U$ and thus $[A] \hyper+\ [-A] = [0].$  \pars
(iii) This follows in the same manner as (i).\pars
(iv) Let $[A] \not= [0].$ Then $\{n\mid A_n =0 = 0_n\}=U \notin \cal U.$ Hence, $\nat - U = \{n \mid A_n \not= 0\} \in \cal U$ since $\cal U$ is an ultrafilter. Define
$$B_n = \cases{A^{-1}_n;&if $A_n \not= 0$\cr 0;& if $A_n =0$\cr}.$$ 
Notice that $\{n\mid A_n\cdot B_n = 1= 1_n\} = \{n\mid A_n \not= 0\} \in \cal U.$ Hence, $[A]\hyper \cdot\ [B] = [1].$ \pars
(v) By definition, $[A]+ [B] = [C]$ iff $\{n\mid A_n + B_n = C_n \} \in \cal U.$ However, $\{n\mid A_n + B_n = D_n \} = \nat \in \cal U.$ Hence, $\{n\mid A_n + B_n = C_n\} \cap \{n\mid A_n + B_n = D_n\} = \{n\mid C_n = D_n\} \in \cal U\}.$ Thus, $[C] = [D].$ In like manner, the result holds for ``multiplication.''\qed\parm 
Clearly, one can continue Theorem 2.6 and show that $\langle \hyperreal, \hyper +, \hyper \cdot \rangle$ satisfies all of the ``field'' axioms. It should be obvious, by now, how the ``order'' relation for $\hyperreal$ is defined. \parm
{\bf Definition 2.7} ({ \bf Order}) For each $a =[A],b=[B] \in \hyperreal$ define $a \hyper \leq  b$ iff $\{n\mid A_n \leq B_n \} \in \cal U.$ \parm
I won't show that this relation is well-defined at this time since I'll do it later for all such relations. But, we might as well show that this $\hyper \leq$ is, indeed, a total order and for $\langle \hyperreal, \hyper +,\hyper \cdot, \hyper \leq \rangle$ as a {\it binary relation only} $\hy \leq $ behaves like the $\leq$ behaves for $\real.$ \parm
{\bf Theorem 2.8.} {\it The structure $\langle \hyperreal, \hyper +,\hyper \cdot, \hyper \leq \rangle$ is a totally ordered field.}\pars
Proof. First, notice that $\{n\mid A_n \leq A_n\} = \nat \in \cal U.$ Thus, $\hyper \leq$ is reflexive.\pars
Next, this relation needs to be anti-symmetric. So, assume that $[A] \hyper \leq [B],\ [B] \hyper \leq [A]$. Then $\{n\mid A_n \leq B_n\} \cap \{n\mid B_n \leq A_n \} \subset \{n\mid A_n = B_n\} \in \cal U.$ Hence, $[A] = [B].$\pars
For transitivity, consider $[A] \hyper \leq [B],\ [B] \hyper \leq [C].$ Then 
$\{n\mid A_n \leq B_n\} \cap \{n\mid B_n \leq C_n\} \subset \{n\mid A_n \leq C_n\} \in \cal U.$ Thus, $[A] \hyper \leq [C].$ 
(Notice that the same processes seem to be used each time. They are that $\cal U$ is closed under finite intersection and supersets.)\pars
Next to the notion of ``totally.'' Let $[A], [B] \in \hyperreal.$ Suppose that $[A] \hyper {\not \leq} [B].$ Thus from the trichotomy law for $\real$, $\{n\mid A_n > B_n\} \in \cal U.$ Hence, $[A] \Hyper > [B]$ or $[A] \hyper < [B]$ or $[A] = [B].$ To show that it is a totally ordered ``field'' all that's really needed is to show that it satisfies two properties related to this order and the $\hyper +,\ \hyper  \cdot$ operators. So, let $[A],\ [B],\ [C] \in \hyperreal$. Let $[A] \hyper \leq \ [B].$ Then $\{n\mid A_n \leq B_n \} \subset \{n\mid A_n + C_n \leq B_n + C_n \} \in \cal U.$ Thus $[A] \hyper +\ [C] \hyper\leq [B] \hyper +\ [C].$ Now suppose that $[0] \hyper \leq\ [A],[B].$ Then $\{n\mid 0\leq A_n\} \cup \{n\mid 0 \leq B_n\} \subset \{n\mid 0\leq A_n B_n \} \in \cal U.$ \qed\parm
By the way, using repeatedly the ultrafilter properties to establish the above results is actually unnecessary when a more general result from Mathematical Logic is used. It's this general result that gives me complete confidence that these theorems can be established directly. Indeed, the very definition for the $\hyper {}$ operators comes from this more powerful approach. A major one of these Mathematical Logic results I'll introduce shortly.  \pars
There is often introduced into this subject certain concepts from abstract algebra and abstract model theory. I've decided to avoid this as much as possible for this simplified version. But, now and then, I need to simply state that something holds due to results from these two areas and you need to have confidence that such statements are fact. \pars
What happens next is to ``embed'' the structure $\langle \real, +, \cdot, \leq \rangle $ into $\langle \hyperreal, \hyper + ,\hyper \cdot,\hyper \leq \rangle$ so that the relations $+,\cdot,\leq$ can be considered as but the relations $\hyper +, \hyper \cdot, \hyper \leq$ restricted to $\real.$ 
All one does is to define a function $f$ that takes each $x \in \real$ and gives the unique $[R],$ where $\{n\mid R_n = x\} \in \cal U.$ Notice that one such representation for $[R]$ is the sequence $X_n = x$ for each $n \in \nat.$ Then $\{n\mid X_n=x\} = \nat \in \cal U.$ This is called the {{\bf constant sequence}} representation for $x$ in $\hyperreal.$ This function determines what is called a {{\bf model theoretic isomorphism}} when the relations $\hyper + ,\hyper \cdot,\hyper \leq$ are restricted to the $[X]$ and is what is used to embed  $\langle \real, +, \cdot, \leq \rangle$ into $\langle \hyperreal, \hyper + ,\hyper \cdot,\hyper \leq \rangle.$ One of the big results from abstract model theory states that if one expresses the properties of $\langle \real, +, \cdot, \leq \rangle$ in the customary mathematicians' way (as a first-order predict statement with constants), then every theorem that holds true in $\langle \real, +, \cdot, \leq \rangle$ will hold true when interpreted within this embedding. It's important to note that the real numbers $\real$ are constructed within our basic set theory. Hence, the object $\real$ has a lot of properties. It's assumed that all such properties that can be properly expressed using our present or future defined operations or relations also hold for the structure $\langle \real, +, \cdot, \leq \rangle$.  Thus, simply consider $\langle \real, +, \cdot, \leq \rangle$ as a piece (a substructure) of the structure $\langle \hyperreal, \hyper + ,\hyper \cdot,\hyper \leq \rangle.$ Under this embedding, the notation can be simplified somewhat, by dropping the $\hyper {}$ from the relations $\hyper + ,\hyper \cdot,\hyper \leq$ always keeping in mind that the structure $\langle \real, +, \cdot, \leq \rangle$ is formed by simply restricting these relations to members of the embedded $\real.$ As mentioned each object with which we work and that becomes part of this embedding will be called  a {{\bf standard}} object. All other objects discussed are {{\bf nonstandard}} objects. \pars
At this point, I could go onto some abstract algebra and show without any doubt that the structures $\langle \real, + ,\cdot,\leq \rangle$ and $\langle \hyperreal,  + , \cdot, \leq \rangle,$ although they are both totally order fields, are not the same. But, let's just show that there is a property that $\langle \real, + ,\cdot,\leq \rangle$ has that 
 $\langle \hyperreal, + ,\cdot, \leq \rangle$ does not have.\parm
{\bf Theorem 2.9.} {\it A field property holds for $\langle \real, + ,\cdot,\leq \rangle$ that does not hold for $\langle \hyperreal, + , \cdot,\leq \rangle$.}\pars 
Proof. There is a property of $\langle \real, + ,\cdot,\leq \rangle$ that states that for each $0\leq r \in \real$ there exists an $n \in \nat$ such that $r<n.$ Now the set $\nat$ is a subset of $\real$ and in the embedded form (not yet introduced) $\nat \subset \hyperreal.$ Consider the sequence $A_n = n.$ Then $[A]\in \hyperreal.$ The ultrafilter $\cal U$ is free and does not contain any finite sets. Thus, for each $m \in \nat,$ $\{n\mid A_n \leq m\} \notin \cal U.$ Hence, $\{n \mid A_n > m\} \in \cal U.$ This means that $[A] > [M]=m.$ Since $m$ is arbitrary, then $[A] > [M],$ for each $m \in \nat.$ Hence, at least for the ordinary embedded $\nat$, this field property for   
 $\langle \real, + ,\cdot,\leq \rangle$ does not hold for  $\langle \hyperreal, + , \cdot, \leq \rangle$. \qed
\noindent For those that understand the terminology, the field $\hyperreal$ is also not complete.\parm
From our definition, the $A_n = n$ used to establish Theorem 2.9 would be a nonstandard object. Now let's add a vast number of additional relations $\Phi_i$ to our structures. This will allow us to apply these notions to analysis. The next idea is to ``carve out'' from our set theory some of the important set-theoretic objects used throughout nonstandard analysis.\parm
{\bf Definition 2.10.} ({ \bf Hyper (*) Extensions of standard objects.}) Let $\cal U$ be a free ultrafilter. For any $C \subset \real$ (a 1-ary relation), let 
$b = [B] \in \Hyper C,$ iff $\{n \mid B_n \in C\}\in {\cal U}.$ Let $\Phi$ be any k-ary $(k > 1)$ relation. Then $(a_1,\ldots,a_k) = ([A_1],\cdots, [A_k])\in \Hyper \Phi$ iff $\{ n\mid (A_1(n),\ldots, A_k(n)) \in \Phi \}\in {\cal U}.$ This extension process can be continued for other mathematical entities as required.  \parm
 
Now if it's shown, in general, that these definitions are well-defined, then we can add to our structure additional $n$-ary relations $\Phi_i$ and get the structures $\langle \real, +, \cdot, \leq, \Phi_i\rangle$ and $\langle \hyperreal, +, \cdot, \leq, \Hyper \Phi_i\rangle.$ In which case, as before, we would have that $\Phi \subset \Hyper \Phi$ because all of these relations are actually defined in terms of members taken from $\real.$\parm
{\bf Theorem 2.11.} {\it The hyper-extensions defined in 2.10 are well-defined.}\pars
Proof. In general, for any $[B] \in \hyperreal,$ let $[B] = [B']$. That is let $B' \in \real^\snat$ be any other member of the equivalence class. Let $C \subset \real.$ Then 
$$ \{n\mid B_n = B_n'\} \subset \{n\mid (B_n \in C) \ {\rm if \ and \ only \ if} \ (B_n' \in C)\},$$
$$\{n\mid B_n \in C\} \cap  \{n\mid (B_n \in C) \ {\rm if \ and \ only \ if} \ (B_n' \in C)\}\subset \{n\mid B'_n \in C\},$$  
$$\{n\mid B'_n \in C\} \cap  \{n\mid (B_n \in C) \ {\rm if \ and \ only \ if} \ (B_n' \in C)\}\subset \{n\mid B_n \in C\}.$$ 
The result for this case follows.\pars
For the other k-ary relations, proceed as just done but alter the proof by starting with 
$$\{n\mid B_1(n)= B'_1(n)\} \cap \cdots \cap \{n\mid B_1(n)= B'_1(n)\} \subset$$ 
$$\{n\mid (B_1(n),\ldots B_k(n)) \in \Phi\ {\rm if\ and\ only\ if}\  
(B'_1(n),\ldots B'_k(n)) \in \Phi\}.$$
This completes the proof. \qed\parm
{\bf Definition 2.12.} ({\bf Standard objects operator $\sig {}$.}) I'm using symbols such as $x,y,z,w$ to represent members of $\real$ or for $n >1$ as members of $\real^n = \real \times\cdots \times \real,$ with ``n'' factors. Later, the Roman font for ``variables'' in formal expressions is used. For each $x \in \real$, let $\hyper x =[X]\in \hyperreal,$ where $\{n\mid X_n = x \} =\nat$ (the constant sequence). Then for $X\subset \real$, let $\sig X = \{ \hyper x\mid x \in X\}\subset \hyperreal.$ For $n >1$ and each $x =(x_1,\ldots,x_n) \in \real^n$, let $\hyper x = (\hyper x_1,\ldots, \hyper x_n) \in \Hyper (\real^n).$ 
For $X \subset \real^n,$ $\sig X = \{\hyper x\mid x \in X\} \subset \Hyper (\real^n).$ Each such $\hyper x$ and $\sig X$ is called a {{\bf standard}} object. Thus, $\sig \real$ is the set of embedded real numbers.\parm
What Definition 2.12 does is to identify within $\langle \hyperreal, +, \cdot, \leq, \Hyper \Phi_i\rangle$ the embedded $\langle \real, +,\cdot, \leq, \Phi_i\rangle$ objects. 
For this structure, it's significant that not all useful objects can be hyper-extended by the above, actually necessary, ultrafilter defined extension process. Indeed, because we are only using sequences with range values in $\real$, various members of $\power {\power {\real}}$ cannot be extended. Further, there's a problem if the membership relation $\in$ is extended. Nonstandard analysis exists as a discipline only because the structures $\langle \real, + ,\cdot,\leq, \Phi_i \rangle$ and $\langle \hyperreal, \hyper + ,\hyper \cdot,\hyper \leq, \Hyper \Phi_i \rangle$ can be analyzed externally since they exist as objects in the model of the set theory being used for their construction. In formal set theory as it might appear in Jech (1971), you find that the natural numbers have the property that $0\in 1\in 2\in 3\in 4 \cdots$ and $n \notin n.$  The $\in$ relation is said to be {{\bf well founded}} because there are no types of sequences of members of this set theory that have this processed reversed. There are no objects such that 
$\cdots a \in b \in c \in d.$ If, however, the $\in$ is extended to $\hyper \in$ for members of $\real$, then this $\hyper \in$ is not well founded. \parm
{\bf Example 2.13} Suppose that we do define the $\hyper \in$ for appropriate members of $\real$ using definition 2.10 and using a set theory like Jech (1971).  Thus $\hyper \in$ is defined for each $n \in \nat$ as the $\nat$ is defined within this set theory. Now let's define a collection of sequences from $\nat$ into $\real$ as follows, for each $n \in \nat$, let
$$f_n(i)=\cases{0;&if $i \in \nat$ and $i \leq n$\cr i-n;& if $i >n$\cr}$$ 
Here are what some of these sequences look like. 
$$\cases{f_0(0) = 0,f_0(1) =1, f_0(2) = 2, f_0(3)= 3,f_0(4) = 4,\ldots;&\cr f_1(0) = 0, f_1(1) = 0, f_1(2) = 1, f_1(3) = 2, f_1(4) = 3,\ldots;&\cr f_2(0)=0, f_2(1) = 0,f_2(2) = 0, f_2(3) = 1, f_2(4) = 2, f_2(5) = 3,\ldots.\cr}$$
 Thus, the sequences after the ``0'' values have ``shifting" range values. From the definition of $\hyper \in$, it follows that $\cdots [f_2] \hyper \in [f_1] \hyper \in [f_0].$ To see this, take, say $[f_2], [f_1].$ Then $f_2(0) = 0 \notin f_1(0) = 0,\ f_2(1) = 0 \notin f_1(1) = 0,\ f_2(2) = 0 \in f_1(2) = 1,\ f_2(3) = 1 \in f_1(3), \ldots.$ Hence, $\{n \mid f_2(n) \in f_1(n) \} = \{n\mid n >1\} \in {\cal C} \subset \cal U.$ \parm
Thus, when viewed from the external set theory, $\hyper \in$ is not well founded and does not behave in the same manner as does $\in.$ Further, the $\in$ is used to define the ``hyper'' objects. In order to avoid this problem for the most basic level, the set $\real$ is considered a set of  {{\bf atoms}}
(Jech, 1971) or {{\bf urelements}} or {{\bf individuals}} (Suppes, 1960). This means that each member of $\real$ is not considered as a set and a statement such as $x \in y$ where $x,\  y \in \real,$ has no meaning for our set theory. \par
\vfil\eject

\centerline{\bf 3. HYPER-SET ALGEBRA}
\centerline{\bf INFINITE AND INFINITESIMAL NUMBERS}
\parm
Usually, it's assumed that we are working with one specific free ultrafilter. Is this of any significance for our embedding? \parm 
{\bf Theorem 3.1.} {\it Let infinite $X\subset \nat.$ Then there exists a free ultrafilter $\cal U$ such that $X \in \cal U$. Let $[A],[B] \in \hyperreal.$ Then $[A] =[B]$ for all free ultrafilters iff $\{n\mid A_n = B_n\} \in {\cal C}$. }\pars
Proof. Let infinite $X \subset \nat$. Suppose that $A \in \cal C$ and $A \cap X = \emptyset.$ Then $X \subset \nat -A =$ a finite set. Since $\cal C$ has the finite intersection property, this contradiction implies that ${\cal C} \cup \{X\}$ has the finite intersection property. Hence, there is an ultrafilter ${\cal U}$ such that ${\cal C} \cup \{X\} \subset {\cal U}.$ Obviously, if $\{n \mid A_n = B_n\} \in \cal C$, then $[A] = [B]$ for all free ultrafilters. Suppose that $[A]=[B]$ for $\cal U$ and that $\{n\mid A_n =B_n\}\notin {\cal C}.$ Then $X = \{n\mid A_n \not= B_n\}$ is infinite. Hence there is some free ultrafilter ${\cal U}_1$ and $X \in {\cal U}_1.$ Thus for this ultrafilter $[A] \not= [B]$ and the proof is complete. \qed\parm

Later, for Theorem 3.11, I'll use this result to show that nonstandard objects contained in the same defined set may be considerable different if different free ultrafilter are used.  However, the actual results  obtained when this material is applied to real analysis are, unless otherwise stated, free ultrafilter independent.\pars
\hrule\smallskip\hrule\pars
 The objects that appear in each structure are also objects that can be discussed by means of the set theory of which these objects are members. It's possible to extend these structures to include other objects from this set theory. Shortly, the *-transform process is introduced and the structures  will be slightly extended to use this process in a technically correct manner.   \pars\hrule\smallskip\hrule\parm
{\bf Theorem 3.2.} {\it $\hy {}$-Algebra.\pars 
\indent\indent {\rm (i)} $\Hyper \emptyset = \emptyset.$\pars
\indent\indent {\rm (ii)} If $X \subset \real$ \r [resp. $\real^n\r ],$ then $\sig X \subset \hyperreal$ \r [resp. $\Hyper (\real^n)$\r ].  \pars
\indent\indent {\rm (iii)} If $X \subset \real,$ then $\hyper x \in \sig X$ iff $x \in X$ iff $\hyper x \in \Hyper X$ \pars
\indent\indent {\rm (iv)} Let $X,Y \subset \real$. Then $X \subset Y$ iff $\hyper X \subset \Hy Y.$\pars
\indent\indent {\rm (v)} Let $X,Y \subset \real.$ Then $\Hyper (X-Y) = \hyper X - \Hyper Y.$  \pars
\indent\indent {\rm (vi)} Let $X,Y \subset \real$ Then $\hyper (X \cup Y) = \hyper X \cup \Hyper Y.$ Also, $\Hyper (X\cap Y) = \hyper X \cap \Hyper Y.$\pars
\indent\indent {\rm (vii)} Let $X \subset \real.$ Then $X$ is a nonempty and finite iff $\hyper X = \sig X.$\pars
\indent\indent {\rm (viii)} Let $X_1,\ldots, X_n \subset \real$. For the customarily defined n-ary relations,  $\Phi = (X_1\times\cdots\times X_{n-1}) \times X_n$ iff $\Hyper \Phi = \Hy (X_1\times\cdots\times X_{n-1}) \times \hy X_n= (\hy X_1 \times \cdots \hy X_{n-1}) \times \hy X_n.$ Thus, $\Hyper (\real^n) = (\hyperreal)^n.$ \pars
\indent\indent {\rm (ix)} The statements (iii), (iv), (v), (vi) and (vii) hold for $\real^n,\ n > 1.$\pars
\indent\indent {\rm (x)} For $i >1$ and $\emptyset \not=\Phi \subset \real^n,$ let $P_i$ denote the set-theoretic i'th projection map. Then $\Hyper (P_i(\Phi))= P_i(\Hyper \Phi).$}\pars
Proof. (i) If $S = \emptyset,$ then for any $a \in  \real^\snat$, $\{n\mid A_n \in S\} = \emptyset \notin \cal U.$ Thus our hyper-set algebra yields that $\Hyper \emptyset = \emptyset.$\pars
(ii) This is simply a repeat of Definition 2.12.\pars
(iii) By definition, $\hyper x \in \sig X$ iff $x \in X.$ Now assume that $x \in X.$ Then, by definition, $\hyper x = [X_n], \ X_n = x$ for each $n \in \nat.$ Hence, $\{n\mid X_n \in X\} = \nat \in \cal U.$ Thus, $\hyper x \in \hyper X.$  Conversely, assume that $\hyper x \in \hyper X.$ By definition, $\hyper x = [X_n]$ and $X_n = x$ for all $n\in \nat$. Thus $\emptyset \not= \{n\mid X_n =x\} = \nat \in \cal U$. Hence, $x \in X.$ \pars
(iv) Let $X \subset Y \subset \real$. Then $X \subset \real$ and $a\in \hyper X$ iff $\{ n\mid A_n \in X\} \in \cal U.$ But, $\{ n\mid A_n \in X\} \subset \{n\mid A_n \in Y\}$. Thus, $\{n\mid A_n \in Y\}\in  \cal U$. Now assume that $\hyper X \subset \Hyper Y.$ Then for each $x \in X$, $\hyper x \in \hyper X$ by (iii). Thus $\hyper x \in \Hyper Y.$ Again by (iii) $x \in Y.$ Thus $X\subset Y.$\pars
(v) First, notice  that $X-Y \subset \real$. Let $a \in \Hyper (X-Y)$. Then $\{n\mid A_n \in (X-Y) \}\in {\cal U}\} =U \in \cal U.$ But, this implies that $U \subset \{n\mid A_n \in X\} \in \cal U$ and $ U \subset \{n\mid A_n \notin Y\}.$ Thus, $\{n\mid A_n \notin Y\} \in \cal U.$ Hence, $a \notin \Hyper Y.$ Consequently, $a \in \hyper X - \Hyper Y.$ I'm sure you can establish the converse that $a \in \hyper X - \Hyper Y $ implies that $a \in \Hyper (X-Y).$\pars
(vi) The sets $X\cup Y$ and $X \cap Y$ are subsets of $\real$. Now simply notice that the following identity characterizes the intersection
operator. $C = X \cap Y = X - (X -Y).$ Thus, $\Hyper C= \Hyper (X\cap Y) = \hyper X - (\hyper X - \Hyper Y) = \hyper X \cap \Hyper Y.$ Then $a \in \Hyper (X \cup Y)$ iff $\{n\mid A_n \in (X \cup Y)\}= \{n\mid A_n \in X\} \cup \{n\mid A_n\in Y\}.$ Hence, if $\{n\mid A_n \in (X \cup Y)\}\in \cal U$, then either $\{n\mid A_n \in X\}\in \cal U$ or $\{n\mid A_n \in Y\} \in \cal U$. Thus, $\Hyper (X \cup Y) \subset \hyper X \cup \Hyper Y.$ Since $X \subset (X \cup Y)$ and $Y \subset (X\cup Y)$, it follows from (iv) that $\hyper X \cup \Hyper Y \subset \Hyper (X \cup Y)$ and the result follows. \pars
(vii) The first part is established by induction. Let $X = \{x\}.$ Then $\{x\} \subset \real.$ By definition $a \in \hyper X$ iff $\{n\mid A_n \in \{x\}\} = \{ n\mid A_n = x \} \in \cal U.$ Now $\hyper x = [X]$ and $\{n\mid X_n = x\} = \nat \in \cal U.$ Thus, $\{n\mid X_n = B_n\} = \{n\mid X_n = x\} \cap \{n\mid B_n = x\} \in \cal U$ implies that $[X] = [B] = \hyper x.$  Assume the result holds for 
a set with $k$ members. Then $\Hyper {\{x_i, \ldots, x_{k+1}\}} = \Hyper {(\{x_1,\ldots, x_{k}\} \cup \{x_{k+1}\})} = \Hyper {\{x_1,\ldots, x_k\}} \cup \Hyper {\{x_{k+1}\}} = \{\hyper x_1,\ldots,\hyper  x_{k+1}\}$ by the induction hypothesis and (v), the result holds for any $k \geq 1.$\pars
For the converse, let infinite $X \subset\real$ and assume that $\sig X = \hyper X.$ There exists an injection $B\colon \nat \to X$. Hence $\{B_n\mid n\in \nat\}$ is an infinite subset of $X$. Let $\hyper x = [X] \in \sig X.$ Then $X_n = x\in X$ for each $n \in \nat$. But, $\{n\mid X_n = B_n\}$ is finite. Hence, $[X] \not= [B]$ since $\{n\mid X_n \not= B_n\} \in {\cal C}$. Also $\{n\mid B_n \in X\} = \nat \in \cal U$ implies that $b \in \hyper X.$ 
There is no $x \in X$ such that $\hyper x = b\in \hyper X$ implies $\sig X \not= \hyper X.$ 
\pars
(viii) The customarily defined notion of an n-ary relation can be found in Jech (1971). The first idea is that a 1-ary relation is but the  subset of $\real$ and this has been established in (iii). The other cases, $n >1,$ for $x_i \in X_i \subset \real,\ 1< i \leq n$ the Cartesian product $X^n$ is characterized by the statement that $(x_1,\ldots,x_n) \in  (X_1 \times \cdots \times X_{n-1}) \times X_n$ iff $x_i \in X_i, \ 1\leq i \leq n$, where the actual ``Cartesian product'' is defined by induction. That is 
$X_1 \times X_2 \times X_3 = (X_1 \times X_2) \times X_3$ etc. (There are other ways to define the Cartesian product more formally just using 2-tuples and finite sequences.) Note that for any $k >1,\ \{n\mid (A_1(n),\ldots,A_k(n)) \in (X_1 \times \cdots \times X_{k-1})\times X_k\}= \{n\mid A_1(n) \in X_1\}\cap \cdots \cap \{n\mid A_k(n)\in X_k\}.$ The result follows from basic filter properties.\pars 
(ix) Statements (iii), (iv), (v) are proved in the exact same manner for  $(\hyperreal)^n.$ Statements (vi), (vii) are proved by application of the method used in (vii) coupled with the characterization used to establish (Viii).\pars 
(x) Let $a \in \Hyper {(P_i(\Phi))}.$ Then $\{n\mid A_n \in (P_i(\Phi))\} \in \cal U$ iff $\{n\mid {\rm \ there\ exist \ } B_1,\ldots,B_{i-1},B_{i+1},\ldots,$\hfil\break $ B_m\ {\rm \ such\ that\ }  (B_1(n),\ldots,B_{i-1}(n),A(n),$ $B_{i+1}(n),\ldots, B_m(n))\in \Phi \} \in \cal U.$ Hence, there exist $b_j, 1\leq j \leq m,\ i \not= j$ such that $(b_i,\ldots,b_{i-1},a,b_{i+1},\ldots, b_m) \in \Hyper \Phi.$ But from the definition of such a projection this gives that $a \in P_i(\Hyper \Phi)).$ Thus $\Hyper {(P_i(\Phi))} = P_i(\Hyper \Phi)$ because of the equivalence of the two set-theoretic statements. (At present, I have not introduced a more formal way of writing definitions for such sets.) \qed  \parm
{\bf Important.} Theorem 3.2 involves properties about our original structure $\langle \real, +, \cdot, \leq, \Phi_i\rangle$ and the $\langle \hyperreal, +, \cdot, \leq, \Hyper \Phi_i\rangle$ {\bf and} the embedded objects. Although the embedded objects ``behave'' like the original objects, they're still different from these. This observation will come into play when I discuss the notion of *-transform.\parm 

It's about time that I demonstrated that the ``ideal'' numbers used by Leibniz, that did not really exist mathematically until 1961, exist within $\hyperreal.$ These are the infinitesimals which solve this three hundred year old problem.\parm
{\bf Definition 3.3.} ({\bf Infinite and infinitesimal numbers.}) As usual define the absolute value function (i.e. binary relation) for members of $a \in \hyperreal$ by requiring as $\Hyper \vert a \vert = \vert a \vert = b$ iff $\{n\mid \vert A_n \vert = B_n\} \in \cal U.$ Although, I won't show it, but later will show how to establish the fact, this function $\Hyper \vert \cdot \vert$ has the same mathematical properties as does $\vert \cdot \vert$ for members of $\real$. So, I have written it as if its a restriction to our embedded $\sig \real$ of the usual absolute value function. An $a \in \hyperreal$ is {{\bf infinitely large}} or simply {{\bf an infinite number}} or shorter still {{\bf infinite}} iff $\hyper x < \vert a \vert$ for each $\hyper x \in \sig \real.$ (Some might go back to the pre-embedded $\real$ for these definitions, but remember we are thinking of $\sig \real$ as our actual set of real numbers.) A $b \in \hyperreal,$ is an {{\bf infinitesimal}} or as Newton stated {{\bf infinitely small}} iff $0\leq \vert b \vert < \hyper x$ for each $x \in \real^+$, the set of all {\bf positive real numbers}.\parm
Now that these ``new'' types of numbers are defined, do any exist?\parm
{\bf Example 3.4.} Let $A$ be the member of $\real^\snat$ with the property that $A_k = k$ for each $k \in \nat.$ Let $\hyper x \in \sig \real.$ Then there exists some $m \in \nat$ such that $\vert x \vert < m.$ Hence, $\vert x \vert = \vert X_n \vert < A_m =m$ for each $n \in \nat$ implies that $\{n \mid A_n > \vert X_n \vert  
\} \supset \{m,m+1,\ldots \} \in {\cal C} \subset \cal U.$ Thus, $a$ is an infinite number. There are a lot more. \pars
Note that $\Hyper 0$ is the trivial infinitesimal. Consider, the sequence $G_n = 1/n, \ n\in \nat -\{0\}$ and $G_0 = 0.$ Then $g \not= \Hyper 0.$ Now for each $x \in \real^+$ there is some $m \in \nat,\ m \not=0$ such that $0 < 1/m < x.$ Thus $\Hyper 0 < \hyper 1/\Hyper m < \hyper x.$ (Note: So far we would need to establish such statements by ultrafilter properties. But, this is really trivial because by definition each $\hyper x$ is the constant sequence representation.) Now $\nat - \{n \mid  G_n  \geq X_n\}$ is a finite subset of $\nat.$ Hence, $\{n \mid 0 < G_n < X_n\} \in {\cal C}\subset \cal U.$  Thus, $g$ is an infinitesimal. Indeed, once we get one nonzero infinitesimal, we can generate infinitely many. \parm
{\bf Definition 3.5.} A $a\in \hyperreal$ is {{\bf finite}} or {{\bf limited}} iff it's not infinite. That is if there is some $\hyper x \in\sig\real^+$ (positive embedded reals) such that $\vert a \vert \leq \hyper x.$ The {{\bf set of all finite numbers}} is denote by $G(0)$, the {{\bf galaxy}} within our universe $\hyperreal$ in which $\sig \real$ resides. (Note: If $a \in \hyperreal$, then $a \in G(0)$ iff there is some $\Hy y$ such that $\vert a\vert < \Hy y$.) The {{\bf set of all infinitesimals}} is denoted by $\mu (0).$\parm
What is the algebra of the infinitesimals and does this algebra display the exact algebra used by Newton and Leibniz? It's customary to let {{\bf lower Greek letters represent nonzero infinitesimals.}} Then as one would expect {{\bf capital Greek letters represent infinite numbers.}} I'll show later that many real valued functions defined on open intervals about zero preserves infinitesimals. Indeed, if $x > 0$ and $f\colon (-x,x) \to \real$ is continuous at $x = 0$ and $f(0) = 0,$ then $\hyper f(\eps) = \lambda$ or $0.$  (Notation: ``$f$ is a function that takes each and every member of $(-x,x)$ and yields members of $\real.$'')\pars
We know that $\hyperreal$ is a totally ordered field and $\mu (0), G(0) \subset \hyperreal.$ Also, $0 \in \mu(0) \cap G(0).$ The next Theorem, 3.8, gives exactly how the relations $+, \cdot, \leq$ behave when they are restricted to $\mu (0)$ and $G(0).$ It will turn out that both of these sets are {{\bf totally order rings with no zero divisors}}. What does this mean  for the relations $\hy +, \hy \cdot, \hy \leq$? This means that $\mu (0)$ and $G(0)$ behave for these binary relations exactly like the integers $\{\cdots, -3,-2-1,0,1,2,3,\cdots \}$ behave with the one exception that $\Hyper 1 \notin \mu (0).$ By the way, if your interested, the $x,y$ are {{\bf zero divisors}} iff $xy = 0$ implies that $x=0$ or $y = 0.$
  \parm
Certainly, establishing results by using the basic properties of the  ultrafilters is getting a bit tedious. There most be a better way. And, there is. If a statement about $\langle \real, +,\cdot, \leq, \Phi_i \rangle$ is expressed in a special way and that statement holds, then there is a process that's used to show that a altered statement holds in $\langle \hyperreal, +, \cdot, \leq, \Hyper \Phi_i \rangle.$ 
The process is called {{\bf *-transform}}. However, to do this properly it's necessary to extend the structure considerably.\parm
\hrule\smallskip \hrule\smallskip
(It's not really necessarily that you fully understand the contents of my new extended structure. You could just go immediately to Definition 3.6 and simply restate Theorem 3.7 only in terms of the notation $\cal M$ and $\hy {\cal M}$ and without the structures being specified.) Because of the way I have defined the *-extension operator in Definition 2.10, the structure I technically need is $\langle \real,\ldots, \real^n,\ldots,\power {\real},\ldots, \power {\real^n},\ldots,+, \cdot, \leq \rangle.$ Although not actually necessary I have identified the three indicated binary relations used previously. The actual $n$ used in practice is rather small, usually. Let each element of each of the objects in $\{\real,\ldots, \real^n,\ldots,\power {\real},\ldots, \power {\real^n},\ldots \}$ have a ``constant'' name and we use $+,\cdot,\leq$ and the like as the ``names'' for these specific objects.  Also the constant that ``names'' a mathematical object itself will not be differentiated from the mathematical object itself. Let $Cn$ be this set of all of these  constants. Notice that ``constants'' that are members of a n-ary relation $n >1$ like $(x,y,z),$ use the constants $x,y,z$ from $\real.$ These n-tuple forms $(x_1,\ldots, x_n)$ are part of our language.\smallskip \hrule\smallskip\hrule\parm
{\bf Definition 3.6.} (*-transform) Consider any properly formed statement (formally a first-order formula with equality and constants using the atomic formula in the appendix) with bounded quantifiers and only using members of $Cn$. Then the *-transform of this statement is obtained by writing a $\hyper {}$ to the left as a superscript of each constant. Also, there is the reverse process where a statement in terms of the $Cn$ is obtained by removing the $\hyper {}.$\parm
I'm not going to present a course in first-order logic in this monograph. So, you'll simply need to assume that I've expressed the ``formal'' statements in the proper {\bf bounded form.} This means that the ``variable'' that appears to the right of a quantifier, the universal ``for each,'' $\forall \r x$, and the existential ``there exists some,'' $\exists \r x,$ must vary over one of the sets in the standard structure. Of course, it turns out that mathematicians seem to always write their informal sentences in forms that are logically equivalent to these bounded forms. The reason for the bounded form is that in the appendix Theorem A3 establishes the following without using the Axiom of Choice. Only some previous results obtained using ultrafilters are needed. Notice in what follows two new symbols are introduced for the respective structures and the structures are now extended slightly.\parm
{\bf Theorem 3.7.} {\it Let $\cal S$ be any sentence in bounded form that uses only constants in $Cn.$ Then $\cal S$ holds for ${\cal M} = \langle \real,\ldots, \real^n,\ldots,\power {\real},\ldots, \power {\real^n},\ldots,+, \cdot, \leq \rangle$ iff the *-transform of $\cal S$ holds in $\hyper {\cal M}= \langle \hyperreal,\ldots, \hyperreal^n,\ldots,\power {\hyperreal},\ldots, \power {\hyperreal^n},\ldots,+, \cdot, \leq \rangle.$}\parm 
\smallskip\hrule\smallskip\hrule\smallskip
\noindent{\bf WARNING} If someone who has experience with nonstandard analysis reads Theorem 3.7, they might state that the theorem is in error, since I have written $\power {\hyperreal^n},$ etc., in the structure. However, it is correct as shown in the appendix for the language being used. For example, an expression such as $\exists x(x \in \power {\hyperreal})$ is not in the proper *-transfer form. This theorem only applies to $\exists x(x \in \hyper {\cal P}(\hyperreal))$. \pars\hrule\smallskip\hrule\smallskip
I don't suppose that you noticed that *-transform is a reversible process that relates our original structure $\cal M$ and the nonstandard structure $\hyper {\cal M}$ and does not technically mention the embedded objects. Here's one place where there is a different notational approach. In some work, the original structure and the embedded structure are consider as identical. I will probably not do this if there's any possible confusion. 
Also there are actually certain properties of $\real$ that can't be expressed in our formal language. The language could be enriched.  But, if this is done one might as well go all the way to the object called a {{\bf superstructure}}. The idea is to see what can be accomplished without such an enriched language and the additional complications this would produce.\parm 

{\bf Theorem 3.8.} {\it The sets $\mu (0),G(0)$ are totally ordered subrings of $\hyperreal$ with no zero divisors and $G(0)$ has an identity.}\pars
Proof. I start with $G(0)$ since it contains $\mu (0)$ and $G(0) \subset \hyperreal$ to show that it is a totally ordered ring with no zero divisors, all that is needed is to show that it is closed under the operations $+,\ \cdot.$ To do this efficiently Theorem 3.4 is used. Informally, we know that if we are given any two real numbers $x,y,$ then $\vert x + y \vert \leq \vert x\vert + \vert y \vert.$ The formal bounded statement of this fact is 
$$\forall \r x\forall \r y ((\r x \in \real)\land (\r y\in \real) \to \vert \r x + \r y \vert \leq \vert \r x \vert + \vert \r y \vert)$$
holds in  $\cal M$, and, hence, its *-transform holds in $\hyper {\cal M}.$ Thus,$$\forall \r x\forall \r y ((\r x \in \hyperreal)\land (\r y\in \hyperreal) \to \vert \r x + \r y \vert \leq \vert \r x \vert + \vert \r y \vert)$$ 
is a fact about $\hyper {\cal M}.$ (Note: You could have written $\vert \cdot \vert$ as $\Hyper \vert \cdot \vert.$ I also point out that this is actually considered as written in the form for a particular $\Phi_j$, where we have in our language the ordered n-tuple notation. Define $\Phi_j = \{(w,y,z)\mid \vert w \vert \leq \vert y \vert + \vert z \vert\}.$ Then $\vert \r x + \r y \vert \leq \vert \r x \vert + \vert \r y \vert$ is equivalent to $((\r w, \r x, \r y) \in \Phi_j) \land (\r w = \r x + \r y)$.)\pars
Thus, the triangle inequality holds in $\hyperreal. $ So, let $a,b \in G(0).$
Then there are standard $\hyper x, \Hy y \in \sig \real$ such that $\vert a \vert < \hyper x,$ and $\vert b \vert < \Hy y.$ But $\vert a + b\vert \leq \vert a \vert + \vert b \vert < \hyper x + \Hy y= \Hyper (x + y)$ from our definitions and the order properties of 
 of $\hyperreal.$ This gives that $a + b \in G(0)$ which gives us closure under $+$ since $\Hyper 0 \in G(0)$. In like manner, one gets that $ab \in G(0).$ Of course, since $G(0) \subset \hyperreal$ the members have the usual associative, commutative and distributive properties and $\Hyper 0$ is its zero. Now either by *-transform or filter properties $\Hyper 1$ is also an identity in $G(0).$ It now follows immediately that since $\hyperreal$ is a totally ordered field, then $G(0)$ is a totally ordered ring. Further, since $\hyperreal$ has no zero divisors neither does $G(0).$\pars
Next, consider $\mu (0).$ We can apply the method used to establish that $G(0)$ is a totally ordered ring with no zero divisors (but $1 \notin \mu (0)$) to show that $\mu (0)$ is a totally ordered ring with no zero divisors. The only difference in the proofs is that instead of writing that there is some $\hyper x > 0$ such that $\vert a \vert < \hyper x$, we have that $ \eps \in \mu (0)$ iff $\vert \eps \vert < \hyper x$ for all arbitrary $\hyper x >0.$ \qed \parm 
Does $\mu (0)$ have any other \underbar{significant} algebraic properties? The answer is yes and it's this most remarkable property that's needed if its members are to mimic the ``infinitely small'' notion of Newton. What $\mu (0)$ does is to ``absorb'' via multiplication every member of $G(0).$ It has this ``ideal'' property. A $I \subset G(0)$, is an ideal iff it is a subring (which $\mu (0)$ is) and for each $a \in G(0)$ and each $b \in I$, the product $ab \in I.$ An ideal $I \subset G(0)$ is maximum iff for any other ideal $I_1 \supset I,$ 
$I = I_1$ or $I = G(0).$\parm
{\bf Theorem 3.9}. {\it The set of infinitesimals $\mu (0)$ is a proper maximum ideal in $G(0)$.}\pars
Proof. Let $a \in G(0)$ and $\eps \in \mu (0).$ Then there is some $\hyper x \in \sig \real$ such that $\vert a \vert < \hyper x.$ Let $\eps = g = [G].$ Consider arbitrary positive $\Hy y.$ Then $$F=\{n \mid \vert A_n \vert < X_n = x\} \cap \{n\mid \vert G_n \vert < Y_n = y\} \in \cal U.$$ However, $\{n\mid \vert A_nG_n\vert < xy \} \supset F.$ But $xy$ is also arbitrary. Hence, $a\eps \in \mu (0).$ \pars
Let $I$ be any ideal in $G(0)$ such that $\mu(0) \subset I.$ Assume that there is some $b \in I - \mu (0).$ Then $b \not = 0$ and there exists some positive $x\in \real$ such that $\{n \mid \vert B_n \vert \geq r\} \in \cal U.$ Hence, $\{n\mid \vert 1/B_n \vert \leq 1/x \}\in \cal U.$ Consequently, $[B^{-1}] = b^{-1} \in G(0)$ implies that $\Hyper 1 = bb^{-1} \in I.$ This last fact will always force $I = G(0)$ since it's an ideal. Well, take any $0\not= x \in \real.$ Then $\hyper x \not= \Hyper 0$ and $\hyper x \notin \mu (0)$  implies that $\mu (0) \not= G(0).$ Hence, $\mu (0)$ is a proper maximum ideal. \qed\parm
I could go into some other abstract algebra material and use the language of quotient rings, isomorphisms and kernels to show exactly how $G(0)$ and $\mu (0)$ are related to $\sig \real$, but it's unnecessary to do this for this simplified approach. It's enough to say that the properties of $\mu (0)$ exactly match the ``infinitely small'' of Newton and the ``ideal numbers'' of Leibniz. \pars 
It has become customary to drop the $\hyper{}$ from the members of $\sig \real$ when there is no confusion. I'll start doing this in the very important next definition.\parm
{\bf Definition 3.10} (Monads of standard numbers.) Let $x \in \sig \real$. Then the {{\bf monad of (about) $x$}} is the set $\mu (x) = \{x + \eps \mid \eps \in \mu (0).$ The only standard object in $\mu (x)$ is $x$. (Recall that when there's no confusion, I might use $x$ in place of $\hy x.$)\parm
Before showing a remarkable relation between the monads and $G(0)$, I need the next theorem. \parm
{\bf Theorem 3.11.} {\it Let $A_n$ be a sequence of real numbers. Then $[A] \in \mu (x)$ for every free ultrafilter iff $\lim_{n \to \infty}A_n = x.$}\pars
Proof. First, note that, for a fixed free ultrafilter $\cal U$ and its monad $\mu (x)$, $[A] \in \mu (x)$ iff there is some $\eps \in \mu (0)$ such that $[A] = x +\eps,\ [A] - x = \eps$ iff $[A] - x \in \mu (0)$ iff $\{n \mid \vert A_n -X_n \vert < r\}\in \cal U$ for any arbitrary positive $r.$  Let $\cal U$ be any free ultrafilter and assume that $A_n \to x.$ Then for arbitrary positive $r$, we have that $\vert A_n - x \vert < r$ for all but a finite number of $A_n$. Thus, $\{n \mid \vert A_n -X_n \vert <r\} \in {\cal C} \subset \cal U.$ But $r$ is arbitrary implies that $[A] \in \mu (x).$ \pars
Conversely, assume that $A_n \not\to x.$ Then there is a positive $r$ such that $X = \{n \mid \vert A_n -x \vert \geq r\}$ is an infinite set. Any infinite subset of $\nat$ is contained in some free ultrafilter ${\cal U}_1$ by Theorem 3.1. Thus, for this ${\cal U}_1$, $[A] \notin \mu_1 (x)$ since the complement of $X$ is not a member of ${\cal U}_1.$ \parm
{\bf Theorem 3.12} {\it The collection $\{\mu (x)\mid x \in \sig \real\}$ is a partition for $G(0).$}\pars
Proof. Technically, to be a partition of $G(0),$ one must have that $\mu(x) \cap \mu (y) \not= \emptyset$ implies that $\mu (x) = \mu (y)$ and that $\bigcup \{\mu (x)\mid x \in \sig \real\} = G(0).$ For the first part, assume that there exists some $a \in \mu (x) \cap \mu (y).$ Then $a = \eps + x, \ a = \lambda + y.$ But, $\eps + x = \lambda + y$ implies that $\eps - \lambda = y -x.$ This is only possible if $\eps - \lambda = 0$ since $y - x \in \sig \real.$ Thus $x = y$.
Let $a \in \bigcup \{\mu (x) \mid x \in \sig \real\}.$ Then $a = \eps + x$ for some $x \in \sig \real.$ Then $\vert a \vert = \vert \eps +x \vert \leq \vert \eps \vert + \vert x \vert < \vert x \vert +1.$ Hence $a \in G(0)$.
Consequently, $\bigcup \{ \mu (x) \mid x \in \sig \real \} \subset G(0).$\pars

Now assume that $a \in G(0).$ Rather than continue to use the properties of the our free ultrafilter, let's just consider the properties of $<$. Hence, there is some $\hyper x \in \sig \real^+$ such that $a < \hyper x.$ So, consider the set $S = \{y \mid \Hy y < a\}$ This set is nonempty since $-x \in S.$ Also since $a < \hyper x$, $S$ is set of real numbers that's bounded above and as such has a least upper bound $z$. The number $z$ needs to be located. Assume that $\vert z - a \vert$ is not an infinitesimal. Thus there is some $w \in \real$ such that $\vert \Hyper z - a \vert > \Hyper w.$ Suppose that $\Hyper z < a$. Then $a - \Hyper z >  \Hyper w$ implies that $\Hyper z + \Hyper w = \Hyper (z + w)< a $ implies $z + w \in S$ and $z$ is not the least upper bound. So, let $a < \hyper z$. This implies that $a < \Hyper (z -w) < \Hyper z$. But, $z - w$ is an upper bound for the set $S$. This contradicts the least upper bound property for $z.$ Hence, $\Hyper z - a =\eps$ implies that $a \in \mu(z)$. \qed\parm
Of course, this implies, in general, that $\bigcup \{\mu (x)\mid x \in \sig \real\} = G(0)$ is free ultrafilter independent. But, the monads that contain some of the members of $\real^\snat$ cannot be readily determined. \parm
{\bf Example 3.13} Consider the sequence $a = \{1,-1,1,-1,\ldots \}.$ Then as done in the proof of Theorem 3.11, $U_1 = \{n\mid A_n = 1\} \cap U_2= \{n \mid A_n = -1\} =\emptyset.$ The set $U_1$ is a member of the free ultrafilter ${\cal U}_1$ and $U_2$ is a member of the free ultrafilter ${\cal U}_2$ where ${\cal U}_1 \not= {\cal U}_2.$ Further, $a \in \mu_1(1)$ and $a \in \mu_2(-1).$\parm
There are some other useful properties that relate members of $\hyperreal$, and the sets $G(0)$ and $\mu (0)$ and show that they model the older notions of real number {``infinities''} and the {``infinitely small.''} Let $\hyperreal - G(0) = \real_\infty$ be the {{\bf infinite numbers.}} It's immediate from the definition that if $a,b \in \real_\infty,$ then $ab \in \real_\infty.$ If $0 < a $ [resp $ a < 0$] $\in \real_\infty$ and $a < b$ [resp. $b < a$], then $b \in \real_\infty$ since $a > r$ [resp. $a< r$] for each $r \in \real.$ \parm
{\bf Theorem 3.14} {\it \pars
\indent\indent {\rm (i)} If $b \in \real_\infty$, then $1/b \in \mu(0).$\pars
\indent\indent {\rm (ii)} If $0 \not= \eps \in \mu(0)$, then $1/\eps \in \real_\infty.$\pars
\indent\indent {\rm (iii)} Let $\eps \in \mu (0),\ b \in \mu (x).$ Then $\eps + b \in \mu (x)$ and $\eps\, b \in \mu (0).$\pars
\indent\indent {\rm (iv)} If $b \in \real_\infty$, and $\hyper x \not= \Hyper 0,$ then $b\hyper x \in \real_\infty.$ If $a \in G(0) - \mu (0),$ then $ba \in \real_\infty.$ \r (The $\hyperreal_\infty$ almost has the special property associated with an ideal.\r )\pars
\indent\indent {\rm (v)} If $\hyper x < \Hy y$, then $\hyper x + \eps < \Hy y + \lambda$ for any $\eps, \lambda \in \mu(0).$} \pars
Proof. (Most mathematicians would consider these proofs as trivial and would ``leave them to the reader.'' But, I'll do most of them.)\pars
(i) If $b \in \real_\infty,$ then for any $\hyper x\in \sig \real+$, $\hyper x < \vert b \vert$. Thus by field properties, $1/\vert b\vert < \hyper x.$ This says that $1/\vert b \vert \in \mu (0).$ \pars
(ii) Same method as (i).\pars
(iii) Let $\eps \in \mu (0)$ and $b \in \mu (x).$ Then $b = \hyper x + \lambda$ implies that $\eps + b  = \hyper x + \eps + \lambda = \hyper x + \gamma \in \mu (0).$ 
Then $\eps \, b \in \mu(0)$ from Theorem 3.9 or $\eps\, b= \eps\hy x + \eps\gamma = \alpha + \beta \in \mu (0).$ \pars
(iv) Using (i), $0\not= 1/(\hyper x b) \in \mu (0)$. Now use (ii). For the second part, use the fact that if $b \in G(0) - \mu(0),$ then if $a > \Hyper 0$, there is some $\hyper x > 0$ such that $\hyper x < a$ and if $a < 0,$ then there is some $\Hy y$ such that $a < \Hy y .$ Now apply the remark I made just prior to this theorem. \pars
(v) Assume that $0\leq \eps -\lambda \in \mu (0).$ Hence, for $a < b$, $0 \leq \eps -\lambda
< b -a$. Thus $a +\eps < b +\lambda.$ \qed\parm 
The fact that the $\{\mu (\hy x)\mid \hyper x \in \sig \real \}$ forms a partition of $G(0)$ immediately defines for all members of $G(0)$ an equivalence relation of some importance, where this relation is a short hand for a member of $G(0)$ being in a unique $\mu (x).$ \parm
{\bf Definition 3.15.} ({{\bf Infinitely close (near) equivalence relation}}.) Two $a,b \in G(0)$ are {\bf infinitely close} iff $a -b \in \mu (0).$ This relation is written as $a \approx b.$ \parm 

We almost have enough of the basic machinery to continue with real analysis. But, there is one last major procedure that needs to be introduced, the ``standard part'' operator.\parm
{\bf Definition 3.16.} ({\bf The standard part operator, st.}) Using Theorem  3.12, there is a function $\St$ on $G(0)$ into $\sig \real$ such that, for each $\mu (x)$, $\st {\mu (x)} = \hyper x \iff x.$ Once the properties of $\St$ are obtained, then, usually, one further allows $\st {\mu (x)} = x \in \real.$ The function $\St$ is called the {{\bf standard part operator}}.\parm
 Most of the results in the next theorem would be what one would expect.\parm
 
 {\bf Theorem 3.16.} {\it Let $\St {\colon} G(0) \to \sig \real$ \r ($\real$\r) be the standard part operator. Then for each $a,\ b \in G(0),$\pars 
\indent\indent {\rm (i)} $\st {a \pm b} = \st a \pm \st b.$ \pars
\indent\indent {\rm (ii)} $ \st {ab} = \st a \st b.$\pars
\indent\indent {\rm (iii)} If $a \leq b$, then $\st {a} \leq \st {b}.$ \pars
\indent\indent {\rm (iv)} $\st {\vert a\vert} = \vert \st a \vert,\ \st {\max \{a,b\}}= \max \{\st a, \st b\},$ $\st {\min \{a, b \} } = \min \{\st a, \st b\}.$\pars
\indent\indent {\rm (v)} $\st a = 0$ iff $a \in \mu(0).$\pars
\indent\indent {\rm (vi)} For any $\hyper x$, $\st {\hyper x} = \hyper x.$\pars
\indent\indent {\rm (vii)} The $\st a \geq 0$ iff $\vert a\vert \in \mu ({\st a}).$\pars 
\indent\indent {\rm (viii)} $a \approx b$ iff $a - b \in \mu (0)$ iff $\st a = \st b$ \pars 
\indent\indent {\rm (ix)} If $\st a \leq \st b,$ then either $a - b\in \mu (0)$ or $a \leq b.$\pars
\indent\indent {\rm (x)} If $\Hyper 0 < c$ \r [resp. $c < \Hyper 0$\r ] and $c \in \hyperreal_\infty,$ then for $a \geq \Hyper 0,$ $\Hyper 0 < c + a \in \hyperreal_\infty$ \r [resp.  $a \leq \Hyper 0,\ c + a \in \hyperreal_\infty$ \r ].
}\pars
Proof. I'll do (iii) and leave the others to the reader. Let $a,b \in G(0), \  a \leq b.$ Then $a \in \mu (\st {a}), b \in \mu (\st {b})$ implies that $a = \st {a} + \eps,\ b = \st {b} + \gamma$ implies that $0\leq \st {b} -\st {a} + \gamma - \eps= \st {b - a} +\gamma - \eps$ implies that $\st {a} \leq \st {b}$ since the monads are disjoint.\qed\parm 
Is it clear that if $\nat_\infty = \hypernat -\sig \nat$, then $\nat_\infty \subset \real_\infty$? \parm
{\bf Theorem 3.17.} {\it The set of infinite natural numbers $\nat_\infty \subset \real_\infty$ and, for each $n \in \nat$ $\Hy n < \Lambda$ for each $\Lambda  \in \nat_\infty$. }\pars
Proof. In example 3.4, the infinite number defined is actually a member of $\nat^\snat \subset \hyperreal.$ Thus, $\hypernat - \sig \nat \not= \emptyset.$ For each $m \in \sig \nat$, there is the $m +1 \in \sig \nat \subset \sig \real.$ Hence, $\sig \nat \subset G(0).$ Let $a\in \hypernat - \sig \nat$ and $a \in G(0).$ Then since $\Hy0 \leq b$ for each $b \in \hypernat,$ there is some $r \in \sig \real$ such that $\Hy 0 \leq \vert a \vert = a < r.$ But, we know there is some $m \in \nat$, hence, $\hyper m \in \sig \nat$ such that $\Hy 0 \leq a < \hyper m.$ *-transform of the statement 
``for each $x$, for each $y,$ for each $z,$ if $x \in \nat$ and $y \in \nat$ and $z \in \nat$ and $x \leq y$, then $z \in [x,y]$ iff $0\leq z\leq y$'' or formally  $\forall x\forall y\forall z((x \in \nat)\land (y \in \nat)\land (z \in \nat)\land (x \leq y) \to (z \in [x,y] \iff (0\leq z \leq y)))$ holds in $\hyper {\cal M}$ and characterizes the set $\Hyper [\Hy 0, \hyper m].$ Thus, $a \in \Hyper [\Hy 0, \hyper m].$ But, since $[0,m]$ is a finite set, then Theorem 3.2 (vi) implies that $a = \hyper n$ for some $\hyper n \in \sig \nat.$ This contradiction implies that $\nat_\infty \subset \real_\infty$ as one would expect. The second part follows immediately. \qed  
 \vfil\eject

\centerline{\bf 4. BASIC SEQUENTIAL CONVERGENCE}
\parm
One intuitive statement about sequences of real numbers states something like ``all the convergence properties are determined by the behavior of the {infinite tails}.'' In fact, for elementary converges, we have ``Well, the values of the sequence get {nearer, and nearer, and nearer and stay near} to the limit no matter how far you go out in the series.'' Does the nonstandard theory of sequential convergence model both ``getting nearer, and nearer'' and ``staying near'' simultaneously? Indeed, you'll find out that, for convergence, the infinite tails are {\bf all} members of $G(0)$. Moreover, each nonstandard characteristic based directly upon a definition is stated in, at least, one less quantifier. {G\"odel} considered that just {removing one quantifier} from any characterization is a major achievement within mathematics. By the way, all of the results presented in the remainder of this book are free ultrafilter independent. Also, many of the definitions and proofs presented are easily generalized to the multi-variable calculus. (I'll use the notation $n \in \hypernat$ where there is no confusion as to the location of the $n$. Usually one might write this as $a \in \hypernat.$) \parm
{\bf Theorem 4.1.} {\it A sequence $S\colon \nat \to \real$ is bounded iff $\hyper S(n)\in G(0)$ for each $n\in \hypernat$ iff \r ($\hyper S[\hypernat] \subset G(0).$\r )}\pars
Proof. Let $S$ be bounded. Recall what this means. There exists some $x \in \real^+$ such that ``for each $n \in \nat$, $\vert S(n) \vert < x$'' or $\forall \r y((\r y \in \nat)\to (\vert S(\r y) \vert < x))$ holds in $\cal M.$ By *-transform, $\forall \r y((\r y \in \hypernat)\to (\hyper \vert \hyper S(\r y) \hyper \vert < \hyper x))$ holds in $\hyper {\cal M}.$ (Note: We can consider with respect to our embedding that $\vert \cdot \vert$ is but a restriction of $\Hyper \vert \cdot \Hyper \vert$ and we need not use the $\hyper {}$ there, although this is but a notational simplification.) Hence, for each $n \in \hypernat,$ $\hyper S(n) \in G(0).$\pars
Conversely, for each $n \in \hypernat$, let $\hyper S(n)\in G(0).$ We know that there is a $b \in \hyperreal^+_\infty \subset \hyperreal^+$ such that for each $c \in G(0),$ $\vert c \vert < b.$ Hence, $\exists \r x ((\r x\in \hyperreal^+ )\land \forall \r y ((\r y \in \hypernat) \to \vert \hyper S(\r y) \vert < \r x))$ holds in $\hyper {\cal M}$. Thus, the statement $\exists \r x ( (\r x\in \real^+ )\land \forall \r y ((\r y \in \nat) \to \vert  S(\r y) \vert < \r x))$, obtained by dropping the $\hyper{}$, holds in $\cal M$ and the sequence is bounded. \qed \parm
What about the ``near to $L$'' and ``stays near'' intuitive notion and it's relation to the ``true'' infinite part of the tail? \parm
{\bf Theorem 4.2} {\it A sequence $S \colon \nat \to \real$ converges to $L \in \real$ \r ($S_n \to L$\r ) iff $\hyper S(\Lambda) - L \in \mu (0)$ for each $\Lambda \in \nat_\infty$   iff $\hyper S(\Lambda)\in \mu (L)$ for each $\Lambda \in \nat_\infty$ iff $\st {\hyper S(\Lambda)} = L$ for each $\Lambda \in \nat_\infty$ iff \r ($\hyper S[\nat_\infty] \subset \mu (L).$\r )}\pars
Proof. Let $S \colon \nat \to \real$ converge to $L$. Let $y \in \real^+.$ Then we know that there exists some $m \in \nat$ such that for each $k \in \nat$ where $k \geq m$, $\vert S(k) - L  \vert <x.$ Hence, the statement
$$\forall \r x ((\r x \in \nat)\land ( \r x > m) \to (\vert S(\r x) - L \vert < y))$$
holds in {\cal M}; and, hence, in $\hyper {\cal M}.$ In particular, by *-transform, for each $\Lambda \in \nat_\infty$, $\vert \hyper S(\Lambda)-L \vert < \Hyper y.$ Since, $y$ is arbitrary, then $\hyper S(\Lambda) - L \in \mu (0)$ for each $\Lambda \in \nat_\infty.$ Hence, $\hyper S(\Lambda) \in \mu (L)$ and $\st {\hyper S(\Lambda)} = L$ for each $\Lambda \in \nat_\infty.$\pars
Conversely, assume that $(\hyper S(\Lambda) -L)\in \mu (0)$ for each $\Lambda \in \nat_\infty.$ Let $y \in \real$. Since $\nat_\infty \not= \emptyset$, then by Theorem 3.17, the sentence
$$\exists \r z((\r z \in \hypernat) \land \forall \r x((\r x \in \hypernat)\land (\r z < \r x) \to ( \vert \hyper S(\r x) - L \vert < \Hyper y)))$$
holds in $\hyper {\cal M}.$ Thus, it holds in $\cal M$, by reverse *-transform. But, this is the standard statement that $S_n \to L.$ All the remaining ``iff'' are but restatements of $\hyper S(\Lambda )- L\in \mu (0)$ for each $\Lambda \in \nat_\infty.$ \qed\pars  
{\bf Corollary 4.3} {\it All the basic limit theorems for sums, products, etc. all follow from Theorem 4.2 and the properties of the ``$\St$'' operator.}\parm 
{\bf Examples 4.4}\pars
(i) $(1/n)^p \to 0,\ n,p > 0,\ p \in \nat.$ We know that for each nonzero $\Lambda \in \nat_\infty,\ (1/\Lambda)^p \in \mu(0)$ and the result follows. (If we had result that the continuous function $f(x) = x^p,\ p >0,$ preserves infinitesimals, then we could extend this to any $p > 0.$ But, maybe it's better to use sequences to motivate continuity.)\pars
(ii) $x^n \to 0,\ 0 < \vert x \vert <1.$ In general, for any $n, m \in \nat$ such that $n < m$, we have that $(1/\vert x\vert)^n < (1/\vert x \vert)^m$. For any $y \in \real$, there is some $n \in \nat$ such that $\vert y\vert < (1/\vert x\vert)^n$. Hence, for each $\Lambda \in \nat_\infty$ $(1/\vert x \vert)^\Lambda \in \real_\infty.$ Consequently, $x^\Lambda \in \mu (0)$ for each $\Lambda \in \nat_\infty.$ \pars
(iii) Let $0 < x,\ x \not= 1.$ Then $x^{1/n} \to 1,\ n > 0.$ Consider that case that $x > 1$ and $S_n = x^{1/n} -1.$ Then  $x = (S_n +1)^n.$ Hence, $x > nS_n$ for each $n > 0.$ Thus, by *-transform, $x > \Lambda \hyper S_\Lambda$ for each $\Lambda \in  \nat_\infty.$ Consequently, $0 < \hyper S(\Lambda) < (x/\Lambda) \in \mu (0)$ for each $\Lambda \in \nat_\infty.$ Thus $\hyper S(\Lambda) \in \mu (0),\ \Lambda \in \nat_\infty$ and result follows in this case.\pars
Now, if $0 < x < 1$, then $1<1/x$ and, as just shown, $(1/x)^{1/\Lambda} \in \mu (1),\ \Lambda \in \nat_\infty.$ Thus $(1/x)^{1/\Lambda} -1 = \eps \in \mu (0).$ Hence, $1 - x^{1/\Lambda} = \eps(x^{1/\Lambda}) \in \mu (0),$ since by *-transform $0 < x^{1/\Lambda} < 1.$ Therefore, $x^{1/\Lambda}\in \mu (1)$ for this case also and the complete result follows.\pars
(iv) $n^{1/n} \to 1,\ n > 0.$ Consider again the sequence $S_n = n^{1/n} -1.$ Then $n = (1 +S_n)^n = \sum_{k = 1}^n \left(\matrix{n\cr k\cr}\right)S^k_n \geq \left(\matrix{n\cr 2\cr}\right)S^2_n,\ n > 1.$ Thus, 
$0 \leq S_n \leq ({{2}\over{n-1}})^{1/2},\ n >1.$ By *-transform, and in particular, 
 $0 \leq \hyper S_\Lambda \leq ({{2}\over{\Lambda-1}})^{1/2},\ \Lambda \in \nat_\infty.$ But, $({{2}\over{\Lambda-1}})^{1/2}\in \mu (0).$ Hence $\hyper S(\Lambda) \in \mu (0)$ and the result follows from the definition of $S_n.$\parm
It seems that some of the above {algebraic manipulations} are what one might do if these limits were established without using nonstandard procedures. There are major differences, however, in the number of quantified statements one needs for the standard proofs as compared to the nonstandard. Let's establish a standard result by nonstandard means.\parm
{\bf Theorem 4.5.} {\it Every convergent sequence of real numbers is bounded.}\parm
Proof. Let $S_n\to L \in \real.$ Then $\hyper S(\Lambda) \in \mu (L) \subset G(0)$ for each $\Lambda \in \nat_\infty.$  Since $\hyper S[\sig \nat] \subset G(0),$ the result follows from Theorem 4.1. \qed\parm
{\bf Theorem 4.6.} {\it A set of real numbers $B$ is bounded iff $\Hyper B \subset G(0).$}\pars
Proof. If $B$ is finite, then it's immediate that $\hyper B \subset G(0).$ If $B$ is infinite, then there is some real number $x$ such that for each $y \in B$, $\vert y \vert \leq x.$ By *-transform of the obvious expression any $a \in \Hyper B$ has the property that $\vert a \vert \leq \hyper x.$ Consequently, $\Hyper B \subset G(0).$ \pars
Conversely, if $B$ is not bounded, then for $n \in \nat$ there is some $x \in B$ such that $\vert x \vert > n.$ Hence, by *-transform, there is some $p\in \Hy B$ such that $\vert p \vert \geq  \Lambda, \ \Lambda \in \nat_\infty.$ From the remark made prior to Theorem 3.14, $p \notin G(0)$ and the converse follows. \qed\parm
One of the first big results one encounters in sequential convergence theory 
is a sufficient condition for convergence. Recall that a sequence is {monotone} iff it is either an increasing or decreasing function. The following characterization is what would be expected, that for monotone sequences only one infinite number is needed for convergence.\parm\vfil\eject 
{\bf Theorem 4.7.} {\it If $S \colon\nat \to\real$ is monotone and there exists some $\Lambda \in \nat_\infty$ such that $\hyper S(\Lambda) \in G(0),$ then $S_n \to \st {\hyper S(\Lambda)}.$} \pars
Proof. Simply assume that $S\colon \nat \to \real$ is increasing since the decreasing case is similar. I first note that $\Hyper y=\st {\hyper S(\Lambda)} \in \sig \real.$ By *-transform, the extension $\hyper S \colon \hypernat \to \hyperreal$ is increasing. Thus for $\Lambda \in \nat_\infty$ and for each $\hyper m \in \sig \nat,\  \hyper S(\hyper m) \leq \hyper S(\Lambda)$ and, since $\hy S(\Lambda) \in G(0),$ the $\st {\hyper S(\hyper m)} = \Hyper (S(m))\in \real.$ Consequently, for  $\real$, then following sentence
$$ \forall \r x((\r x\in \nat) \to (S(\r x) \leq y))$$
holds in ${\cal M}$; and, hence, holds in  $\hyper {\cal M}.$ So, let $\Omega \in \nat_\infty.$ Then $\hyper S(\Omega) \leq \Hyper y =\st {\hyper S (\Lambda)};$ which implies that for each $\Omega \in \nat_\infty,\ \hyper S(\Omega) \in G(0)$.  Thus $\st {\hyper S(\Omega)} \in \sig \real$ for each such $\Omega$ and $\st {\hyper S(\Omega)} \leq \Hyper y.$ Let $\Omega > \Lambda$. Then $\hyper S(\Lambda)\leq \hyper S(\Omega);$ which implies that $\Hyper y \leq \st {\hyper S(\Omega)} .$  
 But, since the above statement still holds for such a $\Omega,$ then 
$\hyper S(\Omega)\leq \Hyper y$ implies that $\st {\hyper S(\Omega)} \leq \Hyper y$. Let $\Omega < \Lambda.$ Then $\hy S(\Omega) \leq \hyper S(\Lambda)$; implies $\st {\hyper S(\Omega)} = z$ and the above statement holds for $z.$. Thus, $\st {\hyper S(\Lambda)} \leq \st {\hyper S(\Omega)}.$ Hence, $\st {\hyper S(\Omega)} = \st {\hyper S(\Lambda)}$ for each $\Omega \in \nat_\infty.$ Consequently, 
$\hyper S(\Omega) - \hyper S(\Lambda) \in \mu (0)$ for all $\Omega \in \nat_\infty$ implies that $\hyper S(\Omega) \in \mu (\st {\hyper S(\Lambda)})$ for each $\Omega \in \nat_\infty$ and the result follows. \qed\pars
{\bf Corollary 4.8.} {\it A bounded monotone sequence converges.}\pars
Proof. By Theorem 4.1.\parm
Please note that if $a - b \in \mu (0),$ and $b \in G(0),$ then the intuitive statement that $a \in \mu (\st {b})$ does, indeed, hold. In a slightly more general mode recall that for a sequence $S\colon \nat \to \real$ a real number $w$ is an {{\bf accumulation point}} or {{\bf limit}} point for $S$ iff for each $r \in \real^+$ and for each $n \in \nat$, there is some $m \in \nat$ that $m > n$ and 
$\vert S_m - w \vert < r.$ This definition allows $1$ to be an accumulation point of sequences such as $\{1,1/2,1,1/3,1,1/4,1,\ldots \}$ where both $1$ and $0$ are accumulation points. This definition does not correspond to most of the accumulation point definitions for  point-sets. However, there will be a another term used in chapter 8, that does so correspond. \parm
{\bf Theorem 4.9.} {\it \pars
\indent\indent {\rm (i)} A $w \in \real$ is an accumulation point for $S \colon \nat \to \real$ iff there exists some $\Lambda \in \nat_\infty$ such that $\hyper S(\Lambda) \in \mu (\Hyper w)= \mu (\st {\hyper S(\Lambda)}.$\pars
\indent\indent {\rm (ii)} A sequence $S \colon \nat \to\real$ has an accumulation point iff there exists some $\Lambda \in \nat_\infty $ such that $\hyper S(\Lambda) \in G(0).$}\pars
Proof (i) Let $w \in \real$ be an accumulation point for $S$. Then the sentence
$$\forall \r x\forall \r y((\r x\in \real^+)\land (\r y \in \nat)\to \exists \r z((\r z	\in \nat) \land(\r z > \r y)\land (\vert S(\r z) - w \vert < \r x)))$$ holds in $\hyper {\cal M}$ by *-transform. So, let $0 < \eps \in \mu (0)$ and $\Omega \in \nat_\infty.$ Then there exists some $\Lambda \in \hypernat$ such that $ \Lambda > \Omega$ and $\vert \hyper S(\Lambda)-\Hyper w\vert < \eps.$ Hence, $\hyper S(\Lambda) \in \mu (\Hyper w).$ Clearly, $\Lambda \in \nat_\infty.$ \pars
Conversely, assume that there exists some $\Lambda \in\nat_\infty$ such that $\hyper S(\Lambda) \in \mu (\Hyper w),\ w \in \real.$ Note that $\mu (\Hyper w) \subset \Hyper (w - y, w + y)$ for each $ y \in \real^+$ and that $\Lambda > \hyper n$ for all $\hyper n \in \sig \nat$. Hence, for given $\Hyper w, \Hyper y > \Hyper 0$ and a given $\hyper m$, we have that
$$\exists \r x((\r x \in \nat)\land(\r x > \hyper m) \land (\vert \Hyper S(\r x) - \Hyper w \vert < \hyper y))$$
holds in $\hyper {\cal M}$; and, hence, in $\cal M$ by reverse *-transform. The result follows. \pars 
(ii) This follows ``immediately'' from (i). \qed \par\medskip 
{\bf Theorem 4.10.} {\it {\rm (i)} A sequence $S \colon \nat \to \real$ has a subsequence that converges to $w \in \real$ iff there exists some $\Lambda \in \nat_\infty$ such that $\hyper S(\Lambda) \in \mu (w).$ \pars
\indent\indent {\rm (ii)} A sequence has a convergent subsequence iff there is some $\Lambda \in \nat_\infty$ such that $\hyper S(\Lambda) \in G(0)$ iff \r ($\hyper S[\nat_\infty] \cap G(0) \not= \emptyset$ \r ).}\pars
Proof. (i) Assume that for $\Lambda \in \nat_\infty$ that $\hyper S(\Lambda) \in \mu (\Hyper w).$  Then $w$ is an accumulation point. You start with $n = 0$ and take $y = 1$. Then you have an $S_m$ such that $\vert S_m - w \vert < 1$. Let $S'_0 = S_m$. Now take $y = 1/2$ and consider the next $S_k$ as the one for $k > m$ and $\vert S_k - w \vert < 1/2.$ This idea can be restated in an induction proof for the other $1/n$ with no great difficulty. This subsequence obviously converges to $w$. (Have I used the Axiom of Choice to obtain the $S_k$?) \pars
On the other hand, if $S' \colon \nat \to \real$ is a subsequence of the sequence $S$ and it converges to $L$, then $\hyper S'[\nat_\infty] \cap G(0) \not= \emptyset$ implies, since $\hyper S'[\nat_\infty] \subset \hyper S[\nat_\infty],$ that $L$ is an accumulation point by Theorem 4.9. (ii) is obvious. \qed\par\medskip 
{\bf Theorem 4.11.} {\it A bounded sequence has a convergent subsequence.}\pars
Proof. From Theorems 4.1 and 4.10. \qed\par\medskip
Have I convinced you that the notion of what happens with the truly infinite tail piece of a sequence does determine all that seems necessary for basic convergence? No. Well, let's look at another idea, the {special types of divergence} written as $S_n \to +\infty$ [resp. $-\infty$]. \pars
Recall that a sequence $S_n \to +\infty$ [resp. $-\infty$] iff for each $y > 0$ [resp. $y < 0$] there exists an $m \in \nat$ such that for each $n \in \nat$ such that $n \geq m,$ $S_n \geq y$ [resp. $S_n \leq y$]. How do we intuitively state such stuff as this? One might say that $S$ {{\bf converges to ``plus infinity'' or converges to ``negative infinity.''}} But, in  basic real analysis, the ``numbers'' $\pm \infty$ do not actually exist. \par\medskip
{\bf Theorem 4.12.} {\it For sequence $S\colon \nat \to \real,$ $S_n \to +\infty$ \r [resp. $-\infty$ \r ] iff for each $\Lambda \in \nat_\infty$ $\hyper S(\Lambda) \in \real^+_\infty$ \r [resp. $\real^-_\infty$ \r ], where $\real^+_\infty = \{\Lambda \mid \Hyper 0 < \Lambda \in \real_\infty\}$ \r [resp. $\real^-_\infty = \{\Lambda \mid \Hyper 0 > \Lambda \in \real_\infty\}$\r ], iff \r ($\hyper S[\real_\infty] \subset \real^+_\infty$ \r [resp. $\real^-_\infty$\r ]\r ).} \pars
Proof. Assume that $S_n \to +\infty.$ We can assume that $S_n >0$ for each $n \in \nat$ since it is not true for only finitely many $n$. Suppose that there exists some $\Lambda \in \nat_\infty$ such that $\hyper S(\Lambda) \notin \real^+_\infty.$ Thus, $\hyper S(\Lambda) \in G(0).$ Therefore there is a subsequence of $S$, $S' \colon \nat \to \real$ and $\hyper S'_n \to L,\ L \in \real$ and $L = \st {\hyper S'(\Lambda)}.$ Thus, there exists an $m \in \nat$ such that for all $n\geq m$, $\vert S'(n) - L \vert <1.$ Hence, for each such $n$, $0 < S'(n) < L + 1.$ Thus, considering $y = L +1$ there does not exist a $p \in \nat$ such that for each $n\in \nat,$ where $n\geq p,$
$S_n \geq y.$  \pars
Conversely, suppose that for each $\Lambda \in \nat_\infty$, $\hyper S(\Lambda) \in \real^+_\infty.$ Let $y > 0$. Consider $\Omega \in \nat_\infty.$ If $\Lambda \geq \Omega,$ then $\Lambda \in \nat_\infty$ and under the hypothesis, $\hyper S(\Lambda) > \Hyper y.$ Consequently, the sentence
$$\exists \r x((\r x \in \hypernat) \land \forall \r z ((\r z \in \hypernat)\land (\r z \geq x)\to (\hyper S(\r z) > \Hy y))$$
holds in $\hyper {\cal M}$ and, hence, in $\cal M.$ This result for the positive infinite numbers follows by reverse *-transform. The case for the  negative infinite numbers follows in like manner and the proof is complete. \qed\par\medskip
By the way, notice how easily the next result is established. \par\medskip
{\bf Theorem 4.13.} {\it If $S:\nat \to \real$ converges to $L \in \real$, then $L$ is unique.}\pars
Proof. If $L\not= M \in \real$, then $\mu (L) \cap \mu (M) = \emptyset.$ \qed\parm 
 Let's recap some the significant nonstandard characterizations for a sequence $S\colon \nat \to \real.$ Notice that they are all quantifier ($\forall,\ \exists$) free. \parm
{\bf Theorem 4.14.} {\it Given a sequence $S \colon \nat \to \real.$ Then \pars
\indent\indent {\rm (i)} $S$ is bounded iff $\hyper S[\hypernat ] \subset G(0),$\pars
\indent\indent {\rm (ii)} $S_n \to L$ iff $\hyper S[\nat_\infty] \subset \mu (L)\subset G(0),$\pars
\indent\indent {\rm (iii)} $S$ has a convergent subsequence iff $\hyper S[\nat_\infty] \cap G(0) \not= \emptyset.$\pars
\indent\indent {\rm  (iv)} $S_n \to \pm \infty$ iff $\hyper S[\nat_\infty] \subset \real^{\pm}_\infty$.}\parm
I wonder whether $S \colon \nat \to \real$ has a subsequence that converges to $\pm \infty$ iff $\hyper S[\nat_\infty] \cap \real^{\pm}_\infty \not= \emptyset$? So far, to show that a specific sequence converges we needed to guess at what the limit might be. One of the more important notions was considered by Cauchy, the Cauchy Criterion, that for the real numbers characterizes convergence without having to guess at a limit $L.$ A sequence $S$ is called a {{\bf Cauchy sequence}} iff for each $y \in \real^+$, there is some $m \in \nat$ such that for each pair $p,q \in \nat$ such that $p,q \geq m$, it follows that $\vert S(p) - S(q) \vert < y.$ \par\medskip
{\bf Theorem 4.15.} {\it \r (Nonstandard Cauchy Criterion.\r ) A sequence $S \colon \nat \to \real$ is Cauchy iff
$$\hyper S(\Lambda) - \hyper S(\Omega) \in \mu (0)$$ for each $\Lambda,\Omega \in \nat_\infty.$}\pars
Proof. For the necessity, simply let real $y > 0$, then there exists some $m_y \in \nat$ such that the sentence
$$\forall \r x \forall \r z((\r x\in \nat)\land (\r z \in \nat)\land(\r z 
\geq m_y)\to (\vert S(\r x ) - S(\r z)\vert < y))$$
holds in $\cal M$ and, hence, in $\hyper {\cal M}.$ In particular, if $\Lambda, \Omega \in \nat_\infty,$ then $\Lambda, \Omega > \hyper m_y$  for any such $m_y$ implies that $\vert\hyper S(\Lambda) -\hyper S(\Omega)\vert < \Hyper y$ for any  $y > 0$. Consequently, $\hyper S(\Lambda) -\hyper S(\Omega)\in \mu (0).$ 
\pars
The sufficiency follows in the usual manner since $\mu (0) \subset \Hyper (-y,y)$ for each $y > 0$ and $\nat_\infty \not= \emptyset$ imply that the sentence
$$\exists \r w((\r w\in \nat) \land \forall \r z \forall \r x ((\r z\in \nat)\land(\r x \in \nat)\land (\r x\geq \r w)\land (\r y \geq \r w) \to (\vert S(\r x) - S(\r z) \vert < y))$$ holds in ${\cal M}$ and the proof is complete. \qed \par\medskip
{\bf Theorem 4.16.} {\it A sequence $S\colon \nat \to \real$ converges iff it is Cauchy.}\pars
Proof. Suppose that $S_n \to L\in \real.$ Then for each pair $\Lambda,\ \Omega \in \nat_\infty,$ $\hyper S(\Lambda) - L \in \mu (0),$ and $ \hyper S(\Omega) - L \in \mu (0).$ Hence, $\hyper S(\Lambda) - \hyper S(\Omega) \in \mu(0).$ \pars
For the converse, let $S\colon \nat \to \real$ be Cauchy. Then for $\Lambda,\ \Omega \in \nat_\infty$, we have that $\hyper S(\Lambda) - \hyper S(\Omega) \in \mu (0)$ from Theorem 4.15. Let $\Lambda \in \nat_\infty$ and $\hyper S(\Lambda) \in G(0).$ Then $\hyper S[\nat_\infty] \subset \mu (\st {\hyper S(\Lambda)}= \mu (\hyper L)$ implies that $S_n \to L.$ So, assume the other possibility, that $\hyper S(\Omega) \notin G(0)$ for any $\Omega \in \nat_\infty.$ This implies that $S$ is unbounded. Let $m  \in \nat$ and let $y = \max\{\vert S_m \pm 1\vert, \vert S_0\vert,\ldots, \vert S_m\vert\}$. Then there is some $p \in \nat$ such that $\vert S_m \pm 1 \vert < y < \vert S_p\vert$ and  
 $p >m.$ Thus by *-transform, given any $\Lambda \in \nat_\infty,$ there is some $\Omega \in \nat_\infty$ such that $\vert \hyper S(\Lambda) \pm \Hyper 1 \vert < S(\Omega).$ Notice that (i) $\hyper S(\Lambda) \pm \Hyper 1 \in \nat^+_\infty,$ in which case, since $\hyper S(\Lambda) - \hyper S(\Omega)\in \mu (0)$, it follows that $\hyper S(\Omega) \in \nat^+_\infty$ or (ii) $\hyper S(\Lambda) \pm \Hyper 1 \in \nat^-_\infty,$ in which case $\hyper S(\Omega) \in \nat^-_\infty.$ For case (i), consider $\Hyper S(\Lambda) +1 < \hyper S(\Omega)$; for case (ii), consider $\hyper S(\Omega) <\Hyper S(\Lambda) -1.$ For these two cases, this yields that $ 1 < \vert \hyper S(\Omega)- \hyper S(\Lambda)\vert  \notin \mu (0).$  This contradicts the hypothesis that $\hyper S(\Lambda) - \hyper S(\Omega) \in \mu (0).$ The proof is now complete. \qed \parm\vfil\eject

\centerline{\bf 5. ADVANCED SEQUENTIAL CONVERGENCE}
\parm
Recall that a {{\bf double sequence}} $S\colon \nat \times \nat \to \real$ 
converges to $L \in \real$ iff for each $y \in \real^+$, there is some $p \in \nat$ such that for each pair $n,m \in \nat,$ such that $n,m \geq p$ and  $ \vert S(n,m) - L \vert < y.$ The same nonstandard characteristics hold for such convergence as in the single sequence case.\parm
{\bf Theorem 5.1.} {\it A sequence $S \colon \nat \times \nat \to \real$ converges to $L \in \real$ iff $\hyper S(\Lambda, \Omega) - L \in \mu (0)$ for each $\Lambda,\Omega \in \nat_\infty$ iff $\hyper S(\Lambda, \Omega)\in \mu (L)$ for each $\Lambda,\Omega \in \nat_\infty$ iff $\st {\hyper S(\Lambda, \Omega)} = L$ for each $\Lambda, \Omega \in \nat_\infty$ iff \r ($\hyper S[\nat_\infty\times \nat_\infty] \subset \mu (L).$\r )}\pars
Proof. With but almost trivial alterations, this proof is the same as the one for Theorem 4.2. \qed \parm
{\bf Example 5.2} Let $S(m,n) = {{m}\over{1+mn^2}}.$ Then for each $\Lambda, \Omega \in \nat_\infty,\ {{1 + \Lambda \Omega^2}\over{\Lambda}} = {{1}\over{\Lambda}} + \Omega^2 \notin G(0).$ Hence, ${{\Lambda}\over{1 + \Lambda \Omega^2}} \in \mu (0)$ and, thus, $S(n,m) \to 0.$ \parm
The following results, and many more, for double sequences follow in the same manner as in Chapter 4.\parm
{\bf Theorem 5.3.} {\it Every convergent double sequence is bounded.}\parm
{\bf Theorem 5.4.} {\it \r (Nonstandard Cauchy Criterion.\r ) The sequence $S\colon \nat \times \nat \to \real$ converges to $L \in \real$ iff for each $\Lambda, \Omega, \Lambda', \Omega' \in \nat_\infty,\ \hyper S(\Lambda, \Omega) - \hyper S(\Lambda',\Omega') \in \mu (0).$}\parm
In the theory of double sequences, one of the interesting questions, at the least to most mathematicians, is the role played by the {{\bf iterated sequences}}, (in brief limit notation) $\lim_n(\lim_m s(n,m))$ and 
$\lim_n(\lim_m s(n,m))$. What this notation means is that, taking the first iterated limit, you might have that for each $n,$ $\lim_mS(n,m) = S'(n) \in \real$. Then, maybe, $\lim_n S'(n) \in \real.$ Now for a convergent double sequence, is it always the case that the iterated sequence converges? In the example 5.2, notice that for $n =0,$ $S(n,m)$ diverges. Indeed, take any natural number $a$. Then the sequence $S(n,m)= {{m}\over{1 + m(n-a)^2}}$ will have this same problem for $n = a.$ \parm
{\bf Example 5.5.} Consider the sequence $S(n,m) = {{m+1}\over{m+n +1}}.$ Then for any $n\in \nat,$ $S(n,m) \to 1,$ while for a fixed $m,$ $S(n,m) \to 0.$ This shows that the double sequence does not converge since the $n,m \in \nat$ are arbitrary pairs and as such it should not matter if one is held fixed and the other varies, the limit being unique, as in this single sequence case, must be the same in all cases. As is well know, this behavior for double sequences is simply a reflection of the same problems that occur with multi-variable real valued functions. \parm
The problem displayed by examples like 5.2, does not occur for members of $\nat_\infty$ as indicated by the following rather interesting pure nonstandard result.\parm
{\bf Theorem 5.6.} {\it Let $S(m,n)$ converge to $L.$ Then for any sequence $\Omega_m \in \nat_\infty$, $\lim_m \st {\hyper S(\hyper m,\Omega_m)} =L$ \r [resp. $\lim_n \st {\hyper S(\Omega_n,\hyper n)} =L.$}\pars
Proof. Let $S(m,n)\to L.$ We know that for $y\in \real^+$ there is some $p\in \nat$ such that for each pair $m,n \in \nat$ and $m\geq p$ and $n \geq p,$ $\vert S(m,n) - L \vert < y$. Now $p$ may be assumed fixed for the $y.$ By *-transform, it follows that $\vert \hyper S(\hyper m,b) - L \vert < \Hyper y$ for each $\hyper m \geq \hyper p$, and $b \in \hypernat,\ b\geq \hyper p.$  Hence, in particular for any sequence $\Omega_m \in \nat_\infty, \ \vert \hyper S(\hyper m,\Omega_m) - L \vert < \Hyper y.$ Now taking the standard part operator on each side of this   
inequality implies that $\vert \st{\hyper S(\hyper m,\Omega_m)} -L \vert \leq y$ for each $\hy m\geq \hyper p.$ This statement is sufficient to state that $\st {\hyper S(\hyper m, \Omega_m)} \to L.$ (Note: Technically, the $m$ that appears in the sequence notation $\Omega_m$ should be considered as a member of $\sig \nat.$ But, this does not come from the *-transform of any standard sequence or any allowed formal statement using our simple language.) \qed\parm
The use of a sequence such as $\Omega_m$ changes the double sequence  $S(m,n)$ into a nonstandard type of ordinary sequence. One of the major concerns for double sequences is their relation to the iterated limits, where I'll use abbreviated limit notation. There are, of course, standard theorems that relate convergence of double sequences to the convergence of iterated limits.\parm 
{\bf Theorem 5.7.} {\it Let $S(m,n) \to L \in \real.$ Then 
$\lim_m(\lim_n S(m,n)) = L$ iff $\lim_n S(m,n)$ exists for each $m \in \nat.$}\parm
Proof. The necessity is obvious. So, assume that $\lim_n S(m,n) = r_m \in \real$ for each $m \in \nat.$ Then, by *-transform for each $\Omega \in \nat_\infty,\ \st {\hyper S(\hyper m, \Omega)} = r_m$ for each $m \in \nat.$ All we need to do is to consider any sequence $\Omega_m \in \nat_\infty,$ like the constant sequence $\Omega_m = \Omega$, and obtain $\lim_m(\st {\hyper S(\hyper m,\Omega_m)} = \lim_m r_m = L$ by Theorem 5.6. \qed \parm
A theorem such as Theorem 5.7 holds with an interchange of the $n$ and $m$ symbols. Theorem 5.7 gives a condition under which an iterated limit will converge to the limit of a converging double sequence. But, are there  necessary and sufficient conditions that determine completely when the limit of a double sequence corresponds to the limit of both iterated limits?  
Of course, if there is, it probably is not obvious. We need something special to happen. Consider the limit statement for $\lim_mS(m,n),$ where $\lim_m S(m,n) \in \real.$ Then  $\lim_mS(m,n)$ {{\bf converges uniformly in n}} iff for each $y \in \real^+$ there exists some $p \in \nat$ such that for each $n \in \nat,$ and  $m,m' \in \nat,$ where $m,\ m' \geq p$ it follows that 
$\vert S(m,n) - S(m', n) \vert < y.$ Thus, the $p$ is such that the sequence $S(m,n)$ seems to behave like an ordinary convergent sequence independent from the actual value of $n \in \nat.$ Let's see if this notion has a somewhat simply nonstandard characteristic. Indeed, one that parallels the statement for a sequence being Cauchy.\parm
{\bf Theorem 5.8.} {\it Let $S \colon \nat \times \nat\to \real.$ Then $\lim_m S(m,n)$ converges uniformly in $n$ iff  
$$\hyper S(\Lambda, n) - \hyper S(\Omega, n) \in \mu (0)$$
for each $\Lambda,\Omega \in \nat_\infty$ and for each $n \in \hypernat$.}\pars
Proof. For the necessity, simply consider the *-transform. Use  the fact that  from the definition $\Hyper y$ is arbitrary, and then select particular $\Lambda, \Omega \in \nat_\infty.$ \pars
For the sufficiency, assume that for each $n \in \hypernat$, $\hyper S(\Lambda, n) - \hyper S(\Omega, n) \in \mu (0)$ for each pair $\Lambda, \Omega \in \nat_\infty.$ Let $y \in \real^+.$ Notice that, for a particular $\Lambda, \Omega,$ there's a $\Gamma\in \nat_\infty$ such that $\Lambda, \Omega \geq \Gamma$  and $\vert \hyper S(\Lambda, n) - \hyper S(\Omega,n) \vert < \hyper y.$ Thus, the sentence
$$\exists \r x((\r x\in \nat)\land \forall \r y \forall \r z \forall \r w ((\r y \geq \r x)\land(\r z \geq \r x) \land (\r w\in \nat)\to (\vert S(\r y,\r w) - \hyper S(\r z, \r w)\vert < \Hyper y)))$$ holds in $\hyper {\cal M}$ and, hence, in $\cal M$ by reverse *-transform. This completes the proof. \qed \parm
{\bf Corollary 5.9.} {\it If for each $n \in \nat$, $\lim_mS(m,n)= S_n \in \real$, then $\lim_m S(m,n)$ converges uniformly in $n$ iff for each $n \in \hypernat$, $\hyper S(\Lambda, n) - \hyper S(n) \in \mu (0)$ for each $\Lambda \in \nat_\infty.$}\pars
Proof. This follows immediately from Theorem 4.15, the Cauchy Criterion for convergence. \qed \parm
There is a standard necessary and sufficient condition for the equality of the limit of the double sequence and its iterated limits . \parm\vfil\eject
{\bf Theorem 5.10.} {\it Let $S \colon \nat \times \nat \to \real.$ Then $\lim S(m,n) = \lim_m(\lim_n S(m,n)) = \lim_n(\lim_m S(m,n)) \in \real$ iff\pars
\indent\indent {\rm (i)} $\lim_m S(m,n)$ converges uniformly in $n$ and\pars
\indent\indent {\rm (ii)} $\lim_n S(m,n)$ converges for each $n \in \nat$.}\pars
Proof. For the necessity, it's clear that (ii) follows from the convergence of the iterated limit. Then $\lim_n (\lim_m S(m,n)) \in \real$ implies that $\lim_m S(m,n) \in \real$ for each $n \in \nat.$ Hence, by the Theorem 4.15, $\hyper S(\Lambda, n) - \hyper S(\Omega,n) \in \mu (0)$ for each $n \in \sig \nat$ and each $\Lambda,\Omega \in \nat_\infty.$ By Theorem 5.4, we also have that $ \hyper S(\Lambda, \Gamma) - \hyper S(\Omega, \Gamma) \in \mu (0)$ for each $\Gamma \in \nat_\infty.$ Hence, $\hyper S(\Lambda, n) - \hyper S(\Omega, n) \in \mu (0)$ for each $n \in \hypernat.$  
Therefore, Theorem 5.8 yields that $\lim_m S(m,n)$ converges uniformly in $n$. \pars
For the sufficiency, let $\lim_m S(m,n) -S_n=0$ for each $n \in \nat$. 
 Uniformly in $n$ means that this limit is independent from the $n$ used.  By *-transform, we have that for any $\Omega \in \nat_\infty$, $\lim_m (\st {\hyper S(\hyper m, \Omega)} = \st {\hyper S(\Omega)}.$ Thus, for each $\Lambda \in \nat_\infty,\ \st {\hyper S(\Lambda, \Omega)} = \st {\hyper S(\Omega)}.$ From (ii), we have that, in like manner, $\st {\hyper S(\Lambda, \Omega) }= \st {\hyper S(\Lambda, \Gamma)}$ for each $\Gamma, \Omega \in \nat_\infty.$ Hence,$\hyper S(\Gamma) - \hyper S(\Omega) \in \mu (0).$ Thus, $ \hyper S(\Lambda, \Gamma) - \hyper S(\Gamma) \in \mu (0),\ \hyper S(\Delta,\Omega)- \hyper S(\Omega) \in \mu (0),$ which implies that 
$\hyper S( \Delta,\Omega) - \hyper S(\Lambda, \Gamma) \in \mu (0)$ for all $\Delta, \Omega, \Lambda, \Gamma \in \nat_\infty.$ Hence, $\lim S(m,n) =L = \lim_n(\lim_m S(m,n)) \to L\in \real$ by Theorem 5.4. Now apply Theorem 5.7 and the proof is complete. \qed \parm 
Although Theorem 5.10 is a necessary and sufficient condition, it's often difficult to apply from the knowledge of the iterated limit behavior. There are, as one would expect, special classes of double sequences where convergence of an iterated limit implies that the double limit converges. Many double sequences can be put into a form $S(m,n)\colon \nat \to \real$, where 
$\lim_m S(m,n) \to 0$ for each $n \in \nat$ and for each $m \in \nat$, $S(m,n)$ is decreasing [resp. increasing] in $n.$\parm
{\bf Theorem 5.11.} {\it Let $S(m,n)\colon \nat \times \nat \to \real$ and $\lim_m S(m,n)=0,$ for each $n \in \nat$  
and $S(m,n)$ is decreasing \r [resp. increasing \r ] in $n$ for each $m \in\nat.$ Then $S(m,n) \to 0.$}\pars
Proof. I show this for the decreasing case since the increasing case is established in like manner. We have that $\lim_mS(m,n) = S'(n) = 0$ for each $n \in \nat$. Thus, $\st {\hy S(\Lambda, \hy n)} = S'(n) = 0$ for each  $\Lambda \in \nat_\infty.$  Now $\lim_n S'(n) =0$ implies since, $S'$ is decreasing, that $0 \leq S'(n)  = \st {\hy S(\Lambda, \hy n)}$ for each $n \in \nat$. Thus, in general, either $\Hy 0 < \hy S(\Lambda, \Omega)$ for $\Omega \in \nat_\infty$ or $\hy S(\Lambda,\Omega) \in \mu(0)$. But, if $\Hy 0 < \hy S(\Lambda,\Omega) \leq \hy S(\Lambda, \hy n),$ then $0 \leq \st {\hy S(\Lambda, \Omega)} \leq \st {\hy S(\Lambda, \hy n)} = 0$. Hence, $S(m,n) \to  0$ and the proof is complete. \qed\parm 

The real numbers are complete. Hence, any nonempty set $A \subset \real$ that is {{\bf bounded above}} (i.e. there is some $y \in \real$ such that $x \leq y$ for each $x \in A$) has a {{\bf least upper bound}} that is denoted by $\sup A.$ This means that $\sup A$ is an upper bound and if $y \in\real$ is an upper bound for $A$, then $\sup A \leq y.$ The {{\bf greatest lower bound}} $\inf A$ exists for any nonempty $B \subset\real$ that is {{\bf bounded below}}. These ideas are applied to sequences that have convergent subsequences. Indeed, if $S[\nat]$ (the range of $S$) is bounded above [resp. below], them $\sup S[\nat]$ [resp. $\inf S[\nat]$] is an accumulation point and there is a subsequence that converges to this point.  (I'll show in the proof of Theorem 5.14 (ii) a method that you can modify to establish this result.) \parm 

{\bf Definition 5.12} ({\bf lim, inf, lim, sup.}) Given the sequence $S \colon \nat \to \real$. Let $y \in E$ iff there is a subsequence $S'$ of $S$ that converges to $y.$  The {{\bf lower limit}} (for $S$) $\lim \inf S_n = \inf E$, and the {{\bf upper limit}} (for $S$) $\lim \sup S_n = \sup E.$ \parm

Notice that if $S_n$ has no upper bound [resp. lower bound], then there is a subsequence $S'_n$ such that $S'_n \to +\infty$ [resp. $S'_n \to -\infty$]. In order to consider subsequences that diverge in this $\pm \infty$ special sense, the two symbols $-\infty,+\infty$ are included in the set $E$ and if $-\infty \in E$ [resp. $+\infty \in E$], then, by symbolic definition, let $\inf E = -\infty$ [resp. $\sup E = +\infty$], and no other cases need to be defined for sequences. We know that $S$ has a subsequence that converges to $L\in \real$ iff $\hyper S[\nat_\infty] \cap  \mu (L)\not= \emptyset.$   Further, the answer to the question I asked immediately after Theorem 4.14 is yes. So, because of this, the definition of $\St$ can be extended to the case where a subsequence  diverges to $\pm \infty.$ For $\Lambda \in \nat_\infty$, if $\hyper S(\Lambda) \in \real^{\pm}_\infty,$ let $\st {\hyper S(\Lambda)} =\pm \infty.$ By Theorem 4.10, the following result clearly holds. \parm
{\bf Theorem 5.13} {\it Let $S \colon \nat \to \real$. Then $\lim \inf S_n = \inf \{\st {\hyper S(\Lambda)}\mid \Lambda \in \nat_\infty\}$ and $\lim \sup S_n =\sup \{\st {\hyper S(\Lambda)}\mid \Lambda \in \nat_\infty\}.$}\parm
{\bf Theorem 5.14} {\it Let $S \colon \nat \to \real.$ Then\par
\indent\indent {\rm (i)} $\lim \inf S_n = - \infty$ \r [resp. $\lim \sup S_n =+\infty$\r ] iff there exists some $\Lambda \in \nat_\infty$ such that $\hyper S(\Lambda) \in \real^-_\infty$ \r [resp. $\real^+_\infty$\r ] iff $\hyper S[\hypernat \cap \real^-_\infty] \not= \emptyset$ \r [resp. $\real^+_\infty$\r ];\par
\indent\indent {\rm (ii)} $\lim \inf S_n = L \in \real$ \r [resp. $\lim \sup S_n$\r ] iff there exists some $\Lambda \in \nat_\infty$ such that $\hyper S(\Lambda) \in \mu (L)$ \r (i.e. $\st {\hyper S(\Lambda)} = L$\r ) and for each $\Omega \in \nat_\infty,$ $\hyper S(\Omega) \in \mu (L)$ or $\hyper S(\Omega) > \hyper S(\Lambda)$ \r [resp. $<$\r ].}\pars 
Proof. (i) Let $\lim \inf S_n = -\infty.$ This implies that for each $y \in \real^-$ and for each $m \in \nat$ there exists some $n \in \nat$ such that $n \geq m$ and $S_n < y.$ It should be obvious by now that such a statement means in our nonstandard structure that for any $a \in \real^-_\infty$ and $\Lambda \in \nat_\infty$ there is a $\Omega \in \nat$ such that $\Omega \leq \Lambda$ and, hence, $\Omega \in \nat_\infty$ such that $\hyper S(\Omega) < a.$ \pars
For the sufficiency, let $y \in \real^-,\ m \in \nat$. Then we know that if $a \in \real^-_\infty$, then $a < y.$ The  hypothesis states that the sentence 
$$\exists \r x((\r x \in \hypernat)\land (\r x > \hyper m) \land (\hyper S(\r x) < \Hy y))$$
holds in $\hyper {\cal M};$ and, hence, in ${\cal M}$. In like manner, for the $\sup.$ Thus (i) is established.\pars
(ii)  Since $\lim \inf S_n = L$, the set $E$ contains, at least one real number. I'll show that $L \in E.$ What we do know is that there is a subsequence of $S_n$ that converges to some number $\geq L.$ Hence, there exists $\Lambda \in \nat_\infty$ such that $\hyper S(\Lambda) \in G(0).$ Let $P = \{ \st {\hyper S(a)} \mid (a \in \nat_\infty)\land (\hyper S(a) \in G(0)\}\not= \emptyset.$ Now $\lim \inf S_n = \inf P = L.$ 
From definition of ``inf,'' if real $r > L$, then there is some $p \in P$ such that $0\leq p-L < r-L.$ Thus, let $0 < r - L = 1/(2n),\ n \in \nat, \ n \not= 0$. Then there is some $p(n) \in P$ such that $0 \leq p(n) - L < 1/(2n).$ Since $p(n)$ is the limit of a subsequence $Q$, then there exists some $m \in \nat$ such that $\vert Q(m) - p(n)\vert < 1/(2n).$ Since $Q_n \in S[\nat]$ then by defining $Q_m = Q'_n$, we have that $Q'$ is a subsequence of $S$ such that $\vert Q'(n) -L \vert <1/n$, for each nonzero $n \in \nat$. Hence, $Q'_n \to L$ implies that $L \in P.$ Of course, $P$ = $E$, as $E$ was previously defined.\pars
Now, there exists $\Lambda \in \nat_\infty$ such that $\hyper S(\Lambda) \in \mu (L).$ Assume that $\Omega \in \nat_\infty$ and $\hyper S(\Omega) \notin \mu (L).$  Since $L \not= -\infty$, then (i) implies that $\hyper S(\Omega) \notin \real^-_\infty.$ Further, note that $\hyper S(\Omega) \notin \mu (r)$ for any real $r < L$, since if this was so than there would be a subsequence of $S$ that converges to $r$ and then this contradicts the notion of ``inf.'' 
Thus, in this case, $\hyper S(\Omega) > \hyper S(\Lambda)$ (recall the monads are disjoint). The ``sup'' follows in like manner and the proof is complete. (The sufficiency is left to reader.)\qed \parm
{\bf Corollary 5.15.} {\it For a given $S \colon \nat \times \nat \to \real$, let $E$ contain the limits of each converging subsequence. Then 
$\inf E \in E$ \r [resp. $\sup E \in E$\r ].}\pars
Proof. This is established in the above proof for real valued ``inf'' and ``sup.'' Now obviously by definition, it also follows for the two defined cases of $\pm \infty.$ \parm\vfil\eject
{\bf Example 5.16.} Let $S \colon \nat \to \real.$ \pars
\indent\indent (i) Define $S_n = (-1)^n(1 + 1/(n+1)).$ Let $\Lambda \in \nat_\infty$ be a *-odd number. Then $\hyper S(\Lambda) = -(1 + 1/(\Lambda +1)) \in \mu (-1).$ Then taking a *-even $\Omega$, it's seen that $\hyper S(\Omega) \in \mu (1).$ Since for each $n\in \nat$, $-1 \leq S_n \leq 1$, we have that $\lim \sup S_n = 1,\  \lim \inf S_n = -1.$ From Theorem 5.14, we also know that for each $\Gamma\in \nat_\infty$ that $\hyper S(\Gamma) \in \mu (-1)$ or $ \hyper S(\Gamma) \in \mu (1)$ or $\hyper S(\Lambda) < \hyper S(\Gamma)$ or $\hyper S(\Gamma) < \hyper S(\Omega).$\pars
\indent\indent (ii) Let $S_n$ be the sequence of all the rational numbers. (Yes, technically there is such a sequence.) Then simply from noticing that there exist negative and positive infinite rational numbers, we have that $\lim \inf S_n = -\infty,\ \lim \sup S_n = +\infty.$\pars
\indent\indent (iii) Let $S_n$ and $Q_n$ be any two sequences. Then
$$\lim \inf S_n + \lim \inf Q_n \leq \lim \inf(S_n + Q_n) \leq $$ $$
\lim \sup (S_n + Q_n) \leq \lim \sup S_n + \lim \inf Q_n.$$
Proof. Let $A = \{\st {\hyper S(\Lambda)}\mid \Lambda \in \nat_\infty) \},\ B = \{ \st {\hyper Q(n)}\mid \Lambda \in \nat_\infty\}.$ If $A$ and $B$ are both nonempty, then nonempty $\{\st {\hyper S(\Lambda) + \hyper Q(\Lambda)} = \st {\hyper S(\Lambda)} + \st {\hyper Q(\Lambda)}\mid \Lambda \in \nat_\infty\} = A + B,$ where this ``addition'' definition is obvious. The result now follows from the ``well known'' result (taking into account the $\pm \infty$ possibilities) that $\inf A + \inf B \leq \inf(A + B) \leq \sup (A + B) \leq \sup A + \sup B.$ \pars
\indent\indent (iv) Let $S_n \to L \in \real$ and $Q_n$ be any sequence. Then $\lim \inf (S_n + Q_n) = L + \lim \inf Q_n,\  \lim \sup (S_n + Q_n) = L + \lim \sup Q_n.$\pars
Proof. Let $A = \{\st {\hyper Q(\Lambda)} \mid \Lambda \in \nat_\infty \}.$ By a trivial proof, when we use the symbols $\pm \infty$, we mean that they correspond to any member of $\real^{\pm}_\infty$ it follows symbolically that for any $a \in G(0),\ \pm \infty + a = \pm \infty.$ Recall how the definition of the $\St$ operator has been extended to $\pm \infty$. For any $a \in \hyperreal$, $\st {a} = \pm \infty$ iff $a \in \real^{\pm}_\infty.$ Thus under this definition $A \not= \emptyset.$  This definition also satisfies the usual extend algebra for $\pm \infty.$ Now we know that for each $\Lambda \in \nat_\infty,\ \hyper S(\Lambda) \in \mu (L).$ Under this extended definition, it follows that for each $\Omega \in \nat_\infty\ \st {\hyper S(\Omega) + \hyper Q(\Lambda)} = \st {\hyper S(\Omega)} + \st {\hyper Q(\Lambda)} = L + \st {\hyper Q(\Lambda)}.$ The result follows as in the proof of (ii), that $\inf A + L = \inf B.$ The ``sup'' part follows in like manner. \pars
\indent\indent (v) Let $S(m,n) \to L \in \real.$ Then $\lim_n(\lim_m \inf S(m,n)) = \lim_n(\lim_m \sup S(m,n)) = \lim_m(\lim_n \inf S(m,n)) = \lim_m(\lim_n \sup S(m,n))= L.$\pars
Proof. Now this can be established by nonstandard means. But, it's immediate from the fact that $S(m,n) \to L$ iff every subsequence $S'(m,n) \to L$ by just considering $n$ or $m$ as fixed. \qed 
\vfil\eject
\centerline{\bf 6. BASIC INFINITE SERIES CONCEPTS}
\parm
Sometimes it's useful to simplified the notation for the finite and infinite series.
Let $A(n) = \sum_{k = 0}^n a_k=\sum_{0}^n a_k ,$ there $k \in \nat.$ 
Then by definition this series converges to $L$ iff $A(n) \to L.$ Hence, all of our previous nonstandard characteristics for sequential convergence apply. For example, $A(n) \to L$ iff $\hyper A(\Lambda) \in \mu (L)$ for each $\Lambda \in \nat_\infty.$ Notationally, you will also see this written as $\sum_0^\Lambda \hy a_k \approx L = \hy L.$ These are called {{\bf hyperfinite}} or {{\bf *-finite}} summations. Indeed, any set such as $\{n\mid (\Hy 0 \leq n \leq \Lambda)\land (n \in \hypernat)\},$ where $\Lambda \in \nat_\infty$ is a *-finite set. The reason it's termed *-finite is that *-finite sets satisfy any of the finite set properties that can be presented in our formal language. To show that most of the basic manipulations done with a finite series hold for *-finite series, it's necessary to give a more formal definition for infinite series than is usually presented.\parm
{\bf Definition 6.1.} Let $a\colon \nat \to \real.$ Then the partial sum function 
$A(k)$ is defined inductively. \pars
\indent\indent (i) Let $A(0) = a_0;$\pars 
\indent\indent (ii) then $A(k+1) = A(k) + a_{k+1},\ k \in \nat.$ \pars 
\indent\indent (iii) Further, define for any $n,m \in \nat,\ n< m,\ A(n,m) = A(m) -A(n)= \sum_{n+1}^m a_k$ and $A(-1,0) = a_0$ and if $m = n \not=0,$ then $A(n,n) = 0.$ Notice that $A(n-1,n) = a_n$ in all cases.    
\parm
Observe $A \colon \nat \to \real.$ Thus, there's the nonstandard extension of this function to $\hy A \colon \nat \to \hyperreal.$ Further, we know that for $n \leq m,\ n,k \in \nat$, $A(m) = A(n) + A(n,m)= A(n,m) + A(n).$ This property also holds for $\Lambda \leq \Omega, \ \Lambda,\Omega \in \nat_\infty.$ But not every ordinary mathematical process that can be done will hold in $\cal M$ for $\hy A$. Whatever holds must be expressible in our formal language. This is not always possible. One thing that cannot be so expressed, generally, is the notion of ``any rearrangement'' of the members of a infinite series. What is needed is a specifically stated rearrangement. For each example, define $\Theta_n(k) = n -k$ for $k \in [0,n].$ Now applying this to the finite sequence of terms for our finite series $a_0 + \cdots +a_n$ yields $b_0=a_n + \cdots + b_n=a_0.$ This can be viewed as a new sequence, and by *-transform, it has meaning for any $\Lambda\in \nat_\infty.$  \pars
Because of how the ``term generating'' function is defined, it may be convenient to assume that the first few terms, say $a_0,a_1,a_2,\ldots, a_k,\ k < n,$  all equal zero. I assume that all sequences $a_k$ that are defined for $n \in \nat,$ where $n > k,$ are extended, if necessary, to sequences defined on $\nat$ by letting $a_i = 0,\ 0\leq i \leq k.$ (There will be times when this is not done and the notation will indicate this.) Moreover, it's also clear that removing finitely many terms from a series does not effect whether it converges or not. Thus, given  original $A \colon \nat \to\real,$ to determine whether $A(n)$ converges you can use a different $B\colon \nat \to \real$ obtained by letting $b_n = a_{k+n}$ for any fixed $k >0$ and then investigate the convergence of $B(n).$ Of course, you would need to adjust the two limits if they do converge. 
However, mostly, one is interested in the terms of a series, the $a_k.$ Further, note that, for $\Lambda,\Omega \in \nat_\infty,\ \Lambda < \Omega,$ $\hy A(\Lambda, \Omega) = \sum_{\Lambda +1}^\Omega \hy a_k.$ \parm
{\bf Theorem 6.2} {\it Let $a\colon \nat \to \real$ be a bounded sequence. Then there exists $\hy M \in \sig \real$ such that for each $\Lambda, \Omega \in \nat_\infty,\ \Lambda \leq\Omega$ 
$$\left|{{\sum_\Lambda^\Omega \hy a_k }\over{\Omega - \Lambda + \Hy 1}}\right| \leq \hy M.$$} \parm
Proof. Since $a_k$ is bounded, then there exists some $M \in \real$ such that $\vert a_k\vert \leq M.$ One of the most significant results for finite series is that for $n \leq m,$ $\vert \sum_n^m a_k \vert \leq \sum_n^m \vert a_k\vert.$ By defining the sequence $b_k = \vert a_k \vert$, then by *-transform it follows that for $\Lambda,\Omega \in \nat_\infty,\ \Lambda \leq \Omega,\ \vert \sum_\Lambda^\Omega \hy a_k \vert \leq  \sum_\Lambda^\Omega \Hy b_k = \sum_\Lambda^\Omega \vert \hy a_k \vert.$ Now for the standard series since $\vert a_k \vert \leq M,$ then $\sum_n^m a_k\leq M(m-n +1).$ By *-transform, it follows that for $\Lambda, \Omega \in \nat_\infty,\ \Lambda \leq \Omega,$ $\vert \sum_\Lambda^\Omega a_k\vert \leq \hy M (\Omega-\Lambda + \Hy 1)$ and the result follows. \qed \parm 
I restate some of the previous nonstandard sequence results that now characterize convergence of the series $A(n).$ \parm
{\bf Theorem 6.3.} {\it The series $A(n) \to L$ iff $\hy A(\Lambda) - L \in \mu (0),$ for each $\Lambda \in \nat_\infty$ iff $\hy A(\Lambda) \in \mu (L),$ for each $\Lambda \in \nat_\infty$ iff $\st {\sum_0^\Lambda a_k} = L$ for each $\Lambda \in \nat_\infty$ iff $\hy A[\nat_\infty] \subset \mu (L).$}\parm
{\bf Theorem 6.4.} {\it {\rm (i)} $A(n) \to L$ iff for each $\Lambda, \Omega \in \nat_\infty,\ \Lambda \leq \Omega,\ \sum_\Lambda^\Omega a_k\in \mu (0)$ iff $\hy A(\Lambda,\Omega) \in \mu (0).$ {\rm (ii)} If $A(n) \to L$, then $a_\Omega \in \mu (0),$ for each $\Omega \in \nat_\infty.$ }\pars
Proof. (i) First, note that if $\lambda = \Omega,$ then $\hy A(\Lambda,\Omega) = 0$ and $\hy A(\Omega -1,\Omega) =\hy A(\Omega) - \hy A(\Omega -1) =\sum_\Omega^\Omega a_k = a_\Omega.$ If $\Lambda < \Omega,$ then $\hy A(\Lambda,\Omega) =\hy A(\Omega) - \hy A(\Lambda) =\sum_{\Lambda+1}^\Omega a_k.$ If $\Lambda > \Omega,$ then $\hy A(\Omega) - \hy A(\Lambda) =- \sum_{\Lambda+1}^\Omega a_k.$ The result, in general, comes from the Cauchy Criterion and, clearly, we may assume that $\Lambda\leq \Omega.$ (ii) This is immediate. \qed\parm
Although it's not required in our investigations, the converse of Theorem 6.4 (ii) holds for certain series.
Recall what G\"odel wrote, that removing one quantifiers from a characterization is significance. So far, the nonstandard characterization do just this and often remove all quantifiers.  \parm 
{\bf Theorem 6.5.} {\it If $A(n) \to L,$ then $a_n \to 0.$}\pars
Proof. From Theorem 6.4 (ii). \qed \parm
{\bf Example 6.6.} \pars
\indent\indent (i) Let $a_k = {{1}\over{(k+1)(k+2)}}.$ Let $\Lambda \in \nat_\infty.$ Then $\hy A(\Lambda) = \sum_0^\Lambda {{ 1}\over{(k+1)(k+2)}}=\sum_0^\Lambda {{1}\over{k+1}} - \sum_0^\Lambda {{1}\over{k+2}} = 1 + \sum_1^\Lambda {{1}\over{k+1}} - \sum_1^\Lambda {{1}\over{k+1}} - {{1}\over{\Lambda +2}}= 1 + 0 - {{1}\over{\Lambda + 2}}$ by *-transform of the finite case and I have applied the convention of writing $\hy x = x$ for $\hy x \in \sig \real.$ But, ${{1}\over {\Lambda +2}} \in \mu (0).$ Hence $A(n) \to 1$ or, as is often written, $\sum_0^\infty a_k = 1.$\pars
(ii) Let $A(x) = \sum_0^\infty x^k,\ x\not= 1.$ We know that, in general, $a_k(x) = {{1-x^{k+1}}\over{1-x}}.$ Hence, $\hy a_\Lambda = {{1-x^{\Lambda +1}}\over{1-x}}= {{1}\over{1-x}} + {{x^{\Lambda +1}}\over{1-x}}$ for $x \in \sig \real.$  If $\vert x\vert <1,$ then ${{x^{\Lambda +1}}\over{1-x}} \in \mu (0)$ implies that $\hy a_\Lambda (x) \in \mu ({{1}\over{1-x}}).$ On the other hand, if $\vert x \vert > 1,$ then ${{1}\over{1-x}} - {{x^{\Lambda +1}}\over{1-x}} \notin G(0),$ for $\Lambda \in \nat_\infty.$ Thus, the series diverges. \parm
When compared with a general series, it's often easier to show that a non-negative type series converges or diverges. A series $A(n)$ is  {{\bf non-negative}} iff there is some $m\in \nat$ such that $a_k \geq 0$ for each $k \geq m$.\parm
{\bf Theorem 6.7.} {\it A non-negative $A(n)$ converges iff there is some $\Lambda \in \nat_\infty$ such that $\Hy A(\Lambda) \in G(0)$ iff $\hy A[\hypernat] \subset G(0).$ }\pars
Proof. There is some $m \in \nat$ such that $a_k \geq 0$ for all $k \geq m.$ 
Thus we write $A(n)$ in terms of a new sequence $B(n),$ where $A(n) = B + B(n)$ and $B(n)$ is an increasing sequence.  The result follows from Theorems 4.2 and 4.7 and the fact that $\Hy B(\Lambda) \in G(0)$ iff  $B + \Hy B(\Lambda) \in G(0)$ for any $B \in \real.$ And, an increasing sequence converges iff it is bounded. The result follows for $A$ is bounded iff $B$ is bounded. \qed \parm
{\bf Theorem 6.8.} {\it A non-negative $A(n)$ diverges iff there is some $\Lambda \in \nat_\infty$ such that $\Hy A(\Lambda) \notin G(0)$ iff $\hy A(\hypernat)\not\subset G(0).$ }\pars
 
There are standard arithmetical results that can aid in determining convergence.\parm
{\bf Theorem 6.9.} {\it Given $A \colon \nat\to\real.$ Let $f\colon \nat \to \nat$ have the property that, for each  $A(f(n+1)) -A(f(n)) \geq b_n$ for each $n \geq j$ \r [resp. $\leq$\r ]. Then for each $n \in \nat,\ n \geq j$} $$A(f(n+1)) \geq A(f(j)) + \sum_j^n b_k,\ {\rm [resp.\ \leq\ ].}$$ \pars
Proof. [For $\geq$.] For $n = j,$ clearly, the result holds. So, assume it for $m > j.$ Then consider $m+1.$ Since $A(f(m+2)) - A(f(m+1)) \geq b_{m+1}$, then $A(f(m+2)) \geq A(f(m+1)) + b_{m+1} \geq A(f(j))+ \sum_j^n b_n+ b_{m+1} = A(f(j)) + \sum^{m+1}_jb_k$ the result holds by induction. \qed\parm 
{\bf Example 6.10.} Let's determine convergence or divergence directly for two very well know series. In all cases, we extend any series to include the necessary zero terms although they might be directly mentioned.\pars 
(i) Consider $b_k = 1/k,\ k >0.$ We look at the series $a_0=0,\ a_k = 1/k,\ k \geq 1.$ Define $f\colon \nat \to \nat$ by letting $f(n) = 2^n.$ Then, for $n
\geq 1,$  
$$A(2^{n+1}) - A(2^{n}) = \sum_{2^n+1}^{2^{n+1}}{{1}\over{k}} = \sum_1^{2^n} 
{{1}\over{2^n +k}} \geq {{2^{n}}\over{2^{n+1}}} = 1/2.$$ Applying Theorem 6.9, $A(2^{n+1}) \geq A(1) + \sum_1^n (1/2) = 1/2 + n/2.$ Consequently, for $\Lambda \in \nat_\infty,$
$A(2^{\Lambda +1}) \in \real_\infty$ and the series diverges. \pars
(ii) Now let's look at famous ``p'' series, where $b_k = 1/k^p, \ k >0 .$  (a) First, let $p >1$ and look at the series $a_0=0$ and  $a_k = 1/k^p,\ k \geq 1.$ As done for (i) $A(2^{n+1}) - A(2^{n}) = \sum_1^{2^{n}}{{1}\over{(2^n +k)^p}}.$  I note that each term of this sum is less than $2^{-pn}$ and there are $2^n$ terms.  Hence $\sum_1^{2^{n}}{{1}\over{(2^n +k)^p}} < {{2^{n}}\over{2^{np}}}$ implies from Theorem 6.9 that 
$$0 \leq A(2^{n+1}) \leq 1/2^p + \sum_1^n {{2^k}\over{2^{kp}}}.$$ 
But, since $p>1, \ \hy A(2^{\Lambda +1}) \in G(0)$ and the non-negative series converges.\pars\
Now, for $0 \leq p\leq 1,$ each term of the finite sum ${{1}\over{(2^n +k)^p}} \geq {{1}\over{2^{(n+1)p}}}.$ Hence, as done above the $$A(2^{n+1}) -A(2^{n})=\sum_1^{2^n}{{1}\over{(2^n +k)^p}}\geq {{2^n}\over{2^{(n+1)p}}} \geq {{1}\over{2^p}}.$$ 
Consequently from Theorem 6.9, $A(2^{n+1}) \geq 1/2^p + {{n}\over{2^p}}$ and for this non-negative series $\hy A(2^{\lambda +1}) \notin G(0)$ and the series diverges. It obviously diverges for all $ p <0.$\parm
All of the standard converges or divergence tests can be translated into appropriate nonstandard statements. However, here is an interesting nonstandard comparison test.\parm
{\bf Theorem 6.11.} {\it Let $\sum_0^\infty a_k$ be a non-negative. If non-negative $B(n)$ converges and there is some $c \in \hyperreal$ such that $0\leq c$ and $c \in G(0)$ and, for each $\Lambda \in \nat_\infty,\ \hy a(\Lambda) \leq c(\Hy b(\Lambda))$, then $A(n)$ converges.}\pars
Proof. Assume that $B(n)$ converges. Then for each $\Lambda,\Omega \in \nat_\infty$ such that $\Lambda \leq \Omega, \ \sum_\Lambda^\Omega\Hy b_k \in \mu (0).$ Also there is some $ r \in \real^+$ such that $c\leq \Hy r.$ Hence $\Hy 0 \leq \hy a_\Lambda \leq c (\Hy b(\Lambda))\leq \Hy r (\Hy b(\Lambda)),$ for each $\Lambda \in \nat_\infty.$ By *-transform of the finite case, this implies that $0\leq \sum_\Lambda^\Omega\hy a_k \leq \sum_\Lambda^\Omega \Hy r (\Hy b_k) = \Hy r \sum_\Lambda^\Omega \Hy b_k \in \mu (0).$ This completes the proof. \qed \parm
{\bf Theorem 6.12.} {\it If non-negative $\sum_0^\infty b_k$ diverges and there exists $c >0,\ c \in \hyperreal -\mu (0)$ and $\hy a(\Lambda) \geq c(\Hy b(\Lambda))$ for each $\Lambda \in \nat_\infty,$ then $\sum_0^\infty a_k$ diverges.}\pars
Proof. Since $\sum_0^\infty b_k$ diverges, then there exist $\Lambda, \Omega,\ \Lambda  \leq \Omega$ and $\sum_\Lambda^\Omega \Hy b_k \notin \mu (0).$ We also know that there exists some $r \in \real^+$ such that $\Hy r < c$. Therefore $\Hy r \sum_\Lambda^\Omega  \Hy b_k = \sum_\Lambda^\Omega \Hy r\hy b_k\notin \mu (0).$ Hence, since $\sum_\Lambda^\Omega \hy a_k  \geq \sum_\Lambda^\Omega \Hy r\Hy b_k > 0,$ then $\sum_\Lambda^\Omega \hy a_k \notin \mu (0)$ and the result follows. \qed \parm 
{\bf Example 6.13.} Assume that $\sum_0^\infty a_nu^k$ converges for $u \not=0.$ Then $\sum_0^\infty a_kx^k$ converges absolutely for each $x$ such that $\vert x \vert < \vert u \vert.$\pars
Proof. Let $b = \vert x\vert/\vert u \vert < 1.$ Hence, the geometric series $\sum_0^\infty b^k$ converges for such an $x$. Let $\Lambda \in \nat_\infty$ Then $\hy a(\Lambda)u^\Lambda \in \mu (0)$, since $\sum_0^\infty a_nu^k$ converges. Notice that $$\vert \hy a(\Lambda) x^\Lambda \vert = \left|\hy a(\Lambda)u^\Lambda {{x^\Lambda}\over{u^\Lambda}}\right| = \vert \hy a(\Lambda) u^\Lambda \vert b^\Lambda < b^\Lambda,$$ since $\vert \hy a(\Lambda) u^\Lambda \vert < 1.$ You can apply Theorem 6.11, where $c = 1.$ \parm
In the Chapter ``Series of nonnegative terms'' (1964, p. 55) W. Rudin states that ``One might thus be led to the conjecture that there is a limiting situation of some sort, a `boundary' with all the convergent series on one side, all the divergent series on the other side - at least as far as a series with monotonic coefficients are concerned. This notion of `boundary' is of course quite vague. The point we wish to make is this: No Matter how we make this notion precise, the conjecture is false.'' However, Rudin's statement using the phrase ``No matter how'' in this section on non-negative series is itself false. Theorems 6.7 and 6.8 show that $G(0)$ is just such a ``boundary.''\pars
Here is another example of the usefulness of the nonstandard methods and direct proofs.\parm 
{\bf Example 6.14.} Let each $a_k >0,$ $\sum_0^\infty a_k$ converge and $a_{k+1}\leq a_k$ for all $k\in \nat.$ Then $\lim na_n = 0.$\pars
Proof. Let $\Lambda \in \nat_\infty.$ For each non-negative real number $r$, there exists a unique natural number $[r]$ such that $[r] \leq r < [r] +1.$ This statement and property can be written in our formal language. Thus by *-transform, since $\Lambda /2\in \real_\infty,$ there exists $a > 0,\ a \in \hyperreal$ and $a = [\Lambda/2].$ Now $0\leq [\Lambda/2] - \Lambda/2 < 1.$ Hence, necessarily, $a = \Omega =[\Lambda/2] \in \nat_\infty.$ Then $\hy A(\Lambda) -\hy A(\Omega)\geq (\Lambda-\Omega)\hy a(\Lambda) \geq (\Lambda/2)\hy a(\Lambda) \geq \Hy 0$ since such a statement holds in $\real$. Thus, $\Lambda\hy a(\Lambda) \in \mu (0)$ and the result follows. \parm
\vfil\eject
\centerline{\bf 7. AN ADVANCED INFINITE SERIES CONCEPT}
\parm
Some of the most interesting aspects of the nonstandard theory of infinite series are developed when various {{\bf infinite series product notions}} are probed. But we need the following result {{\bf Abel's summation by parts.}}\parm
{\bf Theorem 7.1.} {\it Let series $A \colon \nat \to \real,\ B\colon \nat \to\real.$ Then for each $p,q \in \nat,\ p\leq q$
$$\sum_p^qa_kb_k = \sum_p^q(A(k) - A(k-1))b_k = \sum_p^qA(k)b_k - \sum_{p-1}^{q-1}A(k)b_{k+1} =$$ $$ \sum_p^qA(k)(b_k-b_{k+1}) - A({p-1})b_p + A(q)b_{q+1},$$ where $A_{-1} = 0.$}\parm  
{\bf Theorem 7.2.} {\it For each $\Lambda,\Omega \in \nat_\infty,\ \Lambda \leq \Omega$
$$\sum_\Lambda^\Omega \hy a_k\Hy b_k =\sum_\Lambda^\Omega A(k)(\Hy b_k-\Hy b_{k+1}) - \hy A({\Lambda-1})\Hy b(\Lambda) + \hy A(\Omega)\Hy b({\Omega+1}).$$}\pars
Proof. By *-transform. \qed\parm
One can immediately induce upon the right-hand side of the equation in Theorem 7.2 various requirements that will force it to be an infinitesimal. 
This will be seen in the proof of Theorem 7.3. But, first, notice that for the collapsing series $\sum_0^\infty (b_k - b_{k+1})$, we have for each $\Lambda,\Omega \in \nat_\infty,\ \Lambda \leq \Omega, $ if $\vert\sum_\Lambda^\Omega \Hy b_k - \Hy b_{k+1} \vert =  \vert \Hy b(\Lambda) - \Hy b(\Omega +1) \vert \leq \sum_\Lambda^\Omega \vert \Hy b_k - \Hy b_{k+1}\vert \in \mu (0),$ then  $\Hy b(\Lambda) - \Hy b(\Omega +1) \in \mu (0),$ for each such $\Omega$ and $\Omega.$ \parm
{\bf Theorem 7.3.} {\it If $\sum_0^\infty (b_k - b_{k+1})$ converges absolutely and $A\colon \nat \to \real$ is bounded, then $\sum_0^\infty a_kb_k$ converges.}\pars
Proof. From Theorem 7.2, $$\vert\sum_\Lambda^\Omega \hy a_k\Hy b_k\vert \leq\sum_\Lambda^\Omega \vert A(k)(\Hy b_k-\Hy b_{k+1})\vert  +\vert - \hy A({\Lambda-1})\Hy b(\Lambda) + \hy A(\Omega)\Hy b({\Omega+1})\vert.$$
Since there is some $r \in \real^+$ such that $\vert \hy A(\Gamma) \vert 
\leq r$ for each $\Gamma \in \nat_\infty,$ then
$$\vert \sum_\Lambda^\Omega \hy a_k\Hy b_k \vert \leq r\left(\left(\sum_\Lambda^\Omega \vert \Hy b_k - \Hy b_{k+1} \vert\right) +\vert \Hy b(\Omega + 1) - \Hy b(\Lambda -1)\vert \right).$$ Since $\sum_0^\infty(b_k - b_{k+1})$ converges absolutely, then $\vert \Hy b(\Omega +1)-\Hy b(\Lambda -1)\vert = \vert \sum_{\Lambda-1}^\Omega (\Hy b_k - \Hy b_{k+1})\vert \leq \sum_{\Lambda-1}^\Omega \vert \Hy b_k - \Hy b_{k+1}\vert\in \mu (0),$ $\sum_\Lambda^\Omega \vert \Hy b_k - \Hy b_{k+1} \vert \in \mu (0)$ and the result follows. \qed\parm
{\bf Corollary 7.4.} {\it Let $A\colon \nat \to \real$ be bounded. If $\sum_0^\infty (b_k -b_{k+1})$ converges and $\{b_k\}$ is decreasing, then $\sum_0^\infty a_kb_k$ converges.}\pars
Proof. Observe that for each $k \in \nat,\ b_k - b_{k+1} \geq 0$ and, hence, $\sum_0^\infty (b_k - b_{k+1})$ is absolutely convergent. The result follows from the previous theorem. \qed \parm
Now let's complete this chapter by investigating the ``Cauchy product'' and show how nonstandard methods aid intuition. I'll ``play around'' with the ``subscript'' notation somewhat and one needs to understand what the double summation symbol is actually trying to indicate. The inner most of the two summations symbols will always indicate a ``finite'' summation where the index limit symbol is considered as fixed. Thus the notation
$\sum_{k=0}^n\left( \sum_{j = 0}^ka_jb_{k-j}\right)$ means that you fixed each $k,\ 0\leq k$ and obtain the value of the finite sum $\sum_{j = 0}^ka_jb_{k-j}$. Then add all of the $n+1$ results together  to get the double summation. I won't go through what some consider to be an ``easy'' proof that for non-trivial $n\geq 1.$  
$$\left(\sum_0^n a_k\right)\left(\sum_0^n b_k\right) = \overbrace{\sum_{k =0}^n\left(\sum_{j=0}^k a_jb_{k-j}\right)}^C + \overbrace{\sum_{k=0}^{n-1}\left(\sum_{j=0}^k a_{n-k+j}b_{n-j}\right)}^{In}.\eqno (7.5)$$ 
In the above expansion, the double sum indicated by the $C$ is the most significant. Indeed, let $c_k = \sum_{j =0}^k a_jb_{k-j}$. This is often called the {{\bf Cauchy product.}} Then you have the sequence (i.e. series) $C(n) = \sum_0^n c_k = \sum_{k =0}^n\left(\sum_{j=0}^k a_jb_{k-j}\right).$ \parm
{\bf Theorem 7.6.} {\it Let $A(n) \to L_a$ and $B(n) \to L_b.$ Then $C(n) \to L_aL_b$ iff, for any $\Omega \in \nat_\infty,\  \sum_{k=0}^{\Omega-1}\left(\sum_{j=0}^k \hy a(\Omega-k+j)\Hy b(\Omega-j)\right) \in \mu (0)$}\pars 
Proof. From the hypotheses, $\left(\sum_0^\Omega\hy a_k)\right)\in \mu(L_a)$ and $\left(\sum_0^\Omega\Hy b_k)\right)\in \mu(L_b)$ for any $\Omega \in  \nat_\infty.$ Hence, $\left(\sum_0^\Omega\hy a_k\right)\left(\sum_0^\Omega\Hy b_k\right) \in \mu (L_aL_b).$ Now $\sum_0^\Omega c_k = \sum_{k =0}^\Omega\left(\sum_{j=0}^k \hy a_j\Hy b_{k-j}\right).$ But, $\sum_{k =0}^\Omega\left(\sum_{j=0}^k \hy a_j\Hy b_{k-j}\right) \in \mu (L_aL_b)$ iff $\sum_{k=0}^{\Omega-1}\left(\sum_{j=0}^k \hy a(\Omega-k+j)\Hy b(\Omega -j)\right)\in \mu (0).$ \qed \parm 
Although Theorem 7.6 indicates what portion of the right-hand side of equation (7.5) must be infinitesimal for the Cauchy product to equal the product of the limits of two converging series, this characterization is not the most useful. Using our previous notation, consider the sequences $A$ and $B$ and $C$. Suppose that $B(n) \to L_b.$ 
You should be able to show that for all $n \in  \nat,\ 
C(n) = A(n)L_b + \sum_0^na_k(B(n-k) -L_b).$ What is needed in the next few theorems is the notion of the maximum member of any nonempty finite set determined by a given sequence $Q \colon \nat \to \real.$ The following sentence holds in $\cal M.$
$$\forall \r x\forall \r y((\r x\in \nat)\land(\r y \in \nat)\land(\r x \leq \r y) \to \exists \r z((\r z \in \nat)\land (\r x \leq \r z\leq \r y)\land $$ $$\forall \r w((\r w\in \nat)\land (\r x\leq \r w\leq \r y) \to Q(\r w) \leq Q(\r z))))\eqno (7.7)$$
(Recall that $Q(x)\leq Q(z)$ is but a short-hand notation for $(Q(x),Q(y))$ being in the $\leq$ binary relation.) For any two $i,j \in \nat,\ i \leq j$, the $Q(w)$ is called the {{\bf maximum value}} in the nonempty finite set $\{Q(x)\mid i\leq x \leq j\}.$ It's denoted by $\max\{Q(x)\mid i\leq x \leq j\}.$ Further, under *-transform such a member of $\hyperreal$ exists for $\Lambda, \Omega \in \nat_\infty,\ \Lambda \leq \Omega.$ We use this to establish 
the following theorem.\parm
{\bf Theorem 7.8.} {\it Let $B(n) \to L_b.$ If $A(n) \to L_a$ absolutely, then 
$C(n) \to L_aL_b$ }\pars
Proof. Since $B(n)- L_b \to 0,$ then $\hy B(n) -\Hy L_b \in G(0)$ for each $n \in \hypernat.$ Also for each $r \in \real^+$ there is some $m \in \nat$ such that for each $n >m,\ n \in \hypernat,$ $\vert \hy B(n) -\Hy L_b \vert < \Hy r.$ Let $L_a = \sum_0^\infty \vert a_k \vert.$  For any $\Omega \in \nat_\infty,$ consider (in simplified notation) 
$$A_1 = \vert \sum_0^\Omega\hy a_k(\hy B(\Omega - k) - L_b)\vert \leq 
  \sum_0^{\Omega -(m+1)}\vert \hy a_k\vert\, \vert \hy B(\Omega-k) - L_b\vert + $$ $$\sum_{\Omega -m}^\Omega \vert \hy a_k\vert\, \vert \hy B(\Omega -k) -L_b \vert <  rL_a + $$ $$(m +1)\max\{\vert \hy a_k\vert\, \vert \hy B(\Omega -k) - L_b \vert \mid \Omega - m \leq k \leq \Omega \} = rL_a + (m+1)\eps,$$ where $\eps \in \mu (0),$ for $\Lambda -m\leq k \leq \Lambda$ implies that $\vert \hy a_k \vert \in \mu (0),$ which implies that $\{\vert \hy a_k\vert\, \vert \hy B(\Omega -k) - L_b \vert \mid \Omega - m \leq k \leq \Omega \} \subset \mu (0).$ (I have used the *-transform of (7.7).) But, $r$ is an arbitrary member of $\real^+$ implies that $A_1 \in \mu (0)$ and the result follows from Theorem 7.6. \qed \parm
What if $A(n) \to L_a,\ B(n) \to L_b,\ C(n) \to L_c,$ then does it follow that $L_c = L_aL_b$? In order to establish this, I establish, by nonstandard means, two special theorems that are useful for many purposes.\parm 
{\bf Theorem 7.9.} {\it If $S(n) \to L$, then $\lim_n{{\sum_1^ns_k}\over{n}} = L = \lim_n{{\sum_0^n a_k}\over{n+1}},$ where $a_k = s_{k+1}.$}\pars
Proof. For each $r \in \real^+$, $\vert\Hy S(\Lambda) - L \vert < r$ since for each $\Lambda \in \nat_\infty,\ \Hy S(\Lambda) - L \in \mu (0).$ So, consider arbitrary $r \in \real^+.$ Let $\Omega \in \nat_\infty,\ \rho = [\sqrt {\Omega}]$ as defined in Example 6.14. Then $\rho \in \nat_\infty.$
Moreover, $1/\Omega \leq 1/\rho^2.$ So, consider $$
\left|{{\sum_1^\Omega}\over{\Omega}} - L\right| \leq {{\sum_1^\Omega \vert \Hy S -L\vert}\over{\Omega}} \leq {{\sum_1^\rho \vert \Hy S_k - L \vert}\over{\rho}}{{1}\over{\rho}} + {{\sum_{\rho +1}^\Omega\vert \Hy S_k - L\vert}\over{\Omega}}\leq $$ $${{r}\over{\rho}} + {{\Omega-\rho}\over{\Omega}}\max\{\vert \Hy S_x - L\vert \mid \rho + 1\leq x\leq\Omega\}$$ I apply my previous discussion on the ``maximum'' object that exist in any such *-finite set.  Since $(\Omega -\rho)/\Omega < 1$ and all the objects in $\{\Hy \vert S_x - L\vert \mid \rho + 1\leq  x\leq\Omega\}$ are infinitesimals and $r/\rho$ is an infinitesimal, then the result follows. \qed \parm
{\bf Theorem 7.10.} {\it If $a_n \to A,\ b_n \to B,$ then 
$$\lim_n{{\sum_0^na_kb_{n-k}}\over{n+1}} = AB.$$}
Proof. For each (i.e. $\forall$), $\Lambda \in \nat_\infty,\ k \in \nat,$  let 
$$D_\Lambda = {{\sum_0^\Lambda \hy a(k)\Hy b(\Lambda - k)}\over{\Lambda +1}} = {{\sum_0^\Lambda \Hy b(\Lambda -k)(\hy a(k) -A)}\over{\Lambda +1}} + {{A\sum_0^\Lambda\Hy b(\Lambda -k)}\over{\Lambda +1}}.$$ From convergence, for some $M \in \real^+,\ \vert \Hy b(d-k)\vert \leq M$ for each $d\in \hypernat,\ d \geq k$. Hence,
 $\vert \Hy b(\Lambda- k)(\hy a(k) -A)\vert \leq M\vert(\hy a(k) -A)\vert,\ \forall\,  \Lambda \in \nat_\infty.$ From this and *-transform of the finite sum case,
$$0\leq E_\Lambda = {{\vert \sum_0^\Lambda \Hy b(\Lambda -k)(\hy a(k) -A)\vert}\over {\Lambda +1}} \leq M {{\sum_0^\Lambda \vert \hy a(k) -A \vert}\over{\Lambda +1}},\ \forall\, \Lambda \in \nat_\infty.$$
Since $a_n \to A,$ implies that $\vert a_n -A \vert \to 0,$ then from Theorem 7.9, 
$${{\sum_0^\Lambda \vert \hy a(k) -A \vert}\over{\Lambda +1}}\in \mu (0),\ \forall\ \Lambda \in \nat_\infty.$$
Therefore, $E_\Lambda \in \mu (0),\ \forall\,\Lambda \in \nat_\infty.$  
Consequently, 
$$D_\Lambda -A{{\sum_0^\Lambda \Hy b(\Lambda -k)}\over{\Lambda +1}} \in \mu (0),\ \forall\, \Lambda \in \nat_\infty.$$
Note that, in general, $\sum_0^nb_k = \sum_0^nb_{n-k},\ \forall \, n \in \nat.$ Thus, by Theorem 7.9,  ${{\sum_0^\Lambda \Hy b(\Lambda -k)}\over{\Lambda +1}} = {{\sum_0^\Lambda \Hy b(k)}\over{\Lambda +1}} \in \mu (B),\ \forall\, \Lambda \in \nat_\infty.$ Hence, $D_\Lambda -AB \in \mu (0),\ \forall \, \Lambda \in \nat_\infty$ and the result follows. \qed\parm   
{\bf Theorem 7.11.} {\it If $A(n) \to L_a,\ B(n) \to L_b, \ C(n)= \sum_{k =0}^n\left(\sum_{j = 0}^ka_jb_{k-j}\right) \to L_c,$ then  $L_c = L_aL_b.$}\pars
Proof. Recall, that $C(n) = A(n)L_b + \sum_0^na_k(B(n-k)-L_b).$ This can be re-expressed as $C(k) = \sum_{j = 0}^ka_jB(j-k).$ Then by re-arrangement of the terms, it follows that, in general, 
$$\sum_0^n c_n = \sum_0^nA(k)B(n-k),\ \forall\, n \in \nat.$$ 
Hence, by the previous two theorems,
$$\lim_n{{\sum_0^n c_k}\over{n+1}} = L_c = \lim_n{{\sum_0^nA(k)B(n-k)}\over{n+1}}= L_aL_b$$ a the result follows. \qed
\vfil\eject
\centerline{\bf 8. ADDITIONAL REAL NUMBER PROPERTIES}
\parm
Since this is supposed to be a monograph covering some of the basic notions in a first course in real analysis (i.e. calculus IV), then one should expect that certain additional real number properties need to be explored. This is  especially the case if a slight generalization of the notion of continuity and the like is investigated. You will discover that, once again, the monad is the nonstandard ``king,'' so to speak, in characterizing these concepts. There are slightly different definitions within the subject of ``point-set topology'' for the set-theoretic ``accumulation point.'' I has chosen to use a definition that makes this notion equivalent to the previous sequence definition.\parm
{\bf Definition 8.1.} Let $A \subset\real.$ Then $p \in \real$ is an {{\bf accumulation point}} of (for) $A$ iff, for every $w \in \real^+,$ the open interval $(-w + p, p +w)\cap A \not= \emptyset.$ A point $p \in \real$ is a {{\bf cluster point}} iff  for every $w \in \real^+,$ the {{\bf deleted open interval}} $= ((-w + p, p +w)- \{p\}) =(-w + p, p +w)' \cap A \not= \emptyset$ iff $(-w +p,p+w)' \cap A =$ an infinite set.\parm 
A cluster point is an accumulation point but not conversely. Consider the set $A = [1,2] \cap \{3\}.$ Then $3$ is an accumulation point, and not a cluster point. Also each member of a nonempty $A$ is an accumulation point.\parm
{\bf Definition 8.2.} The set of all accumulation points is called the {{\bf closure}} of the set $A \subset \real$ and is denoted by $\overline{A}$ or ${\rm cl}A.$\parm
Note that $A \subset {\rm cl}(A).$ \parm
{\bf Definition 8.3.} A point $p \in A\subset \real$ is an {{\bf interior point}} of $A$ iff there exists some $w \in \real^+$ such that $(-w + p, p +w) \subset A.$ \parm
{\bf Definition 8.4.} A point $p \in A  \subset \real$ is an {{\bf isolated point}} of $A$ iff there exists some $w \in \real^+$ such that $(-w + p,p+w) \cap A = \{p\}.$\parm
Notice that if $S\colon \nat \to \real$, then $p$ is an accumulation point for the sequence iff $p$ is an accumulation point for the set $S[\nat]$ (i.e the range). It also follows that $p$ is an accumulation point for $A$ iff there's a sequence $S$ of members of $A$ such that $S(n) \to p.$ Also a point $p \in \real$ is an isolated iff it is an accumulation point and not a cluster point. This last statement characterizes the difference between the notions of the accumulation point and cluster point. Cluster points are accumulation points that are not isolated. Now how do monads characterize this set-theoretic notions?\parm 
{\bf Theorem 8.5.} {\it Let $A \subset \real, \ p \in \real$. Then\pars
\indent\indent {\rm (i)} $p$ is an accumulation point iff $\mu (p) \cap \hy A \not=\emptyset$;\pars
 \indent\indent {\rm (ii)} $p$ is an isolated point iff $\mu (p) \cap \hy A = \{p\}$;\pars
\indent\indent {\rm (iii)} $p$ is a cluster point iff the deleted monad
$ \mu (p) - \{p\}= \mu'(p) \cap \hy A \not= \emptyset$ iff $\mu (p) \cap \hy A=$ an infinite set.}\pars
Proof. These are rather easy to establish and, as usual, depend upon *-transform. (i) Let $p \in \real$ be an accumulation point for $A \subset \real.$ Then the formal sentence, which I'm sure you can obtain form the informal, 
$$\forall \r x ((\r x \in \real^+) \to \exists \r y(\r y \in A)\land \vert\r y - p\vert < \r x))$$
holds in ${\cal M}$; and, hence in $\hy {\cal M}.$ So, let $0 < \eps \in \mu (0).$ Then there exists some $ a \in \hy A$ such that $\vert a - p\vert < \eps;$ which implies that $a \in \mu (p).$ \pars
Conversely, assume that $\mu (p) \cap \hy A \not=\emptyset.$ Obviously $\mu (p) \subset \Hy (-w + p, p +w),\ \forall\, w \in \real^+.$ Hence, letting $b \in \mu (p) \cap \hy A$ and $w \in \real^+$ the sentence
$$\exists \r y ((\r y \in \hy A)\land \vert \r y - \hy p \vert < \hy w))$$ holds in $\hy {\cal M};$ and, hence, in ${\cal M}$ by reverse *-transform and the conclusion follows.\pars
\indent\indent (ii) The sufficiency follows since $p$ is an accumulation point. For the necessity, there exists some $w \in \real^+$ such that $(-w + p,p+w) \cap A = \{p\}.$ Hence, $\Hy (-w + p,p+w) \cap \hy A = \{\hy p\}= \{p\}$, under our notation simplification, and the result follows for $p \in \mu (p) \subset \Hy (-w + p,p+w).$ \pars
\indent\indent (iii) This follows from the observation about accumulation points, cluster points and isolated points and the fact that the only standard number in $\mu (p)$ is $p$. The second iff follows, for if otherwise there would be a ``smallest'' $w_1 \in \real^+$ such that $\Hy (-w_1 +p,p+w_1)' \cap \hy A \not= \emptyset.$ \qed\parm
{\bf Corollary 8.6.} {\it {\rm (i)} A point $p \in \real$ is an accumulation point for $A \subset \real$ iff there exists some $a \in \hy A$ such that $\st a = p.$\pars
\indent\indent {\rm (ii)} A point $p \in \real$ is a cluster point for $A\subset \real$ iff there exists an $a \in \hy A - \sig A$ such that $\st a = p.$} \parm
For $B \subset \hyperreal,$ define the {{\bf standard part of $B$}} as the set $\st B = \{ x\mid (x \in \real)\land \mu (x) \cap B \not=\emptyset\}.$ Of course, you can consider $\st B \subset \sig \real.$ Notice that for any $A \subset \real$, the standard part operator is defined, at the least, for all members of $\sig A.$ Indeed, our definitions and characterizations for these set-theoretic notions are only in terms of monads about standard points.  \parm
{\bf Theorem 8.7.} {\it Let $A \subset \real.$ Then $\st A = {\rm cl}A.$}\parm
{\bf Theorem 8.8.} {\it A point $ p \in \real$ is an interior point iff $\mu (p) \subset \hy A.$}\pars
Proof. I'm sure you can show that $\mu (p) = \bigcap\{\Hy (-w + p,p+w)\mid w \in\real^+\}.$ Hence, the necessity follows.\pars
For the sufficiency, assume that $p$ is not a member of the interior of $A$. Then for each $w \in \real^+,\ (-w + p,p+w) \cap (\real - A) \not= \emptyset.$ Thus, $p \in {\rm cl}(\real - A)$ and $\mu (p) \cap \Hy (\real -A) = \mu (p) \cap (\hyperreal - \hy A) \not=\emptyset$ implies that $\mu (p) \not\subset \hy A$ and the proof is complete. \qed\parm 
{\bf Definition 8.9.} Let $A\subset \real.$ Then the {{derived set $A'$}} for $A$ is the set of all cluster points. Notice that the derived set contains no isolated points. Example, let $A = (1,2) \cup \{3\}.$ Then $A' = [1,2].$ \parm
{\bf Theorem 8.10.} {\it For $A \subset \real,$ the set $A' = \st {\hy A - \sig A},$ \r (using the extended definition for $\St$ ).}\pars
Proof. Theorem 8.5 (ii). \qed \parm
{\bf Theorem 8.11.} {\it For $A \subset \real,$ the set of all isolated point is $A - \st {\hy A- \sig A}.$}\pars
Proof. An isolated point $p$ for $A$ is a member of $A,$ and such a $p$ is isolated iff $\mu (p)\cap \hy A = \{p\}$ iff $\mu'(p) \cap \hy A =\emptyset$ iff $p \in A - \st {\hy A - \sig A}$ (or in simplified notation) iff $p \in A - \st {\hy A - A}.$ \qed\parm
A set $A \subset \real$ is {{\bf closed}} $A = {\rm cl}A = \st {\hy A}.$ The set is {{\bf open}} iff $\mu (p) \subset \hy A,\ \forall\, p\in A.$ Please note that $\emptyset,\ \real$ are open and closed. (Actually, this is not the standard definition for an open nonempty set.  But, I leave it to you to show that this is equivalent to the statement that for each $p \in A$, there exists a $w_p \in \real^+$ such that $(-w + p_p,p +w_p) \subset A.$ Also $A$ is {{\bf perfect}} if it is closed and has no isolated points.\parm\vfil\eject 
{\bf Theorem 8.12} {\it A set $A \subset \real$ is perfect iff $A= A'$.}\pars
Proof. Please note that ${\rm cl}A = A \cup A'$ Hence, a set is closed iff $A' \subset A.$ For the necessity, $A' = \st {\hy A - \sig A} = \st {\hy A -\sig A} \cup A = \st {\hy A -\sig A} \cup \st {\sig A} = \st {(\hy A -\sig A) \cup \sig A} = \st {\sig A} = A.$\pars
The sufficiency is clear and this completes the proof. \qed\parm
  
Much of our interest will be restricted to the derived set. The reason for this is that for every $p \in A'$ there is a sequence $S\colon \nat \to p$ such that $p \notin S[\nat].$ Please consider the following remarkably short proof of the Bolzano-Weierstrass theorem.\parm 
{\bf Theorem 8.13.} {\it If bounded and infinite $A \subset \real$, then $\st{\hy A - \sig A} \not= \emptyset$ (i.e. $A$ has a cluster point).}\pars
Proof. Since $A$ is infinite then $\hy A - \sig A \not= \emptyset.$ Since $A$ is bounded that $\hy A \subset G(0).$ Thus, $\hy A - \sig A \subset G(0)$ implies that $\st {\hy A - \sig A} \not= \emptyset$ and this completes the proof. \qed \parm
One of the most important topological concepts used throughout analysis is the notion of ``compactness.'' Numerous equivalent definitions for this concept exist in  the literature. I select the most important for our purposes. Intuitively, compactness should mean ``closely packet'' or ``close together'' but it's different from the notion of density since density is usually a comparison between two different sets. Often this intuitive understanding for ``compactness'' is not achieved from the definition. I'll give a nonstandard definition that yields this intuitive notion and then show that it's equivalent to one of the usual definitions. \parm
{\bf Definition 8.14.} A set $A \subset \real$ is {{\bf compact}} iff for each $b \in \hy A$ there is some $p \in A$ such that $b \in \mu (p)$ (i.e. $b \approx p$) iff $\hy A \subset \bigcup \{\mu (p) \mid p \in A\}$ iff each  
 $b \in \hy A$ is {{\bf near-standard}} (meaning $\approx$ to a member $p \in A.$) The set $\bigcup \{\mu (p) \mid p \in A\}$ is often denoted by ns(A) (the set of all near-standard points).\parm
Our next, and what is a major, result requires what appears to be a rather long proof. I have not introduced the idea of the $\delta$-incomplete ultrafilter and concurrent relations. For the ultrafilters I am considering and due to real number property discussed in the next paragraph, the sufficiency part of the next theorem can be established in but a few lines using a concurrent relation. In general, this result holds for topological spaces, using a concurrent relation, if a special type of ultrafilter is used (Herrmann 1991). \pars
 A set $\cal G$ of nonempty open sets is said to {{\bf cover}} of (for) $A \subset \real$ iff  $A \subset \bigcup\{G\mid G \in {\cal G}\}.$ One standard definition for ``compactness'' says, that $A$ is compact iff for every open cover $\cal G$ there exists a finite subset (a subcover) ${\cal G}_f\subset \cal G$ such that 
${\cal G}_f$ covers $A$. A set $A$ is said to be {{\bf countable}} iff either $A$ or there exists a one-to-one correspondence from $\nat$ onto $A$. The {{\bf countably  compact sets}} are those that have this covering property but only for countable open covers. For the real numbers, due mainly to the fact that the rational numbers are dense in the reals and for every real $0< r <1$ there is a natural number $n$ such that
$r< 1/n< 1$, if nonempty $G$ is an open set, then there exists a rational number $w \in \real^+$ and a rational number $r \in \real$ such that $p \in I = (-w + r,r + w) \subset G.$ Thus, every open cover $\cal G$  of $A$ can be replaced by a countable open cover $\{I_i\}$ of such open intervals and such that $A \subset\bigcup\{I_i\} \subset \bigcup \{G\mid G \in {\cal G}\},$ where each member of $\cal G$ contains, at least, one member of $\{I_i\}.$ Hence,  replace the covering definition for compactness with  countable open covers by such a collection of open sets. \parm\vfil\eject
{\bf Theorem 8.15.} {\it Let nonempty set $A \subset \real.$ Then $\hy A \subset \bigcup\{\mu(p)\mid p \in A\}$ iff every countable open cover $\{I_i\}$ for $A$ has a finite subcover.}\pars
Proof. Assume that $A$ satisfies the countable covering definition  for compactness but that $\hy A \not\subset \bigcup\{\mu (p) \mid p \in A\}.$ 
There exists some $a \in \hy A$ such that $a\notin \mu (p)$ for any $p \in \sig \real.$ Consequently, there is some open interval $I(p)$ with rational end points about some rational number such that $a \notin \Hy I_p$ and $p \in I(p)$.  Let $\cal G$ be a set of all such intervals $I(p).$ Then $\cal G$ is a countable cover of $A$ and there should exists a finite subcover, say $\{I(p_1),\ldots,I_(p_n)\},$ such that $A \subset I(p_1)\cup \cdots \cup I(p_n).$ Consequently, $\hy A \subset \Hy I(p_1) \cup \cdots \cup \Hy I(p_n).$ Hence, we have the contradiction that $a \in \Hy I(p_i)$ for some $i = 1,\ldots,n.$\pars
For the sufficiency, just assume that there is a countable open cover $\cal G$ of $A$ which has no finite subcover. Our basic aim is to construct by induction from $\cal G$ another cover and do it in such a manner that a sequence of members of $A$ exists which, when viewed from the $\hyperreal$ and with respect to any free ultrafilter ${\cal U}$, the equivalence class containing this sequence is not near to any member of $A$. First, consider the nonempty countable set ${\cal G}' = \{C_i\mid i =1,2,\ldots\} = \{ x \cap A\mid (x \in {\cal G})\land (x \cap A \not= \emptyset)\}.$ Let $D_0 = C_1.$ Now, let $m_1$ be the smallest natural number greater than 1 such that $C_{m_1} \not\subset C_1.$ This unique number exists since $\{C_1\}$ cannot be a cover of $A$ for $C_1 \subset C$ for some $C\in {\cal G}.$  Assume that the $D_k$ have been defined. Let $m_{k+1}$ be the smallest natural number great than $m_k$ such that 
$$C_{m_{k+1}} \not\subset \bigcup \{D_i\mid i = 1,\ldots,k\}.$$
These unique natural numbers continue to exist since $A$ is not covered by any finite subset of sets in $\cal G$. Now define $D_{k+1} = C_{m_{k+1}}.$ The sets $D_n,\ \forall\, n \in \nat$ are defined by induction \pars
Let ${\cal G}_1 = \{D_n\mid n \in \nat\}.$ Since $\cal G$ is a countable cover of $A$, then ${\cal G}_1$ is a countable cover, although not generally an open cover. Further, ${\cal G}_1$ has no finite subcover. By definition $D_0 \not= \emptyset$ and 
$$ D_n - \bigcup\{D_k \mid k = 0,\ldots,n-1\} \not= \emptyset$$ for each positive $n \in \nat$ since $D_n \not\subset \bigcup\{D_k \mid k = 0,\ldots,n-1\}.$ Thus, define $p_0$ to be any point in $D_0$ and for each positive $n \in \nat$, define $p_n$ to be any point in $D_n - \bigcup\{D_k\mid k = 0,\ldots,n-1\}.$ (Did I use the Axiom of Choice or can this be considered an induction definition?) Thus, there is this sequence $P\colon \nat \to A$ such that $P(n) =p_n.$ If the natural number $m > n,$ then $p_m \notin D_i,\ i = 0,\ldots,n.$ Thus, if $p_m \in D_k$ for any $k = 0,\ldots,n$, then $m \leq n.$ This means that for each $n \in \nat$ the set of natural numbers 
$\{x\mid(x\in \nat)\land(P(x) \in D_n\}$ is finite. 
Hence, for each $n \in \nat,$ $\{x\mid(x\in \nat)\land(P(x) \notin D_n\}\in {\cal U}$ for any free ultrafilter $\cal U.$  This yields, in general, that  $[P]\notin \hy D_n$ for each $n \in \nat$ and $[P] \in \hy A.$ For $D_k \in {\cal G}_1,$ there exists some $c_k \in \cal G$ such that $D_k = A \cap c_k.$ Let ${\cal G}_2$ be the set of all such $c_k$. Since $[P] \notin \Hy D_n$ for  $n \in \nat$, then $[P] \notin \hy c_n.$   But, the set ${\cal G}_2$ is an open cover of $A$. Thus, for each $p \in A,$ there is some $c_k \in {\cal G}_2$ such that $\mu (p) \subset \hy c_k$. Consequently, $[P] \notin \bigcup \{\mu (p)\mid p \in A\}$ and the proof is complete. \qed \parm
Next, I present nonstandard proofs of a few additional characteristics for compactness, where trivially a finite $A \subset \real$ is compact.\parm
{\bf Theorem 8.16.} {\it A nonempty $A\subset \real$ is compact iff it is closed and bounded.}\pars
Proof. Assume that $A$ is compact. Since $\hy A \subset \bigcup\{\mu (p) \mid p \in A\}\subset G(0)$, the $A$ is bounded. Now let $\mu (q) \cap \hy A \not=\emptyset$. Then $\mu (q) \cap \mu (p)\not= \emptyset$ for some $p \in A.$ Hence, $q = p.$ Thus, $A = \st {\hy A}$ and $A$ is closed. \pars
Conversely, let $A$ be bounded. Then $\hy A \subset \bigcup\{\mu (x)\mid x \in \real\} = G(0).$ Also, $A \not= \real$. Let any $q \in \real -A.$ Since $A$ is closed, then $\mu (q) \cap \hy A = \emptyset.$ Thus $\hy A \subset \bigcup\{\mu (p)\mid p \in A\}$ and this completes the proof. \qed\parm
{\bf Theorem 8.17.} {\it Let infinite $A \subset \real.$ Then $A$ is compact iff each infinite $B \subset A$ has a cluster point in $A$.}\pars
Proof. Since $A$ is compact, then $A$ is bounded and, hence, $B$ is bounded.    
Thus, by 8.13, $B$ has a cluster point $p$. But, since $A$ is closed, then $p \in A.$\pars
For the sufficiency, assume that $A$ is not compact. Then either $A$ is not bounded or $A$ is not closed. Assume that $\hy A\not\subset G(0).$ Let $r = 1.$ Consider the case, that $A$ is not bounded above. Then there's some $p_1 \in A$ such that $p_1 > 1.$ Let $r = p_1 +1.$ Then there exists some $p_2 \in A$ such that $p_2 > p_1 +1.$ Assume that we have defined $p_k$. Then there is some $p_{k+1}$ such that $p_{k+1} > p_k +1 > p_{k-1} +1 > \cdots > 1.$ Let $p_0 =1.$ Thus, there is a sequence $P \colon \nat \to \real$ such that $\lim_n p_n = +\infty.$ Hence, this sequence has no accumulation point in $A$, which in this case is equivalent to not having a cluster point for the
infinite $P[\nat] \subset A.$ The case where $A$ is not bounded below follows in like manner. \pars
Now suppose that $A$ is not closed. Then there exists some $q \in A' -A.$ Hence, there is an infinite sequence of distinct members of $A$ that converges to $q$. Again, $q$ would be a cluster point for $A$. This completes the proof. \qed \parm \vfil\eject 

\centerline{\bf 9. BASIC CONTINUOUS FUNCTION CONCEPTS}
\parm
{\bf For all that follows in this chapter, $D$ will denote the domain for the real valued function $f$.} Recall that the notation $f\colon D \to \real$ means that $f$ is a real valued function defined on $D$. Of course, in this case, $f$ is also defined on any nonempty subset of $D$. First, let's consider the idea of the limit of $f$ as $x \to s$ or $\lim_{x \to s} f(x)$ or, in abbreviated notation, $\lim_s f(x)$ where I use $s$ so as not to confuse this with the more general notation for the specific case where we look only at sequences and use $n$ or $m$ below the $\lim$ symbol.\pars
Recall that for $f \colon D \to \real,$ $\lim_s f(x) = L$ iff for every $r \in \real^+$, there exists some $w \in \real^+$ such that, whenever $x\in D$ and $0 <\vert x-s\vert < w,$ then $\vert f(x) -L\vert < r.$ Clearly, {\bf $s$ must be a cluster point for $D$. That is $s\in D'$} for this notion to have a significant unique meaning. This is one of the first definitions  that appears in a calculus book and that often gives students some difficulty in its application. But, as will be seen, the nonstandard characteristics, especially (i), for this limit concept are much easier to state and yield the actual intuitive idea. Recall that for each $p \in \real$ $\mu'(p) = \mu (p) -\{p\}$ is the deleted monad about $p$ and if $g\colon B \to \hyperreal, \ A\subset B$, then $g[A] = \{g(x)\mid x \in A\}.$\parm
{\bf Theorem 9.1.} {\it Let $f \colon D \to \real$. Then $\lim_s f(x) =L$ iff\pars
\indent\indent {\rm (i)} $\hy f[\mu'(s)\cap \hy D] \subset \mu (L)$ iff\pars
\indent\indent {\rm (ii)} for each $q \in \mu'(s)\cap \hy D, \st{\hy f(q)} = L$ iff\pars
\indent\indent {\rm (iii)} for each nonzero $\eps\in \mu (0)$ such that $s + \eps \in \hy D,$ then $\hy f(s + \eps) -L \in \mu (0)$ iff \pars
\indent\indent {\rm (iv)} for each $\eps \in \mu (0)^+$ and $x \in  \hy D$ such that $0 < \vert x - s\vert < \eps,$ then $\hy f(x) - L \in \mu (0).$} \pars
Proof. (i) For the necessity, let $\lim_s f(x) = L$ and $r \in \real^+.$ Then there exists some $w \in \real^+$ such that the following sentence
$$\forall \r x((\r x \in D) \land (0 < \vert \r x- s \vert < w) \to (\vert f(\r x) - L\vert < r))$$
holds in ${\cal M};$ and , hence, in $\hy {\cal M}.$ In particular, for each $p \in \mu'(s) \cap \hy D$, $\vert \hy f(p) - L \vert < r.$ Since $r$ is an arbitrary positive real number, and we have that $0<\vert p - s \vert < w$ for all $w \in \real^+,$ it follows that for each $p \in \mu'(s) \cap \hy D,\ \vert \hy f(p) - L \vert \in \mu (0)$ or that $\hy f(p) \in \mu (L).$\pars
For the sufficiency, assume that $r \in \real^+.$ There exists a $q \in \mu' (s) \cap \hy D$ since $s$ is a cluster point of $D$. Thus, $q \not= s$ and, hence, there is some $\eps \in \mu'(0)$ such that $q = s +\eps.$ Consequently, $0 < \vert q-s\vert = \vert \eps\vert \in \mu (0).$ If $p\in \hy D$ such that $0 < \vert p -s\vert < \vert \eps \vert,$ then $p \in \mu' (s) \cap \hy D$ implies that $\hy f(p)\in  \mu (s).$ Consequently, the sentence
$$\exists \r x ((\r x \in \real^+) \land \forall \r y((\r y \in D)\land (0< \vert \r y-s\vert <\r x)\to (\vert f(\r y) - L\vert < r))$$ holds in $\cal M$ by reverse *-transform and this first ``iff'' is established.\pars
All but the last ``iff'' are immediately equivalent to this first one. The necessity of ``iff'' (iv) is clear. The sufficiency of (iv) follows from the above sentence for the sufficiency of (i) and this completes the proof. \qed\parm
{\bf Corollary 9.2.} {\it If $\lim_s f(x) = L,$ then $L$ is unique.}\parm
{\bf Corollary 9.3.} {\it If $T \subset D,\ s \in T'$ and $\lim_s f(x) =L$ with respect to $D$, then $\lim_s f(x) = L$ with respect to $T$.}\parm
{\bf Corollary 9.4.} {\it If $\lim_s f(x) = L,$ then there exists a nonempty  open set
$G$ such that $s \in G$ and $\hy f[(\Hy G - \{s\})\cap \hy D] \subset G(0).$ Note that $\{s\}$ is not an open set.} \parm
{\bf Theorem 9.5.} {\it $\lim_s f(x) = L$ iff there exists a sequence $S$ such that for each $n \in \nat, \ S_n \not= s,\ S_n \in D;\ S_n \to s$ and $\lim_nf(S_n) = L.$}\pars
Proof. Suppose that $\lim_s f(x) = L$ and that $S\colon \nat \to D,\ S_n \to s$ and that for each $n \in \nat,\ S_n\not=s.$ Then for each $\Lambda \in \nat_\infty,\ \hy S(\Lambda) \in \mu (s),\ \hy S(\Lambda)\not= s$ and $\hy S(\Lambda) \in \hy D$ implies that $\hy S(\Lambda) \in \mu'(s) \cap \hy D.$ Hence, $\hy f(\hy S(\Lambda)) = \Hy (fS) (\Lambda) \in \mu (L)$ and the necessity follows.\pars
For the sufficiency, assume that $\lim_s f(x) \not\to L.$ Then there exists some $r \in \real^+$ such that for each $w \in \real^+$ whenever  $x \in D$ and $0 < \vert x- s\vert < w,$ it follows that $\vert f(x) - L \vert \geq r.$ Since $\bigcap\{\Hy (-w+s,s)\mid w\in\real^+\} \cap \hy D \not= \emptyset,$ then for each $w = 1/n,\ 0\not= n \in \nat,$ there exists a sequence such that $S_n \not= s,\ S_n \in D,$ and $0 < \vert S_n -s\vert < 1/n$ and $\vert f(S_n) -L \vert \geq r.$ Consequently, $S_n \to s$, but $f(S_n) \not\to L$ and the proof is complete. \qed \parm
\noindent Of course, it's this ``sequence'' theorem that gives the major intuitive characteristic for such limits. \pars
Modifying the definition for $\lim_s f(x) = L$ yields the {{\bf one-sided}} limits. Recall that the modifications are $f(s\pm) = L$ iff for each $r \in \real^+$ there exists a $w \in \real^+$ such that, whenever $\cases{0<x-s <w&\cr 0< s-x < w&\cr},$ then $\vert f(x) - L \vert < r.$ For these limits, the monads need to be modified in the obvious manner. For each $p\in \real$, let $\mu (p)^+ = \{x\mid (x > p) \land (x \in \mu (p))\}=\{x\mid (x > p) \land (x\approx p)\}= \bigcap\{\Hy (p, p+w)\mid w \in \real^+\},\ \mu (p)^- = \{ x \mid (x< p)\land(x \in \mu (p))\}=\{ x \mid (x< p)\land(x \approx p)\}= \bigcap\{\Hy (-w + p,p)\mid w \in \real^+\}.$ Using these {{positive or negative monads}} our previous theorems and corollaries all hold with the appropriate modifications.\parm 
 {\bf Theorem 9.6.} {\it Let $f \colon D \to \real$. Then $f(s\pm) =L$ iff\pars
\indent\indent {\rm (i)} $\hy f[\mu'(s)^\pm\cap \hy D] \subset \mu (L)$ iff\pars
\indent\indent {\rm (ii)} for each $q \in \mu'(s)^\pm\cap \hy D, \st{\hy f(q)} = L$ iff\pars
\indent\indent {\rm (iii)} for each $\eps\in \mu (0)^\pm$ such that $s + \eps \in \hy D,$ then $\hy f(s + \eps) -L \in \mu (0)$ iff \pars
\indent\indent {\rm (iv)} for each $\eps \in \mu (0)^+$ and $x \in  \hy D$ such that $\cases{0<x-s<\eps&\cr 0<s-x<\eps&\cr},$ then $\hy f(x) - L \in \mu (0).$} \parm
{\bf Corollary 9.7.} {\it If $f(s\pm) = L,$ then $L$ is unique.}\parm
{\bf Corollary 9.8.} {\it If $T \subset D,\ s \in T'$ and $f(s\pm) =L$ with respect to $D$, then $f(s\pm) = L$ with respect to $T$.}\parm
{\bf Corollary 9.9.} {\it If $f(s\pm) = L,$ then there exists a nonempty open 
interval $\cases{I^+=(s,r)\cr I^-=(r,s)\cr}$ such that $\hy f[(\Hy I^\pm)\cap \hy D] \subset G(0).$ Note that $\{s\}$ is not an open set.} \parm
The following is the appropriate modification for Theorem 9.5\parm
{\bf Theorem 9.10.} {\it Let $f \colon D \to \real.$ Then $f(s+) = L$ \r [resp. $f(s-)$\r ] iff there is a  sequence $S$ such that for each $n \in \nat,\ S_n \in D, \ S_n >s,$ \r [resp. $S_n<s$\r ], $S_n \to s$ and $\lim_n f(s_n) = L.$}\pars 
Proof. I prove this for $f(s-)$ since $f(s+)$ is done in like manner. Let $f(s-) = L$ and $S_n \to s,\ \forall \, n \in \nat, S_n \not= s,\ S_n < s,\ S_n \in D$. Then $\forall\, \Lambda \in \nat_\infty,\ \hy S(\Lambda) \in \mu (s)^-$
and $S(\Lambda) <s.$ Thus $\hy f(\hy S(\Lambda))= \Hy(fS)(\Lambda)\in  \mu (L)$ and the necessity follows.\pars
For the sufficiency, the method is similar to that for Theorem 9.5. Assume that $f(s-) \not\to L.$ Then there exits some $r \in \real^+$ such that 
$\forall\, w \in \real^+$ whenever $0 < s -x <w,\ x \in D$, then $\vert f(x) - L \vert \geq r.$ Since $\mu (s)^- = \bigcap\{\Hy (-w +s,s)\mid w \in \real^+\}$, by *-transform, for $0\in \nat$, there is some $s_0 \in (-1+ s,s)\cap D$ and $S_0 < s.$ Assume that for $k \in \nat,\ k \geq 1,$ there are $S_k \in (-1/(k+1) + s,s)\cap D,$ and $S_k < s.$ Now consider $k+1.$ Then since $(-1/(k+2) +s,s)\cap D\not=\emptyset$, there is some $S_{k+1} \in (-1/(k+2) +s,s)\cap D$ and $S_{k+1} < s.$ Thus yields a sequence $S \colon \nat \to D$ such that $\forall \, n \in \nat,\ S_n \in D,\ S_n \to s$ and $S_n < s$, but $\vert f(S_n) -L\vert \geq r$ and the proof is complete. \qed\parm
{\bf Theorem 9.11.} {\it Let $f\colon D \to \real$ and $s \in {\rm int}(D)$ (the set of all interior points). Then $\lim_s f(x) = L$ iff $f(s\pm) = L.$ }\pars
Proof. $\mu' (s) = \mu (s)^+ \cup \mu (s)^-.$ \parm
{\bf Example 9.12.} In the usual calculus text, it's established that 
  $\lim_0{{\sin x}\over {x}} = 1,$ by means of a geometric proof. Although, I won't mention any apparent geometric facts in the following nonstandard proof, it might be necessary to use the geometric definitions to establish the facts I do use. \pars
For each $r \in \real^+$ such that $0 < r < \pi/2,$ since $\sin (r) < r < \tan (r),$ $\cos (r) < {{\sin (r)}\over{r}} < 1.$ Thus, for each $\eps \in \mu (0)^+,$ by *-transform,
$$\hy \cos (\eps) < {{\hy \sin (\eps)}\over{\eps}} < 1.$$
But, since $\vert \sin (r)\vert \leq \vert r \vert,\forall \, r \in \real$, then
$\vert \hy \sin(\eps)\vert \leq \eps$ implies that $\hy \sin (\eps) \in \mu (0).$ This yields
$$1 - \Hy (\cos (\eps))^2 = \Hy (\sin (\eps))^2 \in \mu (0)$$
which implies that $\hy \cos (\eps) \in \mu (1).$ Consequently, 
$1 \leq \st {\hy \cos (\eps)} \leq \st { {{\hy \sin (\eps)}\over{\eps}}} \leq 1$ for each $\eps \in \mu (0)^+.$ Thus, ${{\sin (0+)}\over{0^+}} =1.$ To show that this last equation holds for $\eps \in \mu (0)^-$, simply notice that ${{\hy \sin (-\eps))}\over{-\eps}}= {{\hy \sin (\eps))}\over{\eps}}\,  \forall\, \eps \in \mu'(0).$ Hence, the result follows.\parm
All of the usual limit and one-sided limit algebra for such functions follow from the properties of the standard part operator. Now let's establish the {{Cauchy Criterion for functions.}}\parm
{\bf Theorem 9.13.} {\it \r (Cauchy Criterion.\r ) Let $f\colon D \to \real.$ Then $\lim_s f(x) = L$ iff for each pair $p,q \in \mu'(s) \cap \hy D, \ \hy f(p) - \hy f(q) \in \mu (0).$}\pars
Proof. The necessity follows from Theorem 9.1. \pars
For the sufficiency, assume that there does not exists $w \in \real^+$ such that $f$ is bounded on $(-w +s,s+w)'\cap D.$ Hence, for  $r = 1,$ for each $w \in \real^+$ there are, at least, two distinct $x_1,x_2  \in (-w+s,s+w)' \cap D,\ x_1,x_2\not= s$ and $\vert f(x_1)-f(x_2) \vert \geq 1.$ Consequently, the sentence
$$\forall \r x((\r x \in \real^+ \to \exists \r y\exists \r z ((\r y \in D)\land (\r z \in D)\land (0<\vert s-\r y\vert < \r x)\land$$ $$ (0<\vert s-\r z \vert < \r x)\land (\vert f(\r y) -\vert f(\r z))\geq 1)))$$ holds in $\hy {\cal M}$ by *-transform. So, let $\eps \in \mu (0)^+.$ Then there exists distinct $p,q$ such that $0< \vert s-p\vert <\eps$ and $0 < \vert s-q\vert <\eps$ and $\vert \hy f(p) - \hy f(q) \vert \geq 1.$ But, this contradicts the requirement that $\hy f(p) - \hy f(q) \in \mu (0).$ Thus, there is some $w \in \real^+$ such that $f$ is bounded on $(-w +s,s+w)'\cap D.$ Consequently, $\hy f[\mu'(s)\cap \hy D] \subset G(0).$ Letting $q \in \mu'(s)\cap\hy D$, then for each $p \in \mu'(s) \cap \hy D,\ \hy f(p) \in \mu (\st {\hy f(q)}$ and the result follows where $\st {\hy f(q)} = L.$\parm
{\bf Corollary 9.14.} {\it Let $f\colon D \to \real$. Then $f(s\pm) = L$ iff for each pair $p,q \in \mu'(s)^\pm \cap \hy D,\ \hy f(p) - \hy f(q) \in \mu (0).$} \parm
{\bf Theorem 9.15.} {\it Let $f\colon (a,b) \to \real$ and $a< c< d <b$. 
If $f$ is increasing \r [resp. decreasing\r ], then 
$$f(c-) =\sup \{f(x)\mid a < x < c\} \leq f(c) \leq f(c+) = \inf\{f(x)\mid c < x < b\}.$$
$$\r [resp. \ f(c+) =\inf \{f(x)\mid a < x < c\} \leq f(c) \leq f(c-) = \sup\{f(x)\mid c < x < b\}.\r ]$$
Further, $f(c+) \leq f(d-)$ \r [resp. $f(d-) \leq f(c+)$\r ].}\pars 
Proof. I show this only for an increasing function $f$. Clearly, $\sup \{f(x)\mid a< x < c\} = L,\ L \leq f(c).$ For any real number $r < L,$ there exists some $p \in (a,c)$ such that $r < f(p) \leq L.$ Thus, let $\eps \in \mu (0)^-.$ Then $ a< c+\eps < c$ implies that $r < \hy f(c+ \eps) \leq L$ since $\hy f$ is increasing on $\Hy (a,b).$ Therefore, $r < \st {\hy f(c+\eps)} \leq L.$ Since $r < L$ is arbitrary, this implies that $\st {\hy f(c+\eps)} =L$ for each such $\eps,$ and this first part follows from Theorem 9.6. The $\inf$ case, follows in like manner. \pars
Now if $c < d$, then $c + \eps < d + \gamma$ for each $\eps \in \mu (0)^+$ and each $\gamma \in \mu (0)^-.$ Consequently, $\st {\hy f(c + \eps)}= f(c+) \leq f(d-) = \st {\hy f(d + \gamma)}$ and the proof is complete. \qed \parm 
I guess I should mention the other ordinary limit of a function notion used when $D$ is not bounded above or below, the $\infty$. Recall that if $D$ is {{\bf not bounded above}}, then $\lim_\infty f(x) = L$ iff for each $r \in \real^+$ there exists some $w \in \real^+$ such that for each $p \in D$ such that $p >w,\ \vert f(p) -L\vert < w.$ For $D$ that is not bound below, this limit notion is defined in the obvious manner.\parm
{\bf Theorem 9.16.} {\it Suppose that  $f\colon D \in \real$ is not bounded above \r [resp. below\r ]. Then $\lim_\infty f(x) = L$ iff $\hy f(p) \in \mu (L)$ for each $p \in \real_\infty^+ \cap \hy D$ \r [resp. $\hy f(p) \in \mu (L)$ for each $p \in \real_\infty^- \cap \hy D$\r ].}\parm
Proof. Left to the reader.\parm 
{\bf Theorem 9.17.} {\it Suppose that $f\colon \to D \in \real$ is not bounded above \r [resp. below\r ]. Then $\lim_\infty f(x) = L$ iff for each pair $p,q \in \real_\infty^+ \cap \hy D$ \r [resp. $\real_\infty^- \cap \hy D$\r ], $\hy f(p) - \hy f(q) \in \mu (0).$}\pars
Proof. Left to the reader. \qed\parm
Our major interest is to investigate properties of continuous real valued functions defined on $D.$ Since for $f\colon D \to \real$ to be continuous at $s$, all one needs is that $\lim_s f(x) = f(s)$ and, hence, we need $s \in D.$ \parm 
{\bf Theorem 9.18.} {\it Let $s \in {\rm int}(D)$. Then function $f\colon D \to \real$ is continuous at $s$ iff $\lim_s f(x) = f(s) = f(s+)= f(s-).$}\pars
Proof. Note that $\mu (s) = \mu (s)^+ \cup \{s\} \cup \mu (s)^-.$\parm
Each of the previous theorems on the left and right-hand limits, when slightly modified, hold for continuous functions. Also, each monadic characteristic for continuity holds for isolated points. Thus, $s$ need not be a cluster point. The changes are made by replacing the deleted monads with the complete monad and such statements as $0< \vert x -s \vert$ by $\vert x-s\vert$ and the like. The must used result is that $f$ is continuous at $p \in D$ iff $\hy f[\mu (p) \cap \hy D] \subset \mu (f(p)).$ Now let's apply these results
to obtain three highly significance continuous function properties.\pars
{\bf Theorem 9.19.} {\it Let continuous $f\colon D \to \real$ and let $D$ be compact. Then the range, $f[D],$ is compact.}\pars
Proof. Since $D$ is compact, then $\hy D \subset \bigcup \{\mu (p)\mid p \in D\}.$ But, using a property that holds for any function, it follows that 
$$\hy f[\hy D]\subset \bigcup\{\hy f[\mu (p)\cap \hy D]\mid p \in D\}\subset    
\bigcup\{\mu (f(p))\mid f(p) \in f[D]\}$$
and the result follows. \qed \parm
{\bf Theorem 9.20.} {\it \r (Extreme Value Theorem.\r ) Let continuous $f \colon D\to \real$ and $D$ be compact. Then there exists $p_m,p_M \in D$ such that for each $p \in D,$ $f(p_m) \leq f(p) \leq f(p_M).$}\pars
Proof. Since $f[D]$ is compact, then it is closed and bounded. Thus, from boundedness $\sup\{f(x)\mid x \in D\} = p_M$ and $\inf\{f(x)\mid x \in D\} = p_m.$ Since $f[D]$ is closed that $p_m,\ p_M \in D$ and the result follows.\parm
I mention that all such standard theorems can be extended to ``nonstandard statements'' by *-transform. To establish the intermediate value theorem the notion of connectedness is often introduced. But, rather than do this, I'll give a nonstandard proof where connectedness is not mentioned. \parm
{\bf Theorem 9.21.} {\it Let continuous $f\colon [a,b] \to \real$. The for each $d$ such that $f(a)\leq d\leq f(b)$ \r [resp. $f(b) \leq d \leq f(a)$\r ], there is some $c\in [a,b]$ such that $f(c) = d.$}\pars
Proof. The result is immediate if $a=b.$  So, assume that $a< b$ and consider the case where that $f(a)\leq d \leq f(b).$ Let nonzero $n \in \nat$ and $h = (b-a)/n.$ Then we have a finite partition of $[a,b]$ $\{a,a+h, a+2h,\ldots, a+nh =b\}.$ Thus, there exists some $m \in  \nat$ such that $m<n$ and (i) $f(a) \leq f(a +mh) \leq d\leq  f(a + (m+1)h)\leq f(b)$ or (ii) $f(a) \leq f(a + (m+1)h)\leq d\leq f(a +mh)\leq f(b).$ Assume (i). By *-transform, if $\Lambda \in \nat_\infty,$ then $(b-a)/\Lambda \in \mu (0).$ There exists some $m_1 \in \hypernat$ such that $m_1 <\Lambda$ and $f(a)\leq \hy f(a + m_1h) \leq d \leq \hy f(a + (m_1+1)h)\leq f(b).$ (Note the use of simplified notation for such things as $f(a)$, where technically this should be written as $\sig f(\hy a).$) Since $a < a +m_1h \leq b,$ then there is a real $c = \st {a+ m_1h}$ and $a \leq c \leq b$ From the continuity of $f,$ $f(c) = f(\st {a+m_1h}) = \st {\hy f(a + m_1h)}\leq d$ for $a+ m_1h \in \mu (c)\cap \Hy [a,b].$ However, $a + m_1h +h = a + (m_1 +1)h\in \mu (c)$ implies that $f(c) = f(\st {a+ (m_1)h}) =\st {\hy f(a+ m_1 +1)h)} \geq d.$ Therefore, $f(c) = d.$ The other cases follow in a similar manner and the proof is complete. \qed \parm
 The results that the sum and product function and similar processes defined for continuous functions yield continuous functions follows from the properties of the standard part operator. Our last result in this chapter is a nonstandard proof of the composition properties for continuous functions.\parm
{\bf Theorem 9.22.} {\it Let continuous $f\colon D \to \real$ and continuous $g\colon T \to R$ be such that $f[D] \subset T.$ Then the composition $gf\colon D\to \real$ is continuous.}\pars
Proof. Let $p \in D$. Then $\hy f[\mu (p)\cap \hy D] \subset \mu (f(p))\cap \Hy (f[D]) \subset \Hy (f[D]) \subset \Hy T$ imply that $\hy g[\hy f[\mu (p) \cap \hy D]] \subset \hy g[\mu (f(p))\cap \Hy T] \subset \mu (g(f(p))$ and the result follows. \qed \parm
\vfil\eject
\centerline{\bf 10. SLIGHTLY ADVANCED CONTINUOUS FUNCTION CONCEPTS}
\parm
{\bf Unless otherwise specified, for all the follows in this chapter, $D$ will denote the domain for the real valued function $f$.} Here is a result, you may never have seen before, that also implies the intermediate value theorem. The original standard proof and result is due to Bolzano.\parm
{\bf Theorem 10.1.} {\it For continuous $f\colon [a,b] \to \real,$ if $f(a)f(b)<0,$ then there exists some $c \in (a,b)$ such that $f(c) = 0.$}\pars
Proof. First note that the hypotheses require that $a \not= b.$ Assume that $f(c) \not= 0$ for each $c \in (a,b).$ Since $f(a) \not= 0$ and $f(b) \not= 0,$ then $f(c) \not= 0\ \forall \, c \in [a,b].$ I now show that for each nonzero $m \in \nat$ (i.e. $m \in \nat'$) that there exist real numbers $s_m,t_m$ such that $t_m -s_m = (b-a)/m$ and 
$$a\leq s_m < t_m <b,\ {\rm and}\ {{f(t_m)}\over{f(s_m)}} < 0.$$
For $m \in \nat',$ consider the function
$$g(x) = {{f(x +(b-a)/m)}\over{f(x)}};\ a\leq x \leq a + {{m-1}\over{m}}(b-a).$$
Consequently, the product $\Pi_0^{m-1} g(a + k(b-a)/m) = f(b)/f(a) <0.$ Thus, there is some $k\in \nat,\ 0\leq k \leq m-1$ such that $g(a +k(b-a)/m) <0.$ Hence, $f(a+(k+1)(b-a)/m)/f(a + k(b-a)/m) <0.$ Let $s_m = a + k(b-a)/m$ and $t_m = a + (k+1)(b-a)/m$ and the conditions required hold for $s_m,\ t_m$. Thus, by *-transform, if $\Lambda \in \nat_\infty,$ there exists $p,q \in \hyperreal$ such that $q - p = (b-a)/\Lambda,\ a \leq p < q\leq b$ and $\hy f(q)/\hy f(p) < 0.$ Since $q-p \in \mu (0),$ then $p \in \mu (\st {q})$ Consequently, using the result that $\st q \leq b$ and the continuity of $f$, $\hy f(p) \in \mu (f (\st q)).$ Therefore, $\st {\hy f(p)} = f(\st p)= f(\st q)$ implies that $f(\st p) /f(\st q) = 1 \leq 0.$ This contradiction yields the result. \qed\parm 
To obtain the immediate value theorem from Theorem 10.1, just consider for the function $f(x) -d,$ if $f(a) \leq d \leq f(b),$ or $d-f(x)$ if $f(b)\leq d\leq f(b)$ for the non-trivial cases $f(a) \not= d$ and $f(b)\not= d.$ A major result characterizes continuity on the entire set $D$ in terms of open sets. The proof is a little long due to the simplified structure I'm using. Let $\cal G$ be a nonempty collection of open subsets of $\real$. Then since for each $p \in \bigcup\{G\mid G \in {\cal G}\},$ $\mu (p) \subset \Hy G$ for some $G \in \cal G$, then $\mu (p) \subset \bigcup \{\Hy G\mid G \in {\cal G}\} \subset \Hy (\bigcup \{G\mid G \in {\cal G}\})$ implies that the {arbitrary union of a collection of open sets is an open set.}\parm 
{\bf Theorem 10.2.} {\it Let $f\colon D\to \real$. Then $f$ is continuous on $D$ iff for each open set $G \subset\real,$ $f^{-1}[G]$ is open in $D.$}\pars
Proof. Note that a set $G_1$ is ``open'' in $D$ iff there exists an open $G_2 \subset\real$ such that $G_1 = G_2 \cap D.$ Assume that $f$ is continuous on $D$. Let $G$ be an open set in $\real$. If $G = \emptyset$, then $f^{-1}[G] = \{p \mid (p \in D)\land (f(p) \in G\} = \emptyset,$ which is open in $D$. The same result would hold if $G \cap f[D] = \emptyset.$ Hence, assume that $G \cap f[D] \not= \emptyset.$
Then let $f(p) \in G \cap f[D].$ Since $G$ is open, then $\mu (f(p))\subset \Hy G$ and, by continuity, $\hy f[\mu (p) \cap \hy D] \subset \mu (f(p)) \cap \hy D \subset \Hy G.$ Since $\mu (p)$ is the intersection of all the 
intervals $\Hy (-r +p, r +p),\ r\in \real^+$, then there exists some $(-r +p,p+r),\ r \in \real^+$, such that $p \in (-r +p,p+r)\cap D \subset f^{-1}[D].$ Since the arbitrary union of open sets is an open set, then using one of these open intervals for each $p \in D$, one gets an open set $G_0 \subset \real$ such that $G_0 \cap D = f^{-1}[D].$ \pars
For the sufficiency, I'll use inverse image, $f^{-1},$ set-algebra. 
Consider $\mu (f(p)) \cap \Hy (f[D]) = \bigcap \{\Hy (-r +f(p),f(p) +r\mid r \in \real^+\} \cap \Hy (f[D])$ implies that $\hy f^{-1}[\mu (f(p))\cap \hy f[\hy D]] = \hy f^{-1}[\mu (f(p))] \cap \hy D = \bigcap\{\Hy(f^{-1}[(-r +f(p),f(p) + r)]\cap D)\mid r \in \real^+\}.$ Now by the hypothesis, each $f^{-1}[(-r +f(p),f(p) +r))]\cap D$ is open in $D$. Hence, for $p \in f^{-1}[(-r +f(p),f(p) +r)]\cap D$, there is some $s\in \real^+$ such that $p \in (-s +p,p +s) \cap D \subset f^{-1}[(-r +f(p),f(p) +r))]\cap D$. However, 
$\mu (p) \cap \hy D \subset \Hy (-s + p,p+s)\cap \hy D$ implies that 
$\mu (p) \cap \hy D \subset \hy f^{-1}[\mu (f(p))] \cap \hy D.$ Therefore,
$$\hy f[\mu (p)\cap \hy D] \subset \hy f\hy f^{-1}[\mu (f(p))] \cap \hy f[\hy  D] \subset \mu (f(p))\cap \hy f[\hy D] \subset \mu (f(p)).$$ The proof is complete. \qed \parm  
Prior to considering the notion of uniform continuity, here is a nonstandard proof of a rather interesting result. A real valued function is {{\bf additive}} if for each $p,q \in \real,$ $f(p + q) = f(a) + f(q).$ Recall that I'm using simplified notation in that rather than write a  statement such as $\hy x \in \sig \real,$ this is often written as $x\in \real.$\parm
{\bf Theorem10.3.} {\it Let $f\colon\real \to \real$ be additive. If $f$ is bounded on some non-empty interval $I$, then $f(x) = xf(1)$ for each $x \in \real$ and is a continuous function on $\real.$}\pars
Proof. Clearly, $ \hy f\colon \hyperreal \to \hyperreal$ is additive. Additivity implies that for any rational $r$ and any $x \in \real$, $f(rx) = rf(x).$ Hence,$\hy f(rx) = r\hy f(x)$ for each $x \in \hyperreal$ and *-rational $r \in \hyperreal$ (i.e. $r \in \Hy Q)$. Let $p$ be in the interior of $I$ (i.e $\mu (p) \subset \Hy I).$ Then from boundedness, $\vert \hy f[\mu (p)] \vert \leq M \in \sig R.$ Consequently, for each $\eps \in \mu (0), \ \vert \hy f(p + \eps)\vert = \vert \hy f(p) + \hy f(\eps)\vert \leq M$ implies that $\vert \hy f(\eps)\vert \leq M + \vert f(p) \vert.$ Now for each $n \in \nat,\ n\eps \in \mu (0)$ implies, by additivity, that $\vert \hy f(n\eps)\vert = n\vert \hy f(\eps)\vert \leq M + \vert f(p)\vert.$ Therefore, for $n \in \nat',\ \vert \hy f(\eps)\vert \leq (M + \vert f(p)\vert)/n.$ This yields that $\hy f(\eps) \in \mu (0),\ \forall\, \eps \in \mu (0).$ From the density of the rational numbers in $\real$, for any $r \in \real^+$ and any $x \in \real$, there is some $q \in Q$ such that $\vert x-q\vert <r.$ By *-transform, we have that for $\eps \in \mu (0)^+$  and $x \in \real,$ there is $q \in \Hy Q$ such that $\vert x -q \vert < \eps.$ Hence, $x - q \in \mu (0)$ implies that there is some $\gamma \in \mu (0)$ such that $x= q + \gamma.$ Therefore 
$$f(x) = \hy f(q + \gamma) = \hy f(q) + \hy f(\gamma) = \hy f(q\cdot 1)+ \hy f(\gamma) = q(f(1)) + \hy f(\gamma).$$
Thus, $q(f(1))\in \mu (f(x)).$ Finally, $f(x) =\st {q(f(1))} = (\st q)\st{f(1)}= xf(1)$ and, obviously, $f$ is continuous. \qed\parm
{\bf Theorem 10.4.} {\it If $r \in \real,$ then there exists a hyperrational $r \in \Hy Q$ and some $\eps \in \mu (0)$ such that $r = q + \eps.$}\parm
You might try showing from Theorem 10.4 that if $f\colon \real \to \real$ and $\forall\, x,y \in \real\ f(x+y) = f(x)f(y),\ \hy f[\Hy Q] \subset G(0), \lim_0 f(x) = 0,$ then $f(x) = 0,\ \forall\, x \in \real.$\pars
Recall that $f\colon D \to \real$ is {{\bf uniformly continuous}} on $D$ if for each $r \in \real^+$ there exists some $w \in \real^+$ such that whenever $x,y \in D$ and $\vert x-y\vert <w$, then $\vert f(x) -f(y)\vert < r.$ Of course, this concept is highly significant in series work and integration theory. The follow characteristic follows in the usual manner, where the big difference between this and Corollary 9.14, extended to continuity, is that the points are not restricted to a particular $\mu (s).$ \parm
{\bf Theorem 10.5.} {\it The function $f\colon D \to \real$ is uniformly continuous on $D$ iff for each $p,q \in \hy D$ such that $q -p \in \mu (0),\ \hy f(p) - \hy f(q) \in \mu (0).$}\parm
{\bf Theorem 10.6.} {\it Let $f \colon D \to \real$ be continuous on compact $D$. Then $f$ is uniformly continuous.}\pars
Proof. Let $p,q \in \hy D$ and $p- q \in \mu (0).$ Since $\hy D \subset \bigcup\{\mu (p) \mid p\in D\},$ then  $p,q \in \mu (s)$ for some $s \in D.$ Thus, from continuity, $\hy f(p), \ \hy f(q) \in \mu (f(s))$ implies that $\hy f(p) - \hy f(q) \in \mu (0)$ and the result follows. \qed\parm
Since being uniformly continuous is so important within analysis, I'll present a few more pertinent propositions.\parm
{\bf Theorem 10.7.} {\it For real numbers $a < b$, let $f\colon (a,b) \to \real.$ If $h \in \mu (b)^- \cap (a,b) \ \r (= \mu (b)^-\r )$ and $\hy f(h) \notin G(0),$ then for each $r \in \real^+$ and each $p \in (a,b)$ there exists some $q \in (a,b)$ such that $p < q  < b$ and $\vert f(p) -f(q)\vert > r.$ \r (A similar statement holds for the end point ``a.''\r )}\pars
Proof. Let $h \in \mu (b)^-$ and $p \in (a,b).$ Then $p < h < b.$ Since $\hy f(h) \notin G(0),$ then $\forall \, r\in \real,\ \vert \hy f(h) \vert >r$ and, by reverse *-transform, there exists $q \in (a,b)$ such that $p< q < b$ and $\vert f(q) \vert > \vert f(p) \vert + r$ for the $\vert f(p)\vert +r \in \real.$ Hence, $\vert f(p) - f(q) \vert > r.$ \qed \parm
{\bf Theorem 10.8.} {\it Let $f \colon (a,b) \to \real$ be uniformly continuous. Then $f(b-),\ f(a+) \in \real.$}\pars
Proof. Let $h \in \mu (b)^- \cap \Hy (a,b)= \mu (b)^-.$ Then $h < b$. Assume that $\hy f(h) \notin G(0).$ Since $h \in \Hy (a,b)$. Then, by *-transform of the conclusion of  Theorem 10.7, there exist $q \in \Hy (a,b)$ such that $h < q< b$ and $\vert \hy f(h) -\hy f(q) \vert > 1.$ Since $q \in \mu (b)^-,$ this contradicts Theorem 10.5. Hence, $\hy f(h) \in G(0)$ for each $h \in \mu (b)^-$ implies that, for a particular $h,$  $\st {\hy f(h)} = L.$ Now if $k \in \mu (b)^-,$ then uniform continuity implies that $\hy f(k) \in \mu (\st {\hy f(h)}.$ Thus, $f(b-) = L.$ The result for $f(a+)$ is obtained in a similar manner with a similar Theorem 10.7 for $a$ and this completes our proof. \qed \parm 
Let nonempty $E \subset D$ and $f\colon E \to \real.$ A function $g\colon D \to \real$ is called an {{\bf extension of $f$}} iff for each $x\in E,\ g(x) = f(x).$\parm
{\bf Theorem 10.9.} {\it Let $f\colon (a,b) \to \real$ be uniformly continuous on $(a,b)$. Then there exists an extension $g$ of $f$ such that $g\colon \real \to \real$ is uniformly continuous.}\pars
Proof. Simply use the last theorem and define $g(x) =f(x),$ for each $x \in (a,b)$ and $g(x) = f(a+)$ for all $x \leq a,$ and $g(x) = f(b-)$ for all 
$x \geq b.$ It's clear that $g$ is uniformly continuous on $\real.$ \qed\parm
{\bf Theorem 10.10.} {\it Let $f\colon D \to \real$ be uniformly continuous on each bounded $B \subset D.$ Then $f$ has a unique continuous extension $g\colon {\rm cl}(D) \to \real.$}\pars
Proof. A function like $f\colon D \to \real$ is uniformly continuous on each bounded $B\subset D$ iff for each $p,q \in \hy D \cap G(0)$ such that $p-q \in \mu (0),$ it follows that $f(p) - f(q) \in \mu (0).$ Since $p,q \in \hy D \cap G(0)$ iff $p,q \in \hy D\cap \Hy [-a,a]$ for some $a \in \real.$\pars
Define $g \colon {\rm cl}(D) \to \real$ as follows: let $g(\st x) = \st {\hy f(x)}$ for each $x \in \hy D\cap G(0).$ This function is well defined since if $\st y = \st x$, then $\st {\hy f(x)} = \st {\hy f(y)}$ for each $x,y \in \hy D\cap G(0)$ by uniform continuity. Further, as we know, ${\rm cl} (D) = \{ x \mid \mu (x) \cap \hy D \not= \emptyset\} = \{ \st y \mid y \in \hy D \cap G(0)\}.$\pars
 Now $g$ extends $f$, for if $ x \in D,$ then $g(x) = g(\st x) = \st {f(x)} = f(x).$ Now it's necessary to show that $g$ is continuous for any  $p \in {\rm cl}(D).$     
Let $B = D \cap [-1+p,p+1]$ and $r \in \real^+.$ Then there exists some $w \in \real^+$ such that $w < 1$ and for each $x,y \in B$ such that $\vert x - y \vert < w$ it follows that $\vert f(x) - f(y) < r/2$ by uniform continuity on bounded subsets of $D$. By *-transform, for each $x,y \in \hy B$, and $\vert x - y\vert < w$, it follows that$\vert \hy f(x) - \hy f(y)\vert < r/2.$ Let $b \in {\rm cl}(D)$ and $\vert b - p \vert < w.$ Then $b = \st y,\ p = \st x$ for some $x,y \in \hy B$ and $\vert x-y\vert < w.$ Consequently $\vert \hy f(x) -\hy f(y) \vert < r/2$ implies that $\vert \st {\hy f(x)} - \st {\hy f(y)} \vert \leq r/2.$ Therefore, in the usual manner, we have that  $\vert g(b) - g(p) \vert < r.$ Hence, $g$ is continuous at $p$. Finally, $g$ is unique for if $h \colon {\rm cl}(D) \to \real$ continuously and extends $f$, then for $p \in {\rm cl}(D),\ h(p) = h(\st x) = \st {\hy h(x)} = \st{\hy f(x)} = g(p),$ where $x \in \hy D$ and $p = \st x.$ \parm
I need just one more extension result for the next chapter.\parm
{\bf Theorem 10.11.} {\it Let $f\colon D \to \real$ continously and for each $p \in {\rm cl}(D) -D,$ assume that $\lim_p f(x) \in \real.$ Then $f$ has a unique continuous extension $g \colon {\rm cl}(D) \to \real$.}\pars
Proof. Obviously, for each $ p \in {\rm cl}(D),$ we should let $g(p) = \lim_p f(x)$ and $g\colon {\rm cl}(D) \to \real$ is unique, since the limits are unique and extends $f$ for $\lim_b f(x) = f(b),\ b \in D.$ Let $p \in {\rm cl}(D).$ Then $\mu (p) \cap \hy D \not= \emptyset.$ Let $W= (-w +g(p),g(p) +w)$ be an open interval about $g(p).$ Clearly, $ \mu (g(p)) \subset \Hy W.$ and there exists $r \in \real^+$ such that $(-r + g(p),g(p) + r) \subset [-r +g(p), g(p) + r] \subset W.$ Since $\hy f[\mu (p) \cap \hy D] \subset \mu (g(p)),$ it follows that there exists an $r_1 >0$ such that for each open interval open $I_s= (-s +p,p+s), 0< s\leq r_1$ about $p$ such that  $f[I_s\cap D]\subset (-r+g(p),g(p)+r)$ and $I_s\cap D \not= \emptyset.$ Let $q \in I_s\cap ({\rm cl}(D).$ Then since $\mu (q) \cap \hy D \not=\emptyset$ and $\mu (q) \cap \hy D \subset \Hy I_s \cap \hy D,$ it follows that $\emptyset\not= \hy f[\mu (q) \cap \hy D] \subset \Hy (f[I_s\cap D]) \subset \Hy (-r + g(p),g(p) + r).$ From the definition of $g$, we have that $\mu (g(q))\cap \Hy (-r + g(p),g(p) + r) \not= \emptyset $ implies that $g(q) \in \Hy [-r +g(p),g(p) + r].$ Thus $g(q) \in W.$ This yields that $g[I_s\cap ({\rm cl}(D))] \subset W.$ Since $W$ is an arbitrary open interval about $g(p),$ then $\mu (g(p)) = \bigcap \{ \Hy (-r + g(p),g(p) + r) \mid r \in \real^+\}$ and $\mu (p) \cap \Hy ({\rm cl}(D)) \subset (\Hy I_s \cap (\hy {\rm cl}(D)))$ imply that  $\hy g[\mu (p) \cap \Hy ({\rm cl}(D))] \subset \mu (g(p))$ and the proof is complete. \qed \parm
{\bf Corollary 10.12.} {\it Let continuous $f \colon D \to \real$, $D$ be bounded and for each $p \in {\rm cl}(D) - D$, $\lim_p f(x) \in \real.$ Then $f$ has a unique uniformly continuous extension $g\colon {\rm  
cl}(D) \to \real.$}\pars
Proof. Since ${\rm cl}(D)$ is bounded and closed it is compact. Then the unique continuous extension $g$ defined on ${\rm cl}(D)$ by Theorem 10.11 is           
 uniformly continuous.\pars 
{\bf Example 10.13.} Assume that you have defined for $x > 1$ the exponential function $x^r$ for each rational $r \in Q$ and that you have shown that it is a strictly increasing function. Then due to the fact that $Q$ is dense in $\real,$ it follows that on $f\colon Q \to \real$, where $f(x) =x^r = \sup\{f(x)\mid (x \in Q)\land (f(x)\leq r\},$  $f$ is a continuous function. But, ${\rm cl}(Q) = \real.$ Now from the completeness of the real numbers, given any irrational $r$, then, in $\real,$ 
$\lim_rf(x)=\sup\{f(x)\mid (x \in Q)\land( x< r)\}.$ Thus, by Theorem 10.11, there is a unique continuous extension $g$ of $f$ such that for irrational $r \in \real$, $g(r) = \lim_rf(x)$ and this is the value of this exponential defined at $r.$\par
\vfil\eject
\centerline{\bf 11. BASIC DERIVATIVE CONCEPTS}
\parm
I now come to the most striking difference between nonstandard analysis and the standard approach. Although the intuitive notions of the calculus are based upon ``infinitesimal modeling,'' it was precisely the logical difficulties that occured using the intuitive infinitesimal approach that greatly influenced its abandonment. No such difficulties occur for these nonstandard infinitesimals. For what follows, please notice that $D\cap D'$ is the set of all cluster points that are members of $D$ and as such $\mu'(p) \cap \hy D \not= \emptyset$ and $p \in D$. These are the members of $D$ that are not isolated points. I'll denote $D\cap D' = D_{NI}.$\parm 
{\bf Definition 11.1.} ({\bf The Standard Derivative.}) Let $p \in D_{NI}, \ p + h \in D, \ h \not= 0$ and $f\colon D \to \real.$ Then the derivative at $p$ (denoted by $f'(p)$) is {{\bf finite}} and has value $f'(p)$ iff $\lim_0 (f(p +h) - f(p))/h = f'(p).$ The derivative $f'(p) = \pm \infty$ iff $\lim_0 (f(h +p) - f(p))/h =\pm \infty$ and has geometric applications to the notion of ``vertical'' points of inflection.\parm
In all that follows, I'll use, as was done originally, the symbol $dx$ to denote a member of $\mu'(0).$ This idea of $dx$ being a special type of number was not carried over by Weierstrass when he refined the limit concept. The next theorem follows immediately from our characterizations for the limit notion. \parm
{\bf Theorem 11.2.} {\it Let $f\colon D\to \real$. Then for  $p\in D_{NI}$, $f'(p) = s \in \real$ \r [resp. $\pm \infty$\r ] iff for each $dx \in \mu'(0)$ such that $p + dx \in \hy D$
$${{\hy f(p + dx) - f(p)}\over{dx}} \in \mu (s)\ \r [resp. \ \real^\pm_\infty\r ].$$}\parm
Note that if you let $D = [a,b], \ a<b$ and $f'(a)$ exists, then $f'(a)$ is but the ``right-hand'' one-sided derivative. Clearly, this definition extends slightly the concept as it appears in the usual basic calculus course. The idea for the derivative is that it is a type of rate of change in infinitesimal values. In important physical applications, we need to know how infinitesimal rates of change compare with  ordinary real number rates of change. Let $f\colon D \to \real$, $y = f(x), \  x \in D,\ h \in \real$ such that $x+h \in D.$ Then usually one writes the {{\bf increment of (for) $y$}} at $x,$ and $h$  as $\Delta y = f(x+h) -f(x)= \Delta f(x,h).$ This $f$ generated function $(\Delta f)(p,h)$ is actually a function that determines a hyperfunction by *-transform for any $q \in \hy D,\ k \in \hyperreal$ such that $q + k \in \hy D.$ The *-transform states that $\Hy (\Delta f)(q,k) = \hy f(q +k) - \hy f(q) = (\Delta \hy f)(q,k) = \Delta \hy f(q,k).$ Thus, $f'(p) = s$ [resp. $\pm\infty$] iff $\Delta \hy f(p,dx)/dx \in \mu (s)$ [resp. $\real_\infty^\pm$] for each $dx\in \mu'(0)$ such that $p + dx \in \hy D.$ \parm 
{\bf Theorem 11.3.} {\it Let $f\colon D\to \real.$ Then $f$ is continuous at 
$p \in D$ iff $\Delta \hy f(p,dx) \in \mu (0),$ for each $dx\in \mu (0)$ such that $p + dx \in \hy D.$}\parm
{\bf Theorem 11.4.} {\it If $f\colon D \to \real$ and for $p \in D_{NI}, f'(p) \in \real,$ then $f$ is continuous at $p.$}\pars 
Proof. Assume that $f'(p) \in \real.$ Then for each $dx \in \mu'(0)$, such that $p + dx \in \hy D$
$$ {{\hy f(p +dx) - f(p)}\over{dx}} \in \mu (f'(p)).$$ 
Note that there always exists at least one such $dx.$ Hence, 
$\hy f(p +dx) - f(p) \in \mu (0)$ implies that $\hy f[\mu (p)\cap \hy D] \subset \mu (f(p))$ and the result follows. \qed \parm 
Our next notion is that of the differential. This is where we return to the time of Newton and Leibniz, something that could not be done prior to 1961. I mention that there are different approaches to the notion of the differential, especially for multi-variable functions. \parm 
{\bf Definition 11.5.} ({\bf The Differential}.) Let $f \colon D \to \real,\ f'(p)\in \real,\ dx\in \mu (0), \ p +dx \in \hy D.$ Then the {{\bf differential}} is $df=  f'(p)\, dx\in \mu (0).$
  \parm
{\bf Theorem 11.6} {\it Let $f \colon D \to \real$. If  $p \in D,$ $f'(p) \in \real$, then $f'(p) = \st {{{df}\over{dx}}}$ for each $dx \in \mu'(0)$ such that $p + dx \in \hy D.$}\pars
Proof. Immediate.\parm 
We need a better understanding of when the derivative exists and its relation to the differential and the infinitesimal increment. For this reason, let's call a function $h(p,q)$ defined on $A \times B \subset \hyperreal \times \hyperreal,$ where $\mu (0) \subset B,$ an {{\bf infinitesimal function}} at $p \in A$ iff $h(p,dx) \in \mu (0),\ \forall \, dx \in \mu(0).$ \parm

{\bf Theorem 11.7.} {\it Let $f\colon D\to \real,\ p \in D_{NI}.$ Then $f'(p) \in \real$ iff there exists a unique $t \in \real$ and   an infinitesimal function, $h\colon \{p\} \times \mu(0)\to \hyperreal$ such that for each $dx \in \mu'(0)$, where $p +dx \in \hy D,$ `
$$\Delta \hy f(p,dx)= \hy f(p+dx) - f(p) = (dx)t + (dx)h(p,dx).$$}
Proof. For the necessity, simply define $h(p,dx) = (\hy f(p+dx) - f(p))/dx - f'(p),\ dx \not= 0$ and $h(p,0) = 0.$ Then let $t = f'(p).$ It follows that  $h(p,dx) \in \mu (0),\ \forall\, dx \in \mu (0)$ and that 
$\hy f(p+dx) - f(p) = (dx)t + (dx)h(p,dx),\ \forall\, dx \in \mu'(0).$ 
The fact that $t$ is unique follows from the definition of the derivative and the disjoint nature of the monads. \pars
For the sufficiency, let $\hy f(p+dx)-f(p) = (dx)t + (dx)h(p,dx),\ dx \in \mu'(0)$ for each $dx \in \mu'(0)$ such that $p + dx \in \hy D,$ then $(\hy f(p+dx) - f(p))/dx - t = h(p,dx) \in \mu (0)$ implies that $t= f'(p).$ \qed\parm
Note that Theorem 11.7 holds in all cases including the case that $f$ is constant on some interval about $p$. The significance of Theorem 11.7 is that there are collections of infinitesimals called {{\bf order ideals}} that give a type of measure as to how well the differential approximates the infinitesimal increment. For example, the facts are that for a fixed $dx > 0,$ say, and, $f'(p) \not= 0,$ then the set $o(dx)=\{\gamma(dx)\mid \gamma \in \mu (0)\}$ generates an ideal that's a subset of $\mu (0)$ with a lot of properties. Obviously, $dx\, h(p,dx) \in o(dx).$ One says that $df$ is a {{\bf first-order}} approximation for $\Delta \hy f(p,dx)$ for each $dx.$\pars
This notion of infinitesimal approximation is exactly how ``curves'' were viewed in the time of Newton and Leibniz. From Theorem 11.7 we have specifically that $\hy f(p + dx)= f(p) + df + dx\, h(p,dx)$ holds $\forall\, dx \in \mu (0)$. Thus within $\mu (p)$, the monadic neighborhood about $p$, the  
*-line segment $g(dx) = f(p) + dx\, f'(p),\ dx \in\mu (0)$ is a first-order approximation for any $dx$ to the *-graph $y = \hy f(p+dx).$ Of course, this can be phrased in terms of *-range values. One of the original definitions for a curve was that it is an infinite collection of infinitely small line segments. So, once again, we have a rigorous formulation for the original intuitive idea. And, yes, under certain circumstances there are ``higher order'' approximations. \par
Although it's obvious from limit theory that the sum and product of functions $f,\, g$ that are differentiable at $p$ (i.e. this means that $f'(p),\ g'(p) \in \real$) are differentiable at $p$, the following two theorems demonstrate how easily the derivative ``formula'' and the chain rule are obtained.\parm
{\bf Theorem 11.8.} {\it Let $f,\ g \colon D \to \real$ and $f'(p),\ g'(p) \in \real.$ Then\pars
\indent\indent {\rm (i)} if $ u = (f)(g),$ then $u'(p) = f(p)g'(p) + f'(p) g(p);$ \pars
\indent\indent {\rm (ii)} if $g(p) \not= 0$ and $u = f/g$, then $$u'(p)={{g(p)f'(p) - g'(p)f(p)}\over{g(p)^2}}$$.}\pars
Proof. (i) For $dx \in \mu'(0)$ such that $p + dx \in \hy D,$ $[\hy f(p+dx)][\hy g(p+dx)] = [f(p) + dx\, f'(p) + dx\; h(p,dx)][g(p) + dx\, g'(p) + dx\; k(p,dx)] = f(p)g(p) + (f(p)g'(p) + g(p)f'(p))dx + \gamma \, dx$ where $\gamma \in \mu (0)$, and the result follows.\pars
(ii) For $dx \in \mu'(0)$ such that $p + dx \in \hy D,$ 
$$\Delta \hy u(p,dx) = {{\hy f(p+dx)}\over{\hy g(p+dx)}} - {{f(p)}\over{g(p)}}= {{\Delta \hy f(p,dx) + f(p)}\over{\Delta \hy g(p,dx) + g(p)}}-{{f(p)}\over{g(p)}}=$$
$${{g(p)\Delta \hy f(p,dx) - f(p) \Delta \hy g(p,dx)}\over {g(p)(\Delta \hy g(p,dx) + g(p))}}.$$ Thus, for $dx \not= 0,$ 
$${{\Delta \hy u(p,dx)}\over{dx}} = {{g(p){{\Delta \hy f(p,dx)}\over{dx}} - f(p){{\Delta \hy g(p,dx)}\over{dx}}}\over{g(p)(\Delta \hy g(p,dx) + g(p))}}.$$ 
The result follows by taking the standard part operator and using the fact that $\st {\Delta \hy g(p,dx)} = 0.$ \qed\parm
{\bf Theorem 11.9.} {\it Let $f\colon D \to \real,\ p \in D_{NI},\ g\colon f[D] \to\real, \ f(p) \in f[D]_{NL}.$ If $f'(p), \ g'(f(p)) \in \real$, Then  for the composition $(gf)(x)= g(f(x)),\ x \in D, \ (gf)'(p) \in \real$ and $(gf)'(p) = g'(u)f'(p), \  u = f(p).$ }\pars
Proof. Let $dx \in \mu'(0),\ p + dx \in \hy D.$ Then $\hy f(p+dx) = f(p) + k,\ k \in \mu'(0)$ by continuity. Hence, $\hy g(\hy f(p+dx)) - g(f(p)) = k(g'(f(p)) + k(h_g(f(p),k))$ by Theorem 11.7, which also holds if $f$ is a constant in any interval about $p$. Consequently, $\hy g(\hy f(p+dx))- g(f(p)) = (\hy f(p+dx) -f(p))g'(f(p))+ (\hy f(p +dx) - f(p))h_g(f(p),k)= f'(p)g'(f(p))dx +g'(f(p))dx\,h_f(p,dx) + f'(p)dx\,h_g(f(p),k) + \gamma\, h_f(p,dx)h_g(f(p),k),\ \gamma \in \mu (0).$ However, $g'(f(p))dx\, h_f(p,dx) + f'(p)h_g(f(p),k) + \gamma h_f(p,dx)h_g(f(p),k) \in \mu (0)$, for each $dx \in \mu'(0)$ and the result follows. \qed \parm
{\bf Theorem 11.10.} {\it Suppose that $f\colon (a,b) \to \real,\ a < b,$ has a derivative for each $p \in (a,b)$ and both $f,\ f'$ are uniformly continuous on $(a,b).$ Then there is an uniformly continuous extension $g \colon \real \to \real$ that extends $f$ and $g'$ is a uniformly continuous extension $f'.$}\pars
Proof. We know that $f'(b-),\ f'(a+),\ f(b-),\ f(a+)$ exist. The result follows by defining
$$g(x) =\cases{f(x)&$x\in (a,b)$\cr f(b-) + f'(b-)(x-b)&$x\geq b$\cr
f(a+) + f'(a+)(x-a)&$x \leq a$\cr}.$$
 Let $p-q \in \mu (0)$. If $p, q \in \Hy (a,b),$ then the result follows from the hypothesis. If $p,q \in \Hy [b, +\infty),$ then $\hy g(p) -\hy g(q) =f(b-) + f'(b-)(p -b) - f(b-) - f'(b-)(q -b) = f'(b-)(p-q) \in \mu (0)$ and in like manner if $p,q \in \Hy (-\infty, a].$ Let $p\in \Hy (a,b),\ q \in \Hy [b,+\infty), q\approx b$. Then $q\approx b$ implies, since $p\approx q,$ that $p \approx b$ and $\hy g(p) = \hy f(p) \approx f(b-).$ Now $\hy g(q) = f(b-) + f'(b-)(q-b) \approx f(b-)$, since $q-b \in \mu (0).$ Hence, $\hy g(p) - \hy g(q) \in \mu (0).$ In like manner, for $(-\infty,a]$ and for $g'$. The fact that both $g$ and $g'$ are uniformly continuous follows from Theorem 10.5 and the proof is complete. \qed \parm 
The basic calculus I idea of the local (relative) maximum or local minimum point requires in the definition quantification over the set of all open  intervals about $p \in 
D.$ The {{\bf interior}} of a set $D$ denoted by {\rm int}$(D)$ is the set of all interior points, where by Theorem 8.8, $p\in D$ is in ${\rm int}(D)$ iff 
$\mu (p) \subset \hy D.$ Theorem 8.8 eliminates one quantifier from the basic definition. Does a similar elimination happen for a local maximum or local minimum? I'm sure you recall the definition relative to the existence of an interval about $p$ that is contained in $D$. The quantifier eliminated is the ``there exists.''\parm
{\bf Theorem 11.11.} {\it Let $f \colon D \to \real$. A point $p \in {\rm int}(D)$ determines a local maximum \r [resp. minimum\r ] iff $\hy f(q) \leq f(p),$ \r [resp. $\geq\r ]\ \forall\, q \in \mu (p).$}\pars
Proof. The necessity follows from the definition and Theorem 8.8.\pars
For the sufficiency, assume that for every  $r \in \real^+$ such that $(-r+p,p+r) \subset D$ there exists $q_r \in (-r +p,p+r)$ such that $f(p) < f(q_r).$ By *-transform, we have that if $r \in \mu (0)^+,$ there is some $q_r\in (-r +p,p +r)$ such that $f(p) < \hy f(q_r).$ However, $q_r \in \mu (p) \subset\hy D;$ a contradiction and this completes the proof for the local maximum. The local minimum is similar and the proof is complete. \qed\parm
Now for the major theorem used to find many of the local maximums or minimums. But, this theorem does not restrict the derivative in the hypothesis to only finite derivatives, although the conclusion will do so. \parm
{\bf Theorem 11.12.} {\it Let $f \colon D \to \real$ and $p \in {\rm int}(D).$ If $f$ is differentiable at $p$ and $p$ is a local maximum or minimum, then $f'(p) = 0.$}\pars
Proof. First, let $p$ be a local maximum and $f'(p) \in \real.$ Then for $dx \in \mu(0)^+$ $\hy f(p +dx) \leq f(p)$ and $\hy f(p-dx)\leq f(p).$ Hence,
$${{\hy f(p+dx) - f(p)}\over{dx}} \leq 0 \leq {{\hy f(p-dx)- f(p)}\over{-dx}}.$$ 
The result follows by taking the standard part of this inequality. \pars
I now show that we cannot have that $f'(p) = \pm \infty.$ Suppose that $f'(p) = +\infty.$ Then for each $dx \in \mu'(0)^+,\ (\hy f(p +dx) -f(p))/dx > 1.$ This gives that $\hy f(p+dx) - f(p) > dx > 0.$ Therefore, 
$f(p + dx) > f(p) +dx \geq \hy f(p+dx) + dx$ from Theorem 11.11. This implies the contradiction that $dx < 0.$ By considering a $-dx$, it also follows that $f'(p) \not= -\infty.$ In similar manner, the result holds for the local minimum and the proof is complete. \qed\parm 
Prior to a generalization of Rolle's theorem, we need the notion of the boundary of a set $D$. First, recall that if $D$ is bounded, then ${\rm cl}(D)$ is bounded. The {{\bf boundary}} of $D$, $\partial D,$ is exactly what you think it should be, $\partial D = {\rm cl}(D) \cap {\rm cl}(R-D).$ 
The boundary of a set is a closed set and a nonstandard characteristic is obvious. For our basic sets, continuity at a boundary may be a one-sided continuity or even continuity at isolated points. I've mostly been giving definitions and even proofs, that are easily generalized to the multi-variable calculus. \parm
{\bf Theorem 11.13.} {\it A point $p \in \partial D$ iff $\mu (p) \not\subset \hy D$ and $\mu (p) \cap \hy D \not= \emptyset.$}\parm 
{\bf Theorem 11.14.} {\it Let $f \colon D \to \real,$ $D$ be bounded, ${\rm int}(D) \not= \emptyset$ and $f$ is differentiable at each $p \in {\rm int}(D).$ Further, if $p\in {\rm int}(D)$ and $f'(p) = \pm \infty,$ then $f$ is continuous at $p.$ Finally, assume that $\lim_a f(x) = L,$ for each $a \in \partial D.$ Then there exists some $q \in {\rm int}(D)$ such that $f'(q) = 0.$}\pars
Proof. Clearly, $f$ is continuous on int$(D).$ Assume there does not exist some $q \in {\rm int}(D)$ such that $f'(q) = 0.$ First, let $D= {\rm cl}(D)$ and $p\notin {\rm int}(D).$ Since ${\rm cl}(D) =D = {\rm int}(D) \cup \partial D,$ then $p \in \partial D$ and $p \in D.$ This implies that $\lim_p f(x) = L.$ Now assume, ${\rm cl}(D) \not= D$ and that $p \in {\rm cl}(D) - D,\ p \notin {\rm int}(D).$ Then again $p \in \partial D.$ In this case by Corollary 10.12, there exists a unique continuous extension $g\colon {\rm cl}(D) \to \real$. Now cl$(D)$ is bounded and closed and, hence, is compact.
Thus, for $f$ [resp. $g$] there is an  $x_m,x_M \in {\rm cl}(D)$ where $f$ [resp. $g$] attains its minimum at $x_m$  and maximum value at $x_M$ . But, since $f(q) \not= 0$ for each $q \in {\rm int}(D)$, then $x_m,x_M \in \partial D.$ However, since $f$ is continuous on ${\rm int}(D)$ and $g$ is a continuous extension of $f$ on $\partial D,$ then $L = g(x_m) \leq g(x)\leq g(x_M) = L,$ for each $x \in {\rm cl}(D)$, implies that $g = f$ is constant on ${\rm int}(D)$. This implies that $f'(q) = 0$ for each $q \in {\rm int}(D)$; a contradiction and the result follows. (Note: The possibility that $f'(p) = \pm \infty$ for some $p \in {\rm int}(D)$ is still valid.)\parm
{\bf Corollary 11.15.} {\it \r (Rolle's theorem\r ) Let $a < b,\ f \colon (a,b) \to \real$ be differentiable at each $p \in (a,b),$ and $f(a+) = f(b-),$ then there exists some $c \in (a,b)$ such that $f'(c) = 0.$}\parm
The following has a rather involved hypothesis. All of the requirements appear necessary for this generalization of the generalized mean value theorem. (Condition (1) holds if $f$ and $g$ are continuous on $\partial D$. Further, the conclusion obviously holds under certain conditions for any $p \in {\rm int}(D)$ where $f'(p)\pm \infty$ and $g'(p) = \pm \infty.$)\parm
{\bf Theorem 11.16.} {\it Let $f\colon D \to \real,\ g\colon D \to \real,$  where $D$ is bounded and has non-empty interior. Let $f$ and $g$ be finitely differentiable at each $p \in {\rm int}(D)$. \pars
\indent\indent {\rm (1)} Let $\lim_a f(x) \in \real,\ \lim_a g(x) \in \real$ for each $a \in \partial D$ and 
$$(\lim_af(x) - \lim_bf(x))(g(a) - g(b)) = (\lim_ag(x) - \lim_bg(x))(f(a) - f(b)),\ \forall\, a,b \in \partial D.$$ Then
for each $a,b \in \partial D,$ then there is some $p \in {\rm int}(D)$ such that 
$$f'(p)(g(a) - g(b)) = g'(p)(f(a) - f(b)). \eqno (11.17)$$}\pars
Proof. Let $a,b \in \partial D$ and consider $F(x) = f(x)(g(a) - g(b)),\ G(x) = g(x)(f(a) - f(b)).$ Now let $h(x) = F(x) - G(x).$ Then $h'(x) = F'(x) - G'(x)\in \real$ for each $x \in {\rm int}(D).$ Clearly, for each $c \in \partial D,$ $\lim_ch(x) = \lim_c f(x)(g(a) - g(b)) - \lim_cg(x)(f(a) -f(b))$ and condition (1) yields that $\lim_ah(x) = \lim_bh(x)$ for each $a,b \in \partial D.$ Thus, by Theorem 11.14, there is some $p \in {\rm int}(D)$ such that $h'(p) = F'(p) - G'(p)=0$ and the proof is complete. \qed \parm
{\bf Corollary 11.18.} {\it \r (Generalized Mean Value.\r ) Let $D = [a,b], \ a\not= b$ and $f,g$ be finitely differentiable on $(a,b)$ and both are continuous at $a$ and $b$. Then there exists some $p \in (a,b)$ such that $f'(p)(g(a) - g(b)) = g'(p)(f(a) - f(b)).$}\pars
Proof. Condition (1) of Theorem 11.16 holds since $f,g$ are both continuous at $a,b.$ \parm
{\bf Corollary 11.19.} {\it Let $D$ be compact and ${\rm int}(D) \not= \emptyset$. Let continuous $f \colon D \to \real$ be finitely differentiable at each $p \in {\rm int}(D).$ Then, for each $a,b, \in \partial D,$ there exists a $p \in {\rm int}(D)$ such that $f'(p)(b-a) = f(b) - f(a).$} \parm
I conclude this chapter on basic derivative concepts, by apply Theorem 11.16 to the theory of strictly increasing [resp. decreasing] functions. \parm
{\bf Theorem 11.20.} {\it If ${\rm int}(D) \not= \emptyset,$ $f\colon D \to \real$ continuously and $f'(p) >0$ \r [resp. $< 0$\r ] for each $p \in {\rm int}(D)$, then $f$ is strictly increasing \r [resp. decreasing\r ] on every $[a,b] \subset D,\ a \not=b.$}\parm
Proof. Let $a \not= b,\ [a,b] \subset D$. Then $[a,b]$ is compact and $(a,b) \subset {\rm int}(D)$. (The ${\rm int}(D)$ is the union of the collection of all open sets that are subsets of $D$.) Hence, the conditions of Corollary 11.19 hold. Thus, for any $x < y ,\ x,y \in [a,b]$ there exists some $p\in (x,y)$ such that $f(y) - f(x) = f'(p) (y-x) >0$ implies that $f$ is strictly increasing on $[a,b].$ The proof for decreasing is similar and this complete the proof. \qed\parm
{\bf Corollary  11.21} {\it If $a < b$ and $f \colon [a,b] \to \real$ continuously and $f'(p) = 0,\ \forall\, p \in (a,b),$ then $f$ is constant on $[a,b].$}\parm
Finally, I remark that each of the previous theorems hold under *-transform and yield some rather interesting conclusions. Here are two  examples with the first a slightly  modified application of Corollary 11.19. Indeed, the major interest in the next result is when $p \approx q$ and this result is used in the next chapter. \parm 
{\bf Theorem 11.22.} {\it Let ${\rm int}(D)\not= \emptyset$ and $f\colon D \to \real$ be finitely differentiable at each $p \in {\rm int}(D).$ Then for each $p \not= q,\ p,q \in \Hy ({\rm int}(D))$ such that $[p,q] \subset \Hy ({\rm int}(D))$ there exists some $c \in \Hy (p,q)$ such that 
$$\hy f'(c) = {{\hy f(p) - \hy f(q)}\over{p-q}}.$$}\pars
Proof. By *-transform.\parm
{\bf Theorem 11.23.} {\it \r (The first L'Hospital Rule.\r ) Assume that $f\colon (a,b) \to \real,\ g\colon (a,b) \to \real$ and for each $c \in (a,b),\ f'(c), \ g'(c)\in \real$ and $g'(c) \not= 0.$ If $f(a+) = g(a+) = 0$ and $\lim_{a+} (f'(x)/g'(x))\in \real\ \r [resp.\ \pm \infty \r ],$ then   
$\lim_{a+}(f(x)/g(x)) = L.$}\pars
Proof. Let $(f'(x)/g'(x)) \to L$ as $x \to a+$ and define $f(a) = g(a) = 0.$ Then $f$ and $g$ are continuous at $a$. Then $f$ and $g$ satisfy the hypotheses of Corollary 11.18. Let $p \in \mu (a)^+$ and consider the *-transform of Corollary 11.18. Then there exists some $t \in \mu (a)^+$ such that $a < t < p$ and $L \approx \hy f'(t)/\hy g'(t) = (\hy f(p) - f(a))/(\hy g(p) -g(a))\approx \hy f(p)/\hy g(p)$ by considering the standard part operator and the fact that $\hy g(p) -g(a) \not= 0.$  Thus $\lim_a(f(a)/g(a)) = L.$ The proof for $\pm \infty$ is similar and the proof is complete. \qed \parm
Obviously, this last result holds for the substitution of $x \to b-$ for $x \to a+$.\parm
{\bf Corollary 11.24.} {\it Assume that $f\colon (c,b) \to \real,\ g\colon (c,b) \to \real$ and for each $x \in ((c,b)- \{a\}),\ f'(x), \ g'(x) \in \real$ and $g'(c) \not= 0,$ where $a \in (c,b).$ If $\lim_a f(x) = \lim_ag(x) = 0$ and $\lim_a(f'(x)/g'(x)) = L$ \r [ resp. $\pm \infty$\r ], then $\lim_a(f(x)/g(x)) = L$ \r [resp. $\pm \infty$ \r ].}\parm
{\bf Corollary 11.25.} {\it Under the hypotheses, of Theorem 11.23 \r [resp. Corollary 11.24\r ], for each $\eps,\gamma \in \mu'(0)^+$ \r [resp. $\mu (0)$\r ], it follows that $\hy f(a + \eps)/\hy g(a + \gamma) \approx \hy f'(a +\eps)/\hy g'(a + \gamma) \approx L.$} \par
\vfil\eject
\centerline{\bf 12. SOME ADVANCED DERIVATIVE CONCEPTS}
\parm
Before starting this chapter one small remainder. I will be working with non-trivial continuity and the derivative at a point $p \in D.$ In all cases, these are defined via non-isolated points. A rather simple observation is that $p\in D$ is not isolated iff there exists some $q \in \mu (p) \cap \hy D$ such that $q \not= p.$ Let's consider the ideas of the ``higher order'' differentials and their relation to ``higher order'' increments, as well as uniform differentiability, and some inverse function theorems.\parm
Recall the standard definition of the {{\bf nth-order increment,}} where it is assumed the function is appropriately defined at the indicated domain members. It's defined by the recursive expression $\Delta^n f(p,h) = \Delta (\Delta^{n-1} f(p, h)),$ where $\Delta^0f(p,h) = f(p),\ \Delta f(p,h) = f(p+h)-f(p).$
For example, $\Delta^2f(p,h)= \Delta (\Delta f(p,h)) = f(p +2h) -f(p+h) -f(p+h) +f(p) = f(p +2h) -2f(p+h) +f(p).$ Then $\Delta^3 f(p,h) = f(p+3h) -2f(p + 2h) + f(p+h) -f(p+2h) +2f(p+h) - f(p) = f(p+3h) -3f(p+2h) + 3f(p+h) -f(p).$ From this we have that for any $n\in \nat$ 
$$ \Delta^n f(p,h) = \sum_0^n(-1)^k\left(\matrix{n\cr k\cr}\right)f(p + (n-k)h)=\sum_0^n(-1)^{(n-k)}\left(\matrix{n\cr k\cr}\right)f(p + kh),$$ 
where $\left(\matrix{n\cr k\cr}\right) = n!/((n-k)! k!), \ 0\leq k\leq n$ is a ``Binomial Coefficient.'' I now consider the nth derivative $f^n.$ \parm
{\bf Theorem 12.1.} {\it For $n \in \nat'$, $b \in \real^+$ and suppose that $f^n\colon [a,a+nb]  \to \real$. Then there exists some $t \in (a,a+nb)$ such that $\Delta^nf(a,b) = f^n(t)b^n.$}\pars
Proof. This is established by induction. For $n = 1$, Corollary 11.19 yields the result. Let $g(x,b) = f(x +b) -f(x).$ Then $g^{n-1}(x,b) = f^{n-1}(x +b) - f^{n-1}(x) \in \real,$ for each $x \in [a,a + (n-1)b].$ Thus, there exists some $t_0 \in (a,a+ (n-1)b)$ such that $\Delta^{n-1}g(a,b) = g^{n-1}(t_0,b)b^{n-1}.$ Observe that $\Delta ^{n-1}g(a,b) = \Delta^nf(a,b).$ Hence, there exists some $t_1 \in (t_0,t_0 +b)$ such that $g^{n-1}(t_0,b) = f^{n-1}(t_0 +b) - f^{n-1}(t_0) = f^n(t_1)b.$ Consequently, $\Delta g^{n-1}(a,b) = g^{n-1}(t_0,b)b^{n-1}= f^n(t_1)b^n = \Delta^nf(a,b)$, where $t_1 \in (a,a +nb).$ The result follows by induction. \qed \parm 
{\bf Corollary 12.2.} {\it Let $f^n\colon [a,b] \to \real.$ Then for each $dx \in \mu (0)^+$ and $c \in \Hy [a,b)$, there exists some $t \in (c,c+ndx)$ such that $\Delta^n\hy f(c,dx) = \hy f^n(t)(dx)^n.$} \pars
Proof. This follows from *-transform and the fact that $[c,c+ndx] \subset \Hy [a,b).$\parm
Observe that Theorem 12.1 and Corollary 12.2 clearly hold for the case that $f^n\colon [a +nb,a] \to \real,$ where $b \in \mu (0)^-.$ Now define the nth order differential at $p$ for $ y = f(x)$ by $d^ny= d^nf(p) = f^n(p)(dx)^n = f^n(p)dx^n.$ Of course, $f^0(p) = f(p).$ \parm 
{\bf Theorem 12.3.} {\it Let $f^n\colon [a,b] \to \real,\ n\in \nat'.$ \pars
\indent\indent {\rm (i)} If $f^n$ is continuous at $a$, then for each $dx \in \mu (0)^+$ and each $p \in \mu (a)\cap \Hy [a,b],$ 
$$\hy f^n(p) \approx \Delta^n\hy f(p, dx)/dx^n.$$\pars
\indent\indent {\rm (ii)} For each $c \in (a,b)$ and each $dx\in \mu (0)^+,$
$$f^n(c) \approx \Delta^n\hy f(c, dx)/dx^n.$$}\pars 
Proof. (i) This obviously holds for $n = 0$ by continuity at $a$. From Corollary 12.2 and assuming that $n\geq 1,$ we have that for $dx \in \mu (0)^+$ and $ p \in \Hy [a,b),$ there is some $t \in (p,p +ndx)$ such that 
$\Delta^n\hy f(p,dx)/dx^n= \hy f(t).$
Now let $p \in \mu (a)\cap \Hy [a,b]$. Then $p\approx a \approx t$ and continuity imply that $f^n(p) \approx \Delta^n\hy f(p, dx)/dx^n.$\pars
(ii) This obviously holds for $n =0,\ 1.$ For $n \geq 2$, in order to show this, I consider only interior points and establish this by induction.
Notice that it's not required that $f^n$ be continuous on $(a,b).$   
Let $w \in \real^+$ such that nonempty $[c, +c + nw]\subset [c,c+(n-1)w] \subset [a,b].$ Such a $w$ always exists. Define $g\colon [c,c+(n-1)w] \to \real$, by $g(y,w) = f(y+w)-f(y),\ y \in [c,c+(n-1)w].$ The function $g(y,w)$ satisfies the requirements of Theorem 12.1. Hence, there exists some  $t \in (c,c +(n-1)w)$ such that $\Delta^{n-1}g(c,w) = g^{n-1}(t)w^{n-1}.$ By *-transform, we have that if $w = dx \in \mu (0)^+,$ then there exists some $t_1 \in (c + (n-1)dx)$ such that $\Delta^{n-1}\hy g(c,dx) = \hy g^{n-1}(t_1)\, dx^{n-1}.$ The definition of $g(y,w)$, and the fact that for $n\geq 2,$ in general, $g^{n-1}(y,w) = \hy f^{n-1}(y+w)- f^{n-1}(y)$ yields, by *-transform, that 
$\Delta^n\hy f(c,dx) =\Delta^{n-1}\hy g(c,dx)$ and $$ \Delta \hy g^{n-1}(c,dx)\, dx^{n-1}= 
{{\hy f^{n-1}(t_1 + dx) - \hy f^{n-1}(t_1)}\over{dx}}\, dx^n.$$ 
Consequently, 
$${{\Delta^n\hy f(c,dx)}\over{dx^n}}= 
{{\hy f^{n-1}(t_1 + dx) - \hy f^{n-1}(t_1)}\over{dx}} = $$
$${{\hy f^{n-1}(t_1+dx) - f^{n-1}(c)}\over{t_1 +dx -c}}{{t_1 +dx -c}\over{dx}} + {{f^{n-1}(c) - \hy f^{n-1}(t_1)}\over{c-t_1}}{{c-t_1}\over{dx}} =$$  
$${{\hy f^{n-1}(c +dx_1) - f^{n-1}(c)}\over{dx_1}}{{t_1 +dx -c}\over{dx}} - {{f^{n-1}(c) - \hy f^{n-1}(c +dx_2)}\over{dx_2}}{{c-t_1}\over{dx}} \approx $$ 
$$f^n(c)\, \St {\left({{t_1+dx -c}\over{dx}}\right)} + f^n(c)\, \St {\left({{c-t_1}\over{dx}}\right)}= $$ 
$$f^n(c)\, \St {\left({{t_1 +dx-c+c-t_1}\over{dx}}\right)} = f^n(c),$$ 
where $dx_1 = t_1 + dx -c,\ dx_2 = t_1 -c$ and $dx_1,dx_2 \in \mu (0)$ and the proof is complete. \qed \parm
Theorem 12.3 holds for the appropriate negative increments and these results relate directly to notion of the {{\bf nth-order}} approximation via the {{nth order ideals.}} This is because (i) yields that $\hy f^n(p)dx^n\approx \Delta^n\hy f(p,dx)$ and (ii) $\hy f^n(c)dx^n\approx \Delta^n\hy f(c,dx).$ I have mentioned the first-order ideal generated by any $dx$. For $n >1,$ the nth-order ideal are generated by the $dx^n$ and is a strict subset of the $dx^{n-1}$ $(n-1)$th order ideal. I mentioned that many of the standard theorems have useful nonstandard statements. One of these is the nonstandard mean value theorem. For any $x,y \in \hyperreal$, the nonstandard interval $[x,y]$ restricted to members of a particular set $A \subset \real$ is easily defined by *-transform. We know that $\forall \r x \forall \r y\forall \r z((\r x\in A)\land(\r y \in A)\land (\r z\in A) \to (z \in [x,y] \iff x\leq z\leq y)).$ Thus, for $p,q \in \hy A,\ p \leq q$, one simply considers the symbol $[p,q]$ for this *-transform. Such an abbreviation occurs in the *-transform of Corollary 11.19. \parm
{\bf Theorem 12.4.} {\it Let $f \colon D \to \real$ by finitely differentiable at each $p\in {\rm int}(D).$ For distinct $p,q \in \Hy({\rm int}(D))$ such that $[p,q] \subset \Hy ({\rm int}(D)),$ there is some $c \in \Hy (p,q)$ such that $\hy f'(c) = (\hy f(p) -\hy f(q))/(p-q).$}\parm 
Theorem 12.4 is useful in the study of derivatives that are also continuous. Indeed, $f\colon D \to \real$ is said to be {{\bf continuously differentiable}} on $D$ iff $f'$ is continuous on $D.$ \parm\vfil\eject
{\bf Theorem 12.5.} {\it Let $f\colon G \to \real$ be continuously differentiable on $G,$ where $G$ be an nonempty open subset of $\real.$ Then for each $c \in G$ and $p,q \in \mu (c)$, $f'(p) \approx (\hy f(q +dx) -\hy f(q))/dx.$}\pars
Proof. Observer that $\mu (c) \subset \Hy G.$ Suppose $f'$ is continuous at $c.$ Let $dx \in \mu (0)^\pm$. We have that for any $p \in 
\mu (c),\ [p,p+dx] \subset \Hy G$ or $[p +dx, p] \subset \Hy G,$ respectively. Thus by Theorem 12.4, there is some $s$ such that, in either case, $\hy f'(s) = (\hy f(p +dx) -\hy f(p))/dx.$ By continuity, for any other $q \in \mu (c), \ \hy f'(q) \approx \hy f'(s) \approx \hy f(p)$ and this completes the proof. \qed\parm 
As is very well know the derivative can exist but not be continuous. Determining when the derivative is continuous is a substantial problem. With this in mind, I show how a slight change in the conclusion of Theorem 11.2 implies the continuity of the derivative. \parm 
{\bf Definition 12.6} ({\bf Uniformly Differentiable.}) Let $f\colon D \to \real,\  c\in D_{NI}$ and $\hy f'(c) \in \hyperreal.$ Then $f$ is said to be {{\bf uniformly differentiable}} at $c$ iff for each distinct $x,y\in \mu (p)\cap \hy D$
$$f'(c) \approx {{\hy f(x) -\hy f(y)}\over{x-y}}.$$\parm
By now you should have no difficulty translating Definition 12.6 into standard terms, where $p \in D_{NI}, f'(p) \in \real$. Then such a translating gives that for any $w > 0$, there's an open interval $(-r + p,p+r)$ about  $p  \in D_{NI}$ such that whenever distinct $x,y \in (-r +p,p + r)\cap D$, then $\vert f'(p) - (f(x) -f(y))/(x-y) \vert <w$ as the equivalent statement. The uniform part is the requirement that $x,y$ be somewhat unrestricted within an interval about $p.$\parm 
{\bf Theorem 12.7.} {\it If for nonempty open $G \subset\real,\  f\colon G \to \real,\ p \in G$ and $f'$ is continuous at $p.$ Then $f$ is uniformly differentiable at $p$.}\pars
Proof. This come from Theorem 12.5 by letting $x - y = dx$. \qed\parm 
{\bf Example 12.8} Uniform differentiability was first investigated rather recently (Bahrens, 1972). It's major contribution is that a major theorem dealing with inverses, which was previously established for continuously differentiable functions, holds true for uniformly differentiable functions. There are many  functions that are uniformly differentiable at a point but not differentiable throughout any open interval about that point and, hence, not continously differentiable at that point. As an example, consider a function constructed as follows on $[-1,1]$. Consider generating a function $f$ in the following manner. For each $n>0$ generate a collection of points by the recursion starting with $x=\pm 1,\ f(\pm 1/n) = 1,\ n = 1.$ Then, for each $x = \pm 1/(n+1),$ let $f(\pm 1/(n+1)) = f(\pm 1/n) - 1/(n^2(n+1)).$  Then consider line segments, connecting successive pairs of these points as end points, as generating the function $f$ defined on $[-1,1]$. The slope of each of these line segments $n > 0$ from $(\pm 1/(n+1),f(\pm 1/(n+1))$ to $(\pm 1/n, f(\pm 1/n))$ is $1/n$. It follows that $f'(0) = 0$ and that $f$ is uniformly differentiable at $p = 0.$ However, any interval $(-r,r),\ r \in \real^+$ about $p = 0$ contains a point where $f'$ does not exist. \pars
I mention that for the real numbers if $I$ is an open set such that real $p 
\in I$, then there always exists some $r \in \real^+$ such that $p \in (-r +p,p +r) = I_p \subset I.$ The $I_p$ is an open interval {\bf about} $p$.\parm     
{\bf Theorem 12.9.} {\it Let $f \colon D\to \real,\ p\in D_{NI}$. If $I$ is an open interval about $p$, and $f$ is uniformly differentiable for each $c \in I\cap D_{NI},$ then $f'$ is continuous at $p.$}\pars
Proof. Observer that since $p \in D_{NI}$ iff $\mu'(p) \cap \hy D \not= \emptyset$ and $p \in D'$ Thus, there are a lot of these open intervals $I$ about $p$ such that $I' \cap D \not= \emptyset.$ We first have that $f'(c)\in \real,\ c \in D_{NI}.$ Let $r \in \real^+.$ Then there exists a $w \in \real^+$ such that for each  $h \in \real$ such that $0 < \vert h\vert <w,
 c + h \in D$ 
$$\left| f'(c)- {{f(c + h) -f(c)}\over{h}}\right| < r.$$
  Let $q \in \mu (p)\cap \Hy I \cap \hy D$ and any $dx \in \mu (0)$ such that $q + dx \in \mu (p) \cap (\Hy I \cap \hy D).$ Hence, considering *-transform for arbitrary $r \in \real^+,$ it follows that 
$\left|\hy f'(q) -{{\hy f(q + dx) - \hy f(q)}\over{dx}}\right| < r.$
Thus, $\hy f'(q)\approx {{\hy f(q + dx) - \hy f(q)}\over{dx}}.$
Uniformly differentiable yields that $f'(p)\approx \hy f'(q).$ The point $q$ was an arbitrary member of $\mu (p)\cap \Hy I \cap \hy D.$ Since $\mu (p) \subset \Hy I,$ this yields that $\hy f'[\mu (p)\cap \hy D] \subset \mu (f'(p))$ and the proof is complete. \qed \parm
{\bf Corollary 12.10.} {\it If nonempty open $G \subset D$ and $f\colon D \to \real$ is uniformly differentiable on $G$ \r (i.e. at each $c \in G$\r ), then $f'$ is continuous on $G$.}\parm
{\bf Corollary 12.11.} {\it Let $f\colon D \to \real,\ p \in D_{NI}$ and $f$ is uniformly differentiable at $p$. If for each $q \in \mu (p) \cap \hy D$ and $dx \in \mu'(0)$ such that $q + dx \in \hy D,\ \hy f(x) \approx (\hy f(x + dx) - \hy f(x))/dx,$ then $f'$ is continuous at $p.$}\parm
Although uniform differentiability at a point does not imply that  the derivative is continuous at that point, what does happen is that it forces $f$ to be continuous on an entire non-trivial set that contains $p$.\parm
{\bf Theorem 12.12.} {\it Suppose that $f\colon D \to \real$ is uniformly differentiable at $p \in {\rm int}(D).$ Then there exists some open interval $I\subset D$ about $p$ such that $f$ is continuous on $I.$}\pars
Proof. Since $p \in {\rm int}(D)$, there are many open intervals $I$ about $p$ such that $I \subset D.$   Let $L=  \vert f'(p)\vert + 1.$ Assume that for each open interval $I\subset D$ about $p$ there is some $y \in I$ such that $f$ is not continuous at $y$. Since $\mu (p) \subset \hy D,$ then by *-transform, each microinterval $I_\gamma= (-\gamma +p, p+\gamma),\ \gamma \in \mu (0)^+$ contains some $y$ such that $\hy f$ is not *-continuous at $y$. This translates to say that there is some $r\in \hyperreal$ such that for all $w \in \hyperreal^+$ such that $\vert x - y\vert < w$ and $\hy f$ is defined at $x$, then $\vert f\hy f(x) - \hy f(y) \vert \geq r.$  Hence, for any such $r$, there is some $x\in \mu (0)$ such that $x\not= y$ and $r/L > \vert x -y\vert.$    Thus, $\vert \hy f(x) - \hy f(y) \vert > L \vert x- \vert >0.$ Hence, $\vert (\hy f(x) -\hy f(y))/(x-y) \vert > \vert f'(p) \vert + 1.$  This contradicts uniform differentiability at $p$ and the result follows. \qed \parm
{\bf Corollary 12.13.} {\it Let nonempty open $G \subset \real$ and $f\colon G \to \real$ be uniformly differentiable at $p \in G.$ Then there exists an open interval $I$ such that $p \in I,$ and $f$ is continuous on $I.$}\parm
I'll shortly use these ideas on uniform differentiability for an investigation of how the inverse function for an appropriate differentiable function behaves.\parm
{\bf Definition 12.14.} ({\bf Darboux Property.}) A function $f \colon D \to \real$ is said to have the {{\bf Darboux property of $D$}} iff for each $a,b \in D$ such that $[a,b] \subset D$, either $[f(a),f(b)] \subset f[[a,b]]$ or $[f(b),f(a)] \subset [a,b].$ Also recall that a function $f$ on $[a,b]$ is {{\bf one-to-one}} or an {{\bf injection}} iff for each distinct $x,y \in [a,b]$ $f(x) \not= f(y).$\parm
{\bf Theorem 12.15.} {\it If $f\colon [a,b] \to \real, a < b,$ is an injection and Darboux, then $f$ is either strictly increasing or strictly decreasing. Further, $f[[a,b]]$ is a non-trivial closed interval with end points $f(a)$ and $f(b)$.}\pars
Proof. Since $f(a) \not= f(b),$ I can simply assume that $f(a) < f(b).$ Let $x,y,z \in [a,b],\ x < z < y,\ f(x) < f(y),$ but $f(x) \not< f(z)$ or $f(z) \not< f(y).$ One-to-one implies that $f(x) > f(z)$ or $f(z) > f(y).$ If $f(z) > f(y) > f(x),$ then the Darboux property implies that there exists some $w$ such that $x < w< z$ and $f(w) = f(y).$ Since $w\not= y$ this contradicts one-to-one. In like manner, for $f(x) \not> f(z).$ Therefore, $f(x) < f(z) < f(y).$ Now let $c,d \in [a,b]$ such that $ a < c < d <  b.$ Then $f(a) < f(c) < f(d)$ and $f(c) < f(d) < f(b).$ Thus, in this case, $f$ is strictly increasing and $f[[a,b]] =[f(a),f(b)].$ A similar argument shows that if $f(b) < f(a),$ then $f$ is strictly decreasing. The Darboux property now implies that $f[[a,b]]$ is a nontrivial closed interval.  \qed \parm
{\bf Corollary 12.16.} {\it Let continuous $f \colon [a,b] \to \real.$ Then $f$ is an injection iff $f$ is either strictly increasing or decreasing on $[a,b].$}\parm
I now show that there are discontinuous functions that have the Darboux property. \parm
{\bf Theorem 12.17.} {\it If $f\colon D \to \real$ is finitely differentiable on $D,$ then $f'$ has the Darboux property.}\pars
Proof. All we need to do is to consider what happens if $a,b \in D$ and $[a,b] \subset D$ and $f'(a) <f'(b)$. Suppose that $f'(a) < k < f'(b)$. Then the function $g \colon [a,b] \to \real$ defined by $g(x) = f(x) -kx$ is finitely differentiable on $[a,b].$ Hence $g$ is continuous. Thus, there is some $c \in [a,b]$ such that $g(c) \leq g(x),\ \forall\, x \in [a,b],$ (i.e. $g(c)$ is the minimum value of $g$ on $[a,b].$  Since $g'(x) = f'(x) - k$, then $g'(b) = f'(b) -k > 0.$ In like manner, $g'(a) = f'(a) -k <0. $ Let $dx \in \mu (0)^-$. Then  $(\hy g(b+dx) - g(b))/dx > 0$ implies that $\hy g(b + dx) < g(b).$ Consequently, there is some $x \in [a,b],$ by reverse *-transform, such that $g(x) < g(b).$ In like manner, there exists some $y \in [a,b]$ such that $g(y) < g(a).$ Hence, $a,b \not= c$. Therefore, $g'(c) = 0$ implies that $f'(c) = k$ and the proof is complete. \qed\parm  
I point out that there are examples of functions finitely differentiable on $[0,1]$ but with uncountable many discontinuities (Burrill and Knudsen, 1969, p. 191.) Finally, in this chapter, I'll investigate various types of inverse function theorems. Let $f \colon (a,b) \to \real$ be continuous. Then $f$ is Darboux on $(a,b).$ Thus, $f$ defined on $[c,d]\subset (a,b)$ is an injection iff $f$ is strictly increasing or decreasing on $[c,d].$ In this case, $f$ has an inverse function $f^{-1}$ such that $f^{-1}\colon f[[c,d]] \to [c,d]$ and $f^{-1}$ is an injection onto $[c,d]$ which is also strictly monotone in the same sense. \parm
{\bf Theorem 12.18.} {\it Let the injection $f\colon D\to \real$ be continuous on $D$ and $D$ is compact. Then the inverse function $f^{-1} \colon f[D] \to D$ is continuous on $f[D]$.}\pars
Proof. Let $f(p) \in f[D].$ We know from one-to-one that $\hy f$ is one-to-one and that for each $p \in D,$ $\hy f[\mu (p) \cap \hy D] = \hy f[\mu (p)] \cap \Hy (f[D]).$ However, for our purposes consider from continuity that $\hy f[\mu (p) \cap \hy D] \subset \mu (f(p)) \cap \Hy (f[D]).$ Let $q \in  \mu (f(p)) \cap \Hy (f[D]).$ Then there is some $s \in \hy D$ such that $\hy f(s) = q.$ From compactness, there is a $p_1 \in D,$ such that $s \in \mu (p_1)$ and $\hy f[\mu (p_1) \cap \hy D] \subset \mu (f(p_1))\cap \Hy (f[D])$ implies that $q \in \mu (f(p)) \cap \mu (f(p_1))$. Thus $f(p_1) = f(p).$ From one-to-one, this gives that $p = p_1.$ Consequently, $q \in \mu (p)$ implies that $q \in \hy f[\mu (p) \cap \hy D].$ Hence, $\hy f[\mu (p) \cap \hy D] = \mu (f(p))\cap \Hy (f[D]).$ One-to-one gives that $\mu (p) \cap \hy D = \hy f^{-1}[\mu (f(p)) \cap \Hy (f[D]).$ Thus, $f^{-1}$ is continuous at $f(p).$ \qed \parm
{\bf Theorem 12.19.} {\it Let $I$ be an interval with more than one point. If the injection $f\colon I \to \real$ is continuous $I$, then $f^{-1}\colon f[I] \to I$ is continuous on $f[I]$.}\pars
Proof. For any $p \in I$, $p \in [a,b] \subset I,$ for some $a,b$ such that $ a\not= b.$ \parm
{\bf Corollary 12.20.} {\it Let non-empty open $G \subset \real$ and the injection $f\colon G \to \real$ is continuous on $G$. Then $f^{-1}\colon f[G] \to G$ is continuous on $f[G].$} \parm
{\bf Theorem 12.21.} {\it For $a < b,$ let the injection $f\colon [a,b] \to \real$ be continuous on $[a,b]$. If non-zero $f'(p) \in\real, \ p \in [a,b],$ then $(f^{-1})'(f(p)) = 1/f'(p).$}\pars 
Proof. Since $f$ is continuous and one-to-one on [a,b], then $f$ is strictly monotone. Assume $f$ is strictly increasing. Then $f[[a,b]] = [f(a),f(b)], \ f(a) < f(b).$ The injection $f^{-1}\colon [f(a),f(b)] \to [a,b]$ is continuous and strictly increasing on $[f(a),f(b)]$. Let $f(p) \in [f(a), f(b)].$ Clearly $f(p)$ is a cluster point. Let $dx \in \mu (0)'$ and $f(p) +dx \in \mu (f(p) \cap \Hy [f(a),f(b)]$ and consider $h = \hy f^{-1}(f(p) +dx) - f^{-1}(f(p)).$ Since $f^{-1}$ is continuous and one-to-one, then $h \in \mu'(0)$. Further, one-to-one also implies that $\hy f(h + p) = f(p) + dx.$ Now 
$$f'(p) \approx {{\hy f(h + p) - f(p)}\over{h}} = {{dx}\over {h}}$$ and $f'(p) \not= 0$ imply, by considering properties of the standard part operator, that 
$${{1}\over{f'(p)}} \approx {{h}\over{dx}} = {{\hy f^{-1}(f(p) -dx) - f^{-1}(f(p))}\over{dx}}.$$ This completes the proof. \qed\parm 
{\bf Corollary 12.22.} {\it Let non-empty open $G\subset \real$ and the injection $f \colon G \to \real$ be continuous on $G$. If for $p \in G,\ 0\not= f'(p) \in \real,$ then $(f^{-1})'(f(p)) = 1/f'(p).$}\parm
{\bf Theorem 12.23.} {\it Let strictly monotone $f\colon D \to \real$ be continuous on compact $D$ and, for $p \in D$, $0\not= f'(p)\in \real$. Then $(f^{-1})'(f(p)) = 1/f'(p).$}\pars
Proof.  Assume that $f$ is strictly increasing and, thus, one-to-one. The $f^{-1}\colon f[D] \to D$ is continuous and strictly increasing on compact $f[D]$.  Since $p \in D_{NI}$, then $\mu (p) \cap (\hy D - D)\not= \emptyset.$ For, $q \in \mu (p) \cap (\hy D - D),$ continuity and strictly increasing imply that $\hy f(q) \in \mu (f(p)) \cap [\Hy (f[D]) - f[D]].$ Thus, $f(p)$ is a cluster point. The proof now follows as in Theorem 12.21. \qed \parm
Example 12.8 can be modified to obtain a strictly increasing function of $[-1,1]$ such that $f'(0) \not= 0.$ By Theorem 12.21, $(f^{-1})'(f(0)) =1/f'(0)$ but it is not continuous since $f'(p)$ does not exist on any open interval that contains $0.$ It is, however, uniformly differentiable. This is why the next result is a recent improvement over all other previous results relative to differentiable inverses. \parm
{\bf Theorem 12.24.} {\it Let $f\colon D\to \real$, where $D$ is compact, $f$ is strictly monotone on $D$ and at $p \in {\rm int}(D)$, $0\not= f'(p)$ is uniformly differentiable. Then $f^{-1}$ is uniformly differentiable at $f(p)$ and $(f^{-1})'(f(p)) = 1/f'(p).$}\pars
Proof. Assume that $f$ is strictly increasing on $D.$ Then $f^{-1}$ exists for $f[D].$ Uniformly differentiable implies that $f$ is continuous on some $I \subset D$, where $I$ is an open interval about $p$. Thus, there exists $[a,b],\ (a \not= b)$ such that $p \in [a,b] \subset I\subset D.$ Since $[a,b]$ is compact, then the result that $f^{-1}$ is differentiable at $f(p)$ and that $(f^{-1})'(p) = 1/f'(p)$ follows from Theorem 12.23. \pars
Assume that $dy \in \mu'(0),\ y \in \mu (f(p)) \subset \hy f[\hy D]$ such that $y + dy \in \hy f[\hy D]$ and 
$$(f^{-1})'(f(p)) \not\approx {{\hy f^{-1}(y +dy) - \hy f^{-1}(y)}\over{dy}}.$$
Then, from properties of the ``st'' operator,  
$$f'(p) \not\approx{{dy}\over{\hy f^{-1}(y + dy) - \hy f^{-1}(y)}}.$$
Observe that $y =\hy f(q)$ for some unique $q \in \mu (p) \cap \hy D.$ Further, continuity of $f^{-1}$ at $f(p)$ and increasing imply that $0\not= \hy f^{-1}(y + dy) - \hy f^{-1}(y) = h \in \mu'(0).$ Now $\hy f^{-1}(y +dy) = q + h\in \hy D$ and, of course, $q + h \in \mu (p).$ Thus, $y + dy = \hy f(q + h)$ yields $dy = \hy f(q + h) -\hy f(q).$ Therefore, 
$$f'(p) \not\approx {{\hy f(q +h) -\hy f(q)}\over{h}};$$ a contradiction of uniformly differentiable  of $f$ at $p.$ Thus, $f'(p)$ is uniformly differentiable at $f(p)$ and the proof is complete. \qed\parm
Finally, I apply some of these previously results to establish a major classical theorem on inverse functions. THE inverse function theorem.\parm
{\bf Theorem 12.25.} {\it Let $G$ be a non-empty open subset $\real.$ Let the $f'\colon G \to \real$ be continuous on $G$. Then at any $p \in G,$ where $f'(p) \not= 0,$ there exist open intervals $I$ and $U$ such that $p \in I \subset G,\ f[I]=U \subset f[G]$ and $f^{-1}$ exists $U$, $(f^{-1})'$ is continuous on U, and $(f^{-1})'(p) = 1/f'(p)$ for each $ p \in U.$}\pars
Proof. Let $p \in G$ and $f'(p) \not= 0.$ I show first that $f(p) \in {\rm int}(f[G]).$ Let $p \in (a,b) \subset G.$ Since $f'(p) \not=0$ and $f'(p)$ is continuous on $(a,b)$ there exists some open interval $I_0$ such that $p \in I_0 \subset (a,b)$ and $f'(x) >0$ for each $x \in I_0$ or $f'(x) < 0$ for each $x \in I_0.$ Thus, $f$ is strictly monotone and continuous on $I_0.$ So, $f$ is one-to-one on $I_0$. Further, there is a closed non-trivial interval $[c,d]$ and the open $(c,d) $ such that $p \in (c,d) \subset [c,d] \subset I_0.$ Consider that case where $f(c) < f(d).$ Since $[c,d]$ is compact, it follows that for $[f(c),f(d)],$ $f^{-1}\colon  [f(c),f(d)]\to [c,d]$ is continuous on $[f(c),f(d)]$ by Theorem 12.18 and , hence, continuous on $(f(c),f(d)).$ Now Theorem 12.15 implies that $f((c,d)) = (f(c),f(d)) \subset {\rm int}(f[G]).$ Theorem 12.7 yields that $f$ is uniformly differentiable on $(c,d)$. Theorem 12.24 gives that $f^{-1}$ is uniformly differentiable on $(f(c),f(d))$ and $(f^{-1})'(x) = 1/f'(y), f(y) = x$ for each $x \in (f(c),f(d)).$ Theorem 12.9 yields that $(f^{-1})'$ is continuous on $(f(c),f(d)).$ In like manner for the case $f(c) > f(d)$ and the proof is complete. \qed  \vfil\eject
\centerline{\bf 13. RIEMANN INTEGRATION}
\parm
Since I'm using an arbitrary free ultrafilter generated nonstandard model for analysis and the somewhat weak structure $\cal M$ one should not expect that $\cal M$ models all aspects needed for Riemann integration. For this reason, a few results use standard proofs. One the other hand, I'll obtain many results relative to the Riemann integral by means of proofs using nonstandard techniques. {\bf Unless otherwise stated, all functions, such as $f$, discussed in the chapter will be \underbar{bounded} and, for $a < b,$ map $[a,b]$ into $\real.$} This is a generalization for a basic analysis course of the usual Calculus I requirement that $F$ is a (bounded) piecewise continuous function defined on $[a,b].$  \pars
One of the problems with Riemann integration is that it may be stated in terms of any partition of $[a,b]$. Nonstandard analysis allows us to eliminate this ``any partition'' notion. For our purposes a {{\bf partition of [a,b]}} is but a finite collection of points $\{a=x_0 < \cdots < x_n\leq x_{n+1} = b\}.$ There are different ways to generate ``simple partitions.'' The one now introduced is considered as a very simple type of partition of $[a,b]$ and allows any positive infinitesimals to generator a nonstandard partition. \parm
{\bf Definition 13.1.} ({\bf The Simple Partition.}) In this chapter, let $\Delta x$ always denote a positive real number. For $\Delta x$, there exists a largest natural number ``n'' such that $a + n(\Delta x) \leq b.$ Define $x_n = a + n\Delta x \leq b$. Then there is a unique partition of $[a,b]$, $P(\Delta x) = \{a=x_0<\cdots <x_n\leq x_{n+1}= b\}$ such that for each $[x_i,x_{i+1}], x_{i+1} - x_i = \Delta x,\ i = 0,\ldots, i = n-1$ and  $x_{n+1}-x_n = b - (a +nx) < \Delta x$ due to the statement dealing with $n$ being the ``largest $n$'' such that $0\leq b - (a + n\Delta x).$\parm 
 It's possible that $x_n= x_{n+1}.$ Indeed, let $\Delta x = b-a$. Then $n = 1$ and $x_1 =x_2 = b.$ The existence of this unique largest $n$ can be expressed in our formal language as follows 
$$\forall \r x((\r x\in \real^+) \to \exists \r y((y \in \nat)\land (a + yx \leq b) \land \forall \r z ((z
\in \nat)(a + zx \leq b) \to (z\leq y)))).$$  
Further, for this unique $n$, there's a function from $[0,n+1] \to [a,b]$ that generates all of the partition points, where $x_0 = a$ and $x_{n+1} =b.$ It's defined by letting $x_k = (a +k\Delta x),\ 0\leq k \leq n,\ x_{n+1} = b.$ Such functions are called {{\bf partial sequences.}} For every $\Delta x,$ there exists such a partial sequence. This allows one to define what is termed as a ``fine partition'' for each positive infinitesimal. For such $\Delta x,$ the partial sequence has a hyperfinite domain since it's not difficult to show that if $\Delta x = \gamma \in \mu (0)^+$, then the unique $n$ is a member of $\nat_\infty.$ For this reason, such partial sequences generated by positive infinitesimals are often called {{\bf hyperfinite sequences.}}\parm 
{\bf Definition 13.2.} ({\bf A Fine Partition.}) Let $dx,dy,dz$ etc. denote members of $\mu (0)^+.$ For any $dx$, there exists a unique $\Lambda \in \nat_\infty$ such that $a + \Lambda dx\leq b$ and $\forall\, k \in \hypernat$ if $a+ k\, dx \leq b,$ then $k \leq \Lambda.$ The hyperfinite sequence $S\colon [0,\Lambda] \to \Hy [a,b]$ such that $x_k = (a + k\, dx),\ k \in [0,\Lambda]$ and $x_{\Lambda +1}=b$ yields a {\bf fine partition} $P(dx)$ of (for) $[a,b],$ where $x_{k+1} - x_k =dx, \ 0 \leq k < \Lambda,\ b = x_{\Lambda +1} - x_\Lambda < dx.$ \parm   
Since I only consider bounded functions, then this investigation is based upon the completeness of the real numbers. So, as usual, for each closed interval $[x_i,x_{i+1}],$ let $m_i = \inf\{f(x)\mid x_i \leq x \leq x_{i+1} \}$ and $M_i = \sup\{f(x) \mid x_i \leq x \leq x_{i+1}\}$. (Recall that ``inf'' is the greatest lower bound of a set, and ``sup'' is the least upper bound.) Of course, $m_i \leq M_i.$ \parm\vfil\eject
{\bf Definition 13.3.} ({\bf Upper and Lower Sums).} For each $\Delta x,$ and the bounded function $f$, two operators are defined as follows:
$$ L(f,\Delta x) = \left(\sum_0^{n-1}m_i \Delta x\right) + m_n(b-x_n)$$  
$$ U(f,\Delta x) = \left(\sum_0^{n-1}M_i \Delta x\right) + M_n(b-x_n).$$  
For the fixed function $f$, the {{\bf lower sum}} $L(f,\cdot) \colon \real^+ \to \real,$ and the {{\bf upper sum}} $U(f,\cdot) \colon \real^+ \to \real$ have the usual nonstandard extensions. \parm
{\bf Definition 13.4.} ({\bf Hyperfinite sums.}) For each $dx$, $\hy L( f,dx)$ and $\Hy U(f,dx)$ are called the {{\bf lower hyperfinite sum}} and {{\bf upper hyperfinite sum}}, respectively. I've used a slight abbreviation in this notation, where the $f$ in the notation is actually $\hy f.$ \parm 
{\bf Theorem 13.5.} {\it For each $dx$ and any $f$,\pars
\indent\indent {\rm (i)} $\hy L(f,dx) \leq \Hy U(f,dx),$\pars
\indent\indent {\rm (ii)} $\hy L(f,dx), \ \Hy U(f,dx) \in G(0).$}\pars
Proof. Since $f$ is bounded on $[a,b],$ then there exist $n,m \in \real$ such that $m\leq f(x)\leq M,\ \forall\, x \in [a,b].$ Consider any $\Delta x$. Then $m\Delta x  \leq f(x) \Delta x \leq M\Delta x$ yields that 
$$m\left(\left(\sum_0^{n-1} \Delta x\right) + (b -x_n)\right) \leq \left(\sum_0^{n-1} m_i\Delta x \right) + m_n(b-x_n) \leq$$
$$\left(\sum_0^{n-1} M_i\Delta x\right) + M_n(b - x_n) \leq M\left(\left(\sum_0^{n-1}\Delta x\right) + (b - x_n)\right),$$
since these are finite summations. Hence, $m(b-a) \leq L(f,\Delta x) \leq U(f,\Delta x)\leq M(b-a).$ Then the sentence
$$\forall \r x((x\in \real^+) \to m(b-a) \leq L(f,x) \leq U(f,x)\leq M(b-a))$$
holds in $\cal M$; and, hence, in $\hy {\cal M}$. By *-transform, the result follows. \qed\parm
In the theory of Riemann integration, refinements of a partition play a significant role. They also present significant intuitive problems, as well.
The next result is similar to a refinement proposition for Riemann sums.\parm
{\bf Theorem 13.6.} {\it For every $\Delta x$ and for every $p\in \hypernat' = \hypernat -\{0\},$ 
$$\hy L(f,\Delta x) \leq \hy L(f,\Delta x/p) \leq \Hy U(f,\Delta x/p) \leq \Hy U(f,\Delta).$$}\pars
Proof. It follows from the definition,  that for each $\Delta x$ and corresponding partition $P(\Delta x)$ generated by $\Delta x$ that $P(\Delta x )\subset P(\Delta x/n), n \in \nat'.$ Now I need to consider a standard argument at this point and direct you to Theorem 10.1 in Burrill and Knudsen (1969, p. 199), where it is established that, for our case, 
$$L(f,\Delta x) \leq L(f,\Delta x/n) \leq U(f,\Delta x/n) \leq U(f, \Delta x),$$
for $n \in \nat'.$ The result follows by *-transform. \qed\parm
I point out that Theorem 13.5 shows that for each positive infinitesimal $dx,$ since $\hy L(f,dx)\leq \Hy U(f,dx)$, then $\st {\hy L(f,dy)}\leq \st {\Hy U(f,dy)}.$ Also, I won't continue to mentioned the fact that $0\leq \Hy U(f,p) - \hy L(f,p)$ for each $p \in \hyperreal^+.$\parm
{\bf Definition 13.7.} ({\bf Integrable Functions.}) The function $f$ is {{\bf (simply) integrable}} iff there is some $dx$ such that 
$$\st {\hy L(f,dx)} = \st {\Hy U(f,dx)}\ {\rm iff}$$
$$\Hy U(f,dx) - \hy L(f,dx) \in \mu (0).$$
If $f$ is integrable, then we denote $\st {\hy L(f,dx)} = \int_a^bf\, dx$ as the {{\bf integral}}, where it's understood that this is the simple integral. \parm
I also use the notation $\int_a^bf\, dx \in \real$ to indicate that $f$ is integrable on $[a,b].$ At the moment, it appears that the value of the integral might depend upon the $dx$ chosen. I'll show, later, that this is not the case. Of course, it's clear from above that if $\int_a^b f\, dx \in \real$, then $\int_a^bf\, (dx/n) \in \real, \ n \in \hypernat'.$ What functions are integrable?\parm 
{\bf Theorem 13.8.} {\it If $f$ is monotone on $[a,b]$, then $\int_a^b f\, dx \in \real.$}\pars
Proof. Assume that $f$ is increasing. For an $\Delta x$ generated partition,
$M_i = f(x_{i+1}),\ m_i = f(x_i),\ i = 0,\ldots,n.$ For $n \in \nat'$, let 
$\Delta x = (b-a)/n.$ Then 
$$U(f(\Delta x) - L(f,\Delta x) = ((b-a)/n)(f(b) - f(a)).$$
By *-transform, for each $\Lambda \in \nat_\infty$, where $dx = (b-a)/\Lambda,$
$$\Hy U(f,dx) - \hy L(f,dx) = dx(\hy f(b) - \hy f(a)) \in \mu (0).$$   
 and the result follows. \qed\parm
What is needed is a general standard characterization for integrability in our sense. \parm
{\bf Theorem 13.9.} {\it The function $f$ is integrable on $[a,b],$ for some $dx \in \mu (0)^+,$ iff, for each $r \in \real^+$, there is some $\Delta x \in \real^+$ such that 
$$U(f,\Delta x) - L(f,\Delta x) < r.$$}\pars
Proof. Assume that $\int_a^b f\, dx \in \real$ and $r \in \real^+.$ Then $dx \in \mu (0)^+$ and $\Hy U(f,dx) - \hy L(f,dx) < r$. The necessity follows by reverse *-transform. \pars
For the sufficiency, assume that $r\in \real^+$ and that there exists some $\Delta x$ such that $U(f,\Delta x) - L(f,\Delta x) < r.$ If $n \in \nat'$, then it also follows that $U(f,\Delta x/n) - L(f,\Delta x/n) < r$ by Theorem 13.6 restricted to $\nat'.$ However, there always exists an $n \in \nat'$ such that $0 < \Delta x/n < r.$ Consequently, the sentence
$$\forall \r x((x\in \real^+) \to \exists \r y ((y \in \real^+)\land(y < x) \land (U(f,y) - L(f,y) < x)))$$
holds in ${\cal M}$; and, hence, in $\hy {\cal M}.$ Letting $\gamma \in \mu (0)^+,$ there exists some $dx$ such that 
$\Hy U(f,dx) - \hy L(f,dx) < \gamma$. The result follows from Definition 13.7 since $\Hy U(f,dx) - \hy L(f,dx)\in \mu (0).$ \qed \parm
{\bf Corollary 13.10.} {\it The function $f$ is integrable on $[a,b]$ iff for each $\gamma \in \mu (0)'$ there exists some $dx$ such that 
$\Hy U(f,dx) - \hy L(f,dx) < \gamma.$}\parm
I'll denote {{\bf Riemann integration}} by the symbol $R\int_a^b f\, dx.$ This form of integration is defined in terms of what appears to be a more general form of partitioning of $[a,b].$ This may be why many students, when they first encounter the complete definition for Riemann integration, find it somewhat difficult to comprehend. I'll show later that the simple integral as defined here for a rather simple type of partitioning is equivalent to the Riemann integral. There are two major but equivalent definitions for the Riemann integral. (A) You consider the same idea of the lower and upper sums, but you do not restrict the partitions. You define these sums in exactly the same way but for any partition $P'$. But, then you must also do the following. You consider the numbers $\underline{R}(f) = \sup \{L(f,P')\mid P'\ {\rm any \ partition\ of\ }[a,b]\}$ and $\overline{R}(f) = \inf \{L(f,P')\mid P'\ {\rm any \ partition\ of\ }[a,b]\}.$ If $\underline{R}(f) = \overline{R}(f)$, then this value is the Riemann integral of $f.$ For the structure I'm working with, the collection of all such partitions is not part of the structure. So, this is why I need to use some results obtained by standard means. Then there is the more familiar equivalent definition. (B) You consider a general partition
 $P' = \{a =x_0 < \cdots < x_n = b\}, \ n> 0,$ where one defines the ${\rm mesh}(P) = \max\{\Delta x_i\mid (\Delta x_i = x_{i} - x_{i-1})\land (i = 0,\ldots,n)\}$.
Then you consider any finite collection $q_i \in [x_{i-1},x_{i}]$ and evaluate the function at these points and consider the Riemann sum $\sum_1^n f(q_i) (x_{i}- x_{i-1}).$ Then a number $R\int_a^b$ is the Riemann integral iff for each $r \in \real^+,$ there exists a $w \in \real^+$ such that for every partition $P'$ such that ${\rm mesh}(P') <w$ and every $q_i \in [x_{i-1}, x_i]$, 
$$\left| \sum_1^n f(q_i)(x_i - x_{i-1}) - R\int_a^bf\, dx \right| < r.$$
Although the notation contains the symbol $dx$, infinitesimals are not mentioned in definitions (A) and (B).  For bounded $f$, I use Definition (A) for the Riemann integral since the only difference is in the collection of partitions needed. For $dx,$ the partition notation $P(dx)$ is an abbreviation for the fine partition that can be explicitly defined for the $dx.$\parm
{\bf Theorem 13.11.} {\it If $\int_a^b f\, dx \in \real$, then $\int_a^b f\, dx = R\int_a^b f\, dx.$}\pars
Proof. Let $U(f,P'), \ L(f,P')$ be the upper and lower Riemann sums for a any general partition $P'$. Here is where I need a standard result about Riemann integration. It states that $R\int_a^b f\, dx \in \real$ iff, for each $r \in \real^+,$ there     
exists a general partition $P'$ such that $U(f,P') -L(f,P')< r.$ Also, $L(f,P') \leq R\int_a^b f\, dx \leq U(f,P')$ (Burrill and Knudson, 1969, p. 202). But, a simple partition $P(\Delta x)$ is a general partition. Indeed, for our partitions $L(f,P') = L(f,\Delta x),\ U(f,P') = U(f,\Delta x)$. Consequently, Theorem 13.9 yields that $R\int_a^b f\, dx \in \real.$ And, further, by *-transform, that for the *-Riemann partition $P(dx),\ \st {\Hy U(f,dx)} = R\int_a^b f\, dx = \int_a^bf\, dx = \st {\hy L(f,dx)}$ and this completes the proof. \qed\parm
{\bf Corollary 13.12.} {\it If $\int_a^bf\, dx,\ \int_a^bf\, dy \in \real$, then $\int_a^bf\, dx = \int_a^bf\, dy.$}\parm
Let's easily establish some of the basic integral properties. \parm 
{\bf Theorem 13.13.} {\it Let bounded $f$ and $g$ be integrable on $[a,b]$ for $dx$. \pars
\indent\indent {\rm (i)} For each $c \in \real,\ f+ g, \ cf$ are integrable on $[a,b]$ and $\int_a^b(f+g)\, dx = \int_a^b f\, dx + \int_a^b g\, dx,\ \int_a^bcf\, dx = c\int_a^bf\, dx,\ \int_a^b \, dx = b-a.$ \pars
\indent\indent {\rm (ii)} If $f(x) \leq g(x),\ \forall\, x \in [a,b]$ then $\int_a^b f\, dx \leq \int_a^bg\, dx.$\pars
\indent\indent {\rm (iii)} If $m \leq f(x) \leq M,\ \forall\, x \in [a,b],$ then $m(b-a) \leq \int_a^bf\, dx \leq M(b-a)$.}\pars
Proof. These are all established by simple observations about finite sums.\pars
(i) Observe that 
$$L(f,\Delta x) + L(g, \Delta x) \leq L(f+g,\Delta x) \leq U(f+g,\Delta x) \leq  U(f,\Delta x) + U(g, \Delta x).$$ 
This result follows by *-transform using the $dx$ and the standard part operator. \pars
Since $L(cf,\Delta x) = c L(f, \Delta x) \leq U(cf,\Delta x) = c\, U(f\Delta x)$, then result follows by using the standard part operator for the given $dx.$\pars
Now observe that $L(1,\Delta x) = U(1,\Delta x) = b-a.$ Thus, for the given $dx,\int_a^b\, dx = b-a.$ \pars
(ii) Clearly, for each $\Delta x,\ L(f,\Delta x) \leq L(g,\Delta x)$ implies that $\st {\hy L(f,dx)}= \int_a^bf\, dx\leq \st {\Hy U(g,dx)}= \int_a^bg\, dx$ for the given $dx.$\pars
\indent\indent {\rm (iii)} Simply apply (i) and (ii) and the proof is complete. \qed\parm 
Monotone bounded functions need not be continuous, but they are integrable. The continuous functions should be integrable or this integral would not be very useful. There are very short nonstandard proofs of the following result but they require a more comprehensive structure than I'm using. The nonstandard proof that establishes the next result is longer than the standard proof since I've defined integration nonstandardly. So, I'll give the usual standard proof that depends upon the standard characterization of Theorem 13.9.\parm
{\bf Theorem 13.14.} {\it If $f$ is continuous on $[a,b]$, then  $\int_a^bf\, dx\in \real$ for some $dx \in \mu (0)^+.$}\pars
Proof. Let $r \in \real$ and let $c = r/(b-a)$. From uniform continuity, there is a $w \in \real^+$ such that for any $x,y \in [a,b]$ such that $\vert x-y\vert < w,$ then $\vert f(x) - f(y) \vert < c.$ Consider any simple partition $P(\Delta x)$ for $[a,b].$ Let $[x_i,x_{i+1}]$ be one of the subdivisions. Then there is $x',y' \in [x_i, x_{i+1}]$ such that $m_i = f(x'_i),\ M_i = f(y'_i).$ Hence, $U(f,\Delta x)) -L(f,\Delta x) = \sum_0^{n-1} (f(y'_i) - f(x'_i))\Delta x + (f(y'_n) - f(x'_n))(b - x_n) < c(b-a) = r.$ The result follows from Theorem 13.9 and the proof is complete. \qed\parm
I mentioned previously that integration as here defined is equivalent to Riemann integration. It's time to establish this. But, due to the weak structure I'm using, I need one more standard result about general partitions.  \parm
{\bf Theorem 13.15.} {\it For each $r \in \real^+,$ there exists some $w \in \real^+$ such that for all partitions $P',$ where ${\rm mesh}(P') < w,$
$$0\leq \underline{R}(f) - L(f,P') < r,\ 0\leq U(f,P') - \overline{R}(f) <r.$$}\pars
Proof. This is establish in a portion of proof of Theorem 10.28 in Burrill and Knudsen (1969, p. 223.)\parm
{\bf Theorem 13.16.} {\it For each $dx$, 
 $$\underline{R}(f) - \hy L(f,dx)\in \mu (0)^+,\ \Hy U(f,dx) - \overline{R}(f) \in \mu (0)^+.$$}\pars
Proof. Assume that there exists some $dx$ such that $\underline{R}(f) - \hy L(f,dx)\notin \mu (0)^+.$ Since $0\leq \underline{R}(f) - \hy L(f,dx)$, then there exists some $r \in \real^+$ such that $\underline{R}(f) - \hy L(f,dx)\geq r.$ Now let arbitrary $w \in \real^+.$ Then, $0 < dx < w.$ 
 But, Theorem 13.15 holds for our simple partitions where ${\rm mesh}(P(\Delta x)) = \Delta x.$ Hence, there exists some $w_1 \in \real^+$ such that for each $P(\Delta x),$ when $\Delta x < w_1,$ then  $\underline{R}(f) - L(f,P(\Delta x))= \underline{R}(f) -L(f,\Delta x)< r.$ By *-transform of this conclusion, it follows that $\underline{R}(f) - L(f,dx)< r$ since $0< dx< w_1;$ a contradiction. Hence, for each $dx$, $\underline{R}(f) - \hy L(f,dx)\in \mu (0)^+.$ In like manner, it follows that $\Hy U(f,dx) - \overline{R}(f) \in \mu (0)^+$  and the proof is complete. \qed\parm
{\bf Theorem 13.17.} {\it If $R\int_a^bf\, dx \in \real$, then $R\int_a^bf\, dx = \int_a^bf\, dx,$ for each $dx.$}\pars
Proof. Let $R\int_a^bf\, dx \in \real$. Consider arbitrary $dx$. Then $\underline{R}(f) = \overline{R}(f) = R\int_a^b f\, dx$ and, from Theorem 13.16, $R\int_a^bf\, dx -\hy L(f,dx) \in \mu (0),\ \Hy U(f,dx) - R\int_a^bf\, dx \in \mu (0)$ yield the result.  \qed \parm
So, all we need in order to study Riemann integration are the simple partitions and the simple integral that's independent from the actual $dx$ that's used. This is a significant simplification. Further, it's obvious that if one wants to generalize the Riemann integral to other functions defined on $[a,b]$ it's clear that either the partitions must be of a different type than the simple partition or the integral must be dependent upon the $dx$ used. The major generalization is called the {{\bf generalized integral}} and how under a definition similar to (B) such a generalized integral is equivalent to the Lebesgue integral. Anyone who studies the Lebesgue integral from the viewpoint of measure theory knows the subject can be difficult. 
It's a remarkable fact that the Lebesgue integral can be viewed as a Riemann-styled integral under definition (B) for specially selected partitions and specially selected values for the function. From the nonstandard viewpoint, one major difference is that the Lebesgue integral is not infinitesimal independent. The infinitesimals needed are generated by objects called L-microgauges (Herrmann, 1993, p. 217.) But, all of this is well beyond the material in this book. \par
It's a useful fact that any $dx$ can be used to obtain the integral. This allows many standard results to be established easily. Recall that by definition $\int^a_af\, dx = 0,\ \int_b^a f\, dx = - \int_a^bf\, dx.$\parm
{\bf Theorem 13.18.} {\it {\rm (i)} Suppose that $\int_a^bf\, dx \in \real$ and $c \in [a,b].$ Then $\int_a^cf\, dx,\ \int_c^bf\, dx \in \real$ and $\int_a^bf\, dx = \int_a^cf\, dx + \int_c^bf\, dx.$\pars
\indent\indent {\rm (ii)} If $c \in [a,b]$ and $\int_a^c f\, dx, \ \int_c^b f\, dx \in \real$, then $\int_a^bf\, dx\in \real$ and $\int_a^bf\, dx = \int_a^cf\, dx + \int_c^bf\, dx.$}\pars
Proof. (i) Clearly, the result holds for $c = a,\ c = b.$ So, let $c \in (a,b)$ and $\Delta x= (c-a)/n, n \in \nat'.$ Then all the points in the simple partition created by $\Delta x$ for $[a,c]$ are points in the $\Delta x$ generated simple partitions for $[c,b]$ and $[a,b].$ Hence, it follows that 
$$L(f,\Delta x, [a,b]) = L(f,\Delta x, [a,c]) + L(f, \Delta x, [c,b]) \leq$$ $$U(f,\Delta x, [a,c]) + U(f, \Delta x, [c,b]) = U(f,\Delta x, [a,b]).$$
Let $\Lambda \in \nat_\infty$ and $dx = (c-a)/\Lambda$. By *-transform and the standard part operator and the fact that $\int_a^bf\, dx \in \real,$ we have that  
$$\st {\hy L(f,dx, \Hy [a,b])} = \st {\hy L(f,dx, \Hy [a,c])} + \st {\hy L(f, dx, \Hy [c,b])} = $$ $$\st {\Hy U(f,dx, \Hy [a,c])} + \st {\Hy U(f, dx, \Hy [c,b])} = \st {\Hy U(f,dx, \Hy [a,b])}.$$
But, $\st {\Hy U(f,dx,\Hy [a,c])} - \st {\hy L (f,dx,\Hy [a,c])} \geq 0,$ 
$\st {\Hy U(f,dx,\Hy [c,b])} - \st {\hy L (f,dx,\Hy [c,b])} \geq 0$ imply that $\st {\Hy U(f,dx,\Hy [a,c])} = \st {\hy L (f,dx,\Hy [a,c])}$ and $\st {\Hy U(f,dx,\Hy [c,b])} = \st {\hy L (f,dx,\Hy [c,b])}$ and the result follows. \pars
(ii) This follows by considering the same type of simple partition as in (i) and applying the standard part operator and the proof is complete. \qed\parm 
Now let's consider a few of the most significant properties of the integral.\parm \vfil\eject
{\bf Theorem 13.19.} {\it Let $\int_a^b f\, dx \in \real.$ For each $y \in [a,b]$, let $F(y) = \int_a^yf\, dx.$ Then $F$ is uniformly continuous on $[a,b].$ If $f$ is continuous at $c \in [a,b]$, then $F'(c) = f(c).$ }\pars 
Proof. 
From Theorem 13.18, $F(y)\in \real.$ Thus, $F$ is a function from $[a,b]$ into $\real.$ Since $f$ is bounded and again by applying Theorem 13.18, it follows that for any $x,y \in [a,b]$, $\vert F(y) - F(x) \vert \leq M\vert y-x\vert,$ where $\vert f(z)\vert \leq M,\ \forall\, z \in [a,b].$ Hence, by *-transform, if $p,q \in \Hy [a,b]$ and $p-q \in \mu (0)$, then 
$\vert \Hy F(p) - \Hy F(q) \vert \leq M\vert p-q \vert \in \mu (0)$ implies by Theorem 10.5 that $F$ is uniformly continuous on $[a,b].$ \pars
Assume that $f$ is continuous at $c \in [a,b]$.  Now considering the integral as the function $F(y) =\int_a^y f\, dx,\ y \in [a,b]$, our previous integral properties can by translated into properties about $F$, where for $z,y \in [a,b],\  \int_z^yf\, dx =  F(y) - F(z).$ By *-transform, these $\Hy F$ function properties are relative to $z,y \in \Hy [a,b].$ Let $p \in \mu (c)$ such that $p + c \in \Hy [a,b].$ First, assume that $p < c$. From the continuity of $f$, $\vert \hy f(p) - f(c) \vert = \gamma \in \mu (0).$ Let $g(x) = f(x) - f(c).$ Then from the hyper-properties for $\Hy G(x)$, we have that $G(c) - \Hy G(p) = F(c) - \Hy F(p) - f(c)(c-p)$. Moreover, $\vert F(c) - \Hy F(p) - f(c)(c-p)\vert \leq \gamma(c-p).$ Consequently, $$\left|{{F(c) - \Hy F(p)}\over{c-p}} - f(c)\right| \in \mu (0).$$
In like manner, for the case that $p > c$, and the result that $F'(c) = f(c)$ follows. \qed \parm      

{\bf Corollary 13.20.} {\it If $\int_a^bf\, dx \in \real,\ p,q \in \Hy [a,b]$ and $p-q \in \mu (0)$, then $\Hy F(p) - \Hy F(q) = \hy \int_q^p \hy f\, dx \in \mu (0).$}\parm
It's beyond the scope of this book to establish a necessary and sufficient for $\int^b_af\, dx$ to exist. The facts are that there are some very unusual functions that are integrable. For example, consider the non-negative rational numbers (in lowest form) $q/p,\ p > 0.$ Define on $[0,1]$ the function $f(p/q) = 0$ and, for each irrational $r,$ $f(r) = r.$ Then $\int^1_0f\, dx$ exists. I leave it to the reader to find the exact value. It is rather easy, however, to use our methods to show that the value of the integral is independent from the value of the bounded function at finitely many points in $[a,b]$. \parm
{\bf Theorem 13.21.} {\it Let $f$ and $g$ be bounded on $[a,b]$ and there exists a non-empty finite set of numbers $D=\{p_0,\ldots,p_n\} \subset [a,b]$ such that $f$ and $g$ only differ on $D$. If $\int^b_af\, dx \in \real$, then $\int^b_af\, dx = \int^b_ag\, dx.$}\pars
Proof. Consider any $dx.$ Without loss of generality, we may assume that, for $[c,d] \subset [a,b],\ c\not= d,$ that $f$ and $g$ differ at most at the end points $\{c,d\}.$ Then $\hy L(f,dx)= m_0\, dx + \sum_1^{\Lambda -2}m_i\, dx + m_{\Lambda -1}\, dx + m_{\Lambda} (b - x_\Lambda).$ Hence, $\st {\hy L(f,dx)} = \st {\sum_1^{\Lambda -2}m_i\, dx} = \st {\hy L(g,dx)}.$ In like manner, $\st {\Hy U(f,dx)} = \st {\Hy U(g,dx)}$ and the result follows.\qed \par     
\vfil\eject
 \centerline{\bf 14. WHAT DOES THE INTEGRAL MEASURE?}
\parm
In this chapter, I present a slightly advanced look at the type of physical properties that the integral will measure. Also I continue to assume that $f$ is bounded on $[a,b].$\parm
{\bf Definition 14.1} ({\bf Additive Function.}) A function $B \colon [a,b] \times [a,b] \to \real$ is {{additive}} if for each simple partition $P(\Delta x)= \{a=x_0<\cdots<x_n\leq x_{n+1}=b\}$ 
$$B(x_i,x_{i+1}) = B(x_i,x) + B(x,x_{i+1}),\ x_i \leq x\leq x_{i+1},\ i = 0,\ldots, n.$$ (This is not the only definition in the literature for this type of additive function.) \parm
I note that in general $B(x_i,x_i) = 0.$ Of course, it's immediate that if $\int_a^bf\, dx\in \real$, then $B(x,y) = F(y) -F(x)= \int_x^yf\, dx$ is additive on $[a,b]$. But, does the converse hold? If you are given a specific additive function $B$, then $\Hy B$ has meaning for any fine partition generated by $dx$ since the definition of $B$ is relative to the partition points and subdivision closed intervals. \parm
{\bf Definition 14.2.} ({\bf Admissible for $f$.}) A function $B$ that is additive on $[a,b]$ is {{\bf admissible}} for $f \colon [a,b] \to \real$ iff there exists some $dx$ and for the fine partition $P(dx) = \{a=x_0< \cdots <x_\Lambda\leq x_{\Lambda+1}=b\}$ it generates, for each $i =0,\ldots, \Lambda -1,$ where $b - x_\Lambda = 0,$ there exists some $p_i \in [x_i,x_{i+1}]$ such that
$${{\Hy B(x_i,x_{i+1})}\over{dx}} - \hy f(p_i)\in \mu(0), \eqno (IC)$$
and if $b-x_\Lambda \not= 0$, then there also exists some $p_\Lambda \in [x_\Lambda, b]$ such that
 $${{\Hy B(x_\Lambda,b)}\over{b-x_\Lambda}} - \hy f(p_\Lambda)\in \mu(0). \eqno (IC)$$\parm
{\bf Theorem 14.3.} {\it Let $B$ be admissible for $f\colon [a,b] \to \real.$ Then for each $r \in \real^+,$ there exist $dx,dy$ such that 
$$-r(b-a) + \hy L(f,dy) < B(a,b) < \Hy Uf,dx) + r(b-a).$$}\pars
Proof. Let $r \in \real^+$. Assume that for each $dx,$ $$B(a,b) \geq \Hy U(f,dx)= \Hy U(f,dx) + r(b-a). \eqno (14.4)$$
 I make the following observation about $B$ where $B(a,b) \geq U(f+r,\Delta x).$ Let $n >1,$ and $B(a,b) = \left(\sum_0^{n-1}B(x_i,x_{i+1} + \Delta x)\right) + B(x_n,b).$ Then there exists some $k \in [0,n-1]$ such that $ B(x_k,x_k + \Delta x)\geq M_k\Delta x$ and if $b-x_n \not=0$ (or $n=1$), then $B(x_n,b) \geq M_n(b-x_n),$ where $M_i = \sup \{f(x) + r\mid x \in [x_i,x_i + \Delta x]\}$ which exists by boundedness. Thus, by *-transform and assuming that (14.4) holds, we have that there exists $k \in [0,\Lambda-1],\ \Hy B(x_k,x_{k+1})\geq M_k\, dx$ and if $b-x_\Lambda \not=0$, then $B(x_n,\Lambda) \geq M_\Lambda(b-x_\Lambda),$ where $M_i = \sup \{f(x) + r\mid x \in [x_i,x_i + dx]\},$ which also all exist from boundedness and I need not consider the $\Hy\sup$ since by definition the $\Hy \sup = \sup.$ 
Consequently, for each $p \in [x_k,x_{k+1}], \ \Hy B(x_k,x_{k+1}) \geq (\hy f(p) + r)dx,\ k \in [0,\Lambda -1]$ and if $p\in [x_\Lambda,b]$, then, $\Hy B(x_\Lambda,b) \geq (\hy f(p) +r)(b-x_\Lambda),$ where $b- x_\Lambda \geq 0.$ This implies that for each $dx, \ k \in [0,\Lambda -1],$
$${{\Hy B(x_k,x_{k+1})}\over{dx}} - \hy f(p) \geq r,\ \forall\, p \in [x_k,x_{k+1}$$
and if $b-x_\Lambda \not= 0$, then
$${{\Hy B(x_\Lambda,b)}\over{b-x_\Lambda}} - \hy f(p) \geq r,\ \forall\, p \in [x_\Lambda,b].$$
This contradicts admissibility. Thus, there exists some $dx$ such that $B(a,b) < \Hy U(f + r,dx).$ In like manner, there exists some $dy$ such that $-r(b-a) + \hy L(f,dy) < B(a,b)$ and this completes the proof. \qed\parm
{\bf Theorem 14.5.} {\it If $B$ is admissible for integrable $f$, then $B(a,b)= \int_a^b f\, dx.$}\pars
Proof. Let $r \in \real^+.$ Then, from Theorem 14.3, there exists $dx,dy$ such that $-r(b-a) + \hy L(f,dy) < B(a,b) < \Hy U(f,dx) + r(b-a).$ The result follows by taking the standard part operator and  the fact that $r$ is arbitrary. \qed \parm
Thus, the integral can be used to calculate the values of an admissible additive function. However, the converse of Theorem 14.5 does not hold. Indeed, due the fact there are many unusual integrable functions, the converse does not hold where you define the additive function by the integral itself. In the rather simple example below, it's shown that there are additive functions, indeed integrals, where $(IC)$ holds but not for all $dx.$ \parm
{\bf Example 14.6.} Define the integrable function $f(x) = 0,\ \forall\, x \in [0,1),\ f(x) = 1,\ \forall\, x \in [1,2].$ Define $B(x,y) = \int^y_xf\, dx$ for all $x,y \in [0,2].$ Then $B(x,y)$ is additive on $[0,2].$ Let $\Lambda \in \nat_\infty,$ and let $dx = 2/\Lambda$ and $\Lambda$ by *-even.  There exists some $k \in [0,\Lambda -1]$ such that for each $p\in \Hy[0,2],\ p < x_k,\ \hy f(p) =0$ and $p\geq x_k, \ \hy f(p) = 1.$ We also know that 
for each $k \in [0,\Lambda -1],\ m_kdx \leq \Hy B(x_k,x_{k+1}) \leq M_k dx,$ where $m_k = \inf \{\hy f(x)\mid x \in [x_k,x_{k+1}]\},\ M_k = \sup  \{\hy f(x)\mid x \in [x_k,x_{k+1}]\}.$ Consequently, we have that for $j \in [0,\Lambda -1] ,\ j < k$
$${{\Hy B(x_j,x_{j+1})}\over{dx}} = 0= \hy f(p),\ \forall\, p \in [x_j,x_{j+1}]$$
and for each $j \geq k$
$$ {{\Hy B(x_j,x_{j+1})}\over{dx}} = 1= \hy f(p),\ \forall\, p \in [x_j,x_{j+1}]$$
Thus, $B$ is admissible. Now let $dy = 2/\Gamma,$ but $\Gamma$ is a *-odd number. Then again we have that $B(x,y) = \int_x^yf\, dx = \int_x^yf\, dy.$ However, there exists $i \in [0,\Gamma -1]$ such that $1$ is the midpoint of $x_i,x_{i+1}$ and $ \Hy B(x_i,x_{i+1}) = dy/2.$ From, this is follows that the (IC) does not hold for this $dy$. \parm
The point $x =1$ in the above example is a point of discontinuity for $f$. If you altered the definition of admissibility to have the $(IC)$ holds for all $dx$ and for all $p_i \in [x_i,x_{i+1}]$ you get a notion I called {{supernearness.}} I show in Herrmann (1993), that an additive function $B$  is supernear to $f$ iff $f$ is continuous. And, of course, the $B$ is equal to the integral. Let's complete this chapter by considering an additional property for an additive function, a property that models various geometric and physical notions.\parm 
{\bf Definition 14.7.} ({\bf Rectangular Property}) An additive  function $A \colon [a,b] \times [a,b] \to \real,$ has the {{\bf rectangular property}} for $f$ iff for any $c,d \in [a,b], c\leq d,\ m(d-c) \leq A(c,d)\leq M(d-c),$ where, as usual, $m = \inf \{f(x) \mid x\in [c,d]\},\ M = \sup\{f(x)\mid x\in [c,d]\}.$ \parm
What does the addition of the rectangular property do for us? By *-transform, consider $dx$. Then for each $k \in [0,x_{\Lambda-1}],$ $ m_kdx \leq \Hy A(x_k,x_{k+1}] \leq M_kdx,$ and if $v_{\Lambda +1} =b, b - x_\Lambda \not= 0,$ then $m_\Lambda(b-x_\Lambda) \leq \hy A(x_\Lambda,b) \leq M_\Lambda(b-x_\Lambda),$ where the $m_i,M_i$ are defined in the usual way. Consequently, in general, for such functions for each $k \in [0,\Lambda-1],$ there is some $p_k \in [x_k,x_{k+1}$ such that 
$$\left|{{\hy A(x_k,x_{k+1})}\over{dx}} - \hy f(p)\right| \leq M_k - m_k$$ and if $b- x_\Lambda \not= 0$, then there exists some $p_\Lambda \in [x_\Lambda,b]$ such that 
$$\left|{{\hy A(x_\Lambda,b)}\over{b-x_\Lambda}} - \hy f(p)\right| \leq M_\Lambda - m_\Lambda.$$ This discussion leads to the following theorem.\parm
{\bf Theorem 14.8.} {\it Assume that $A$ is additive and has the rectangular property for  $f$ and that there exists some $dx$ such that for the simple partition $P(dx)$, whenever $k\in [0,\Lambda -1],$ then $M_k - m_k \in \mu (0)$, where $m_k = \inf\{f(x)\mid x \in [x_k,x_{k+1}]\}, \ M_k = \sup \{f(x)\mid x \in [x_k,x_{k+1}]\}.$ Further, if $b- x_\Lambda \not= 0,$ then $M_\Lambda- m_\Lambda \in \mu (0),\ m_\Lambda = \inf\{f(x)\mid x \in [x_\Lambda,b]\}, \ M_\Lambda = \sup \{f(x)\mid x \in [x_\Lambda,b]\}.$ 
Then $A$ is admissible for $f$.} \parm
Why is the usual application of the integral to functions that are piecewise continuous on $[c,d]$? Well, first of all, since the value of the integral is independent from the value of the function at the end points of the intervals of definition, then all that is needed is to consider why for a specific closed interval $[a,b].$ The next theorem shows why the result in Example 14.6 occurs. \parm
{\bf Theorem 14.8.} {\it A function $f$ is continuous on $[a,b]$ iff
for each $dx$ and, hence, each fine partition $P(dx)$, whenever $k \in [0,\Lambda -1]$, it follows that $M_k - m_k \in \mu (0)$ and if $b- x_\Lambda \not= 0,$ then $M_\Lambda - m_\Lambda \in \mu (0).$}\pars
Proof. Let $f$ be continuous on $[a,b].$ The $f$ is uniformly continuous. So, let's consider any $dx$ and $\vert p - q \vert \leq dx,\ p,q \in \Hy [a,b].$ Then $p-q \in \mu (0)$ implies that $\hy f(p) - \hy f(q) \in \mu (0).$ Consider the simple partition $P(dx).$ Since, for each $h \in [0,\Lambda-1]$ there exist $p,q \in  [x_k,x_{k+1}]$ such that $M_k = \hy f(p),\ m_k = \hy f(q)$ as well as for the case that $b - x_\Lambda \not= 0,$    
 then the necessity follows. \pars
For the sufficiency, let $p-q \in \mu (0),\ p,q \in \Hy [a,b].$ Then there exists $dx$ such that $\vert p - q\vert \leq dx.$ Consider a $P(dx)$ fine partition. First, assume that for some $k \in [0,\Lambda -1],\ p,q \in [x_k, x_{k+1}]$ or that $p,q \in [x_\lambda, b].$ Then since $M_k - m_k \in \mu (0),\ \hy f(p) - \hy f(q) \in \mu (0).$ If there does not exist some $k \in [0,\Lambda]$ such that $p,q \in [x_k,x_{k+1}],$ then $p,q$ are in adjacent intervals by *-transform of the standard case. So consider $2dx = dy$ and apply the first case argument to show that $\hy f(p) - \hy f(q) \in \mu (0).$ Thus, $f$ is (uniformly) continuous on $[a,b]$ and the proof is complete. \qed \parm
{\bf Corollary 14.10.} {\it Let $A$ be additive on $[a,b]$ and have the rectangular property for $f$. If $f$ is continuous on $[a,b],$ then $A(x,y) = \int^y_x f\, dx,\ x\leq y,\ x,y \in [a,b]$ and the function $A$ is unique.}\parm
Thus, if you start with a geometric or physical property that is measured by an additive function $A$ with the rectangular property for a continuous function $f$, then $A$ is uniquely modeled by the integral. On the other hand, although the function $f$ need not be continuous, if $\int_a^b f\, dx \in \real$, then for $x\leq y,\ x,y \in [a,b]$ the function $A(x,y) = \int_x^yf\, dx$ is additive and has the rectangular property. \par\vfil\eject    
 \centerline{\bf 15. GENERALIZATIONS}
\parm
Much of what I've covered can be highly generalized. It should be obvious that this nonstandard approach, although very restricted in the language used, is not depended upon the codomain of the set of sequences used to obtain the equivalence classes with respect to a free ultrafilter. Hence, the set $\real$ and all the additional ones used in $\cal M$ can be replaced with the set $X \cup \real,$ where $X$ is non-empty. The *-transform method holds and many of the general results, such as Theorem 3.2, follow since they are all obtained from the properties of the ultrafilter. It would be better to have a stronger language and a structure where we can use $\in$ over variables. But, even in our restricted language, much can be done. I give just a few brief example. \parm
{\bf Definition 15.1.} ({\bf Real Metric Space}.) A  nonempty set $X$ is called a {{\bf metric space}} iff there exists a function $d \colon X\times X \to \real$ such that
for each $x,y,z \in X,$  \pars
\indent\indent (i) $ d(x,y) = d(y,x) \geq 0,$\pars
\indent\indent (ii) $d(x,y) = 0$ iff $x = y,$\pars
\indent\indent (iii) $d(x,y) \leq d(x,z) +d(z,y).$\parm
{\bf Definition 15.2.} ({\bf General Finite Points and Monads}.) Any $q \in \hy X$ is {{\bf finite}} iff there is some $p \in X$ such that $\hy d(q,p) \in G(0).$ For each $p \in X,$ the {{\bf monad}} of $p$ is $\mu (p)=\{x\mid (x \in \hy X)\land (\hy d(x,p) \in \mu (0))\}= \{x\mid(x\in \hy X)\land \forall r((r\not= 0)\land (r \in \real)\to \hy d(x,p)< \vert r\vert)\}$. The set ns$(\hy X) = \bigcup\{\mu (p)\mid p \in X\}.$ \parm
In this more general case, what was previously the set $G(0)$ is now denoted by ${\rm fin}(\hy X)$, the set of all finite points in $\hy X.$ And, as before, ${\rm ns}(\hy X) \subset {\rm fin}(\hy X)$. These sets are equal in the case that $X = \real$, but for metric spaces in general they are not equal. \pars
 A closed sphere about $p \in X,$ $S[p,r] = \{x\mid d(x,p)\leq r\}.$ A set      $B \subset X,$ for the metric space $(X,d),$ is {{\bf bounded}} iff there is some closed sphere $S[p,r]$ such that $B \subset S[p,r]$. Now I assume that the theorem on the *-transform has been established for our structure.\parm
{\bf Theorem 15.3.} {\it For the metric space $(X,d)$, $B\subset X$ is bounded iff $\Hy B \subset {\rm fin}(\hy X).$}\pars
Proof. For the necessity, the sentence
$\forall \r x ((x \in S[p,r] \to d(x,p)\leq r)$ holds in $\cal M$; and, hence, in $\hy {\cal M}.$ Thus, by *-transform, 
$\forall \r x ((\r x \in \Hy (S[p,r]) \to \hy d(\r x,p)\leq r).$
Consequently, $\Hy B \subset \Hy (S[p,r]) \subset {\rm fin}(\Hy X).$\pars
For the sufficiency, assume that $B \subset X$ is not bounded. Let $p \in X$. Then the sentence
$\forall \r x((x\in \real^+) \to \exists \r y((\r y \in B) \land (d(p,\r y) > \r x))$
holds in $\cal M.$ Thus, by *-transform, letting $\Lambda \in \real_\infty^+$, then there exists some $q \in \Hy B$ such that $\hy d(q,p) > \Lambda.$ Now let $p' \in X$. Then  $\hy d(q,p')$ cannot be a finite *-real number. For if we assume that $\hy d(q,p') \in G(0),$ then since $d(p,p') \in G(0)$ we would have that $\hy d(p,q) \leq d(p,p') + \hy d(p',q) \in G(0);$ a contradiction. Hence, $q\notin \Hy (S[p',r])$ for any $r \in \real^+$ and any $p\ \in X.$ This completes the proof. \qed\parm
{\bf Corollary 15.4.} {\it A sequence $S\colon \nat \to X$ is bounded iff 
$\hy S(\Lambda) \in {\rm fin}(\hy X)$ for each $\Lambda \in \nat_\infty.$}\parm
The following results are obtained immediately in the same manner as the corresponding real number results. \parm
{\bf Theorem 15.5.} {\it For a metric space $(X,d)$, a sequence $S\colon \nat \to X$ converges to $L$ iff $\hy S(\Lambda) \in \mu (L)$ for each $\Lambda \in \nat_\infty.$}\parm
{\bf Theorem 15.6.} {\it For a metric space, every convergent sequence is bounded.}\parm
{\bf Theorem 15.7.} {\it A point $p \in X$ is an accumulation point for a sequence $S \colon \nat \to X$ iff there exists some $\Lambda \in \nat_\infty$ such that $\hy S(\Lambda) \in \mu (p).$}\parm
{\bf Corollary 15.8.} {\it A sequence $S \colon \nat \to X$ has a convergent subsequence iff there exists some $\Lambda \in \nat_\infty$ such that $\hy S(\Lambda) \in {\rm ns}(\hy X).$}\parm
{\bf Theorem 15.9.} {A sequence $S\colon \nat \to X$ is Cauchy iff $\hy d(\hy S(\Lambda),\hy S(\Omega)) \in \mu (0)$ for each $\Lambda, \Omega \in \nat_\infty.$}\parm
It's possible for metric space, including the real numbers, to define monads at points $q \in \hy X - X$ by letting $\mu (q) = \{x\mid \hy d(x,q)\in \mu (0)\}.$ Then it follows that $S \colon \nat \to X$ is Cauchy iff there exists some $q \in \hy X$ such that $\hy S(\Lambda) \in \mu (q),\ \forall\,  \Lambda \in \nat_\infty.$ One of the most significant metric spaces is the {{\bf normed linear (vector) space}}. I consider a linear space over the real numbers for my example. If $V$ is a linear space over the real numbers, then a norm is a  map $\Vert \cdot \Vert \colon V \to \real$ with the properties that, for each $x,y \in V$, (i) $\Vert x \Vert \geq 0.$ (ii)  
For each $r \in \real$, $\Vert rx \Vert = \vert r\vert\, \Vert x \Vert.$ (iii) $\Vert x + y \Vert \leq \Vert x\Vert + \Vert y \Vert.$ \pars
The metric is defined by letting $d(x,y) = \Vert x - y \Vert.$ Then you now apply nonstandard analysis to this space along with its additional linear space properties. For example, we have that $\mu (p) = \{p + \gamma \mid \gamma \in \mu (0)\}$ for such a metric space, in general. Nonstandard analysis has been applied extensively to linear spaces. For the major generalization known as the {{\bf topological spaces}}, where I have established some immediately of the original results, we need a structure more directly related to set-theory and such an appropriate structure is not what I would consider as elementary in character.\par\vfil\eject 
\centerline {\bf APPENDIX}\parm
{\bf Theorem A1.} {\it Let $\cal F$ be any filter on $X$. Then there exists an ultrafilter ${\cal U}_X \supset {\cal F}.$}\pars
Proof. You can either use the set-theoretic axiom that states that this statement holds; or use Zorn's Lemma, which is equivalent to the Axiom of Choice. Let $\cal G$ be the set of all filters that contain $\cal F.$ Suppose that $\cal C$ is a chain with respect to $\subset$ in $\cal G$. I show that $\bigcup {\cal C}$ is a filter that obviously would be a upper bound for this chain and is contained in $\cal G.$ Clearly, $\emptyset \notin \bigcup {\cal C}.$ Let $A \in  \bigcup {\cal C}.$ Then $A \in {\cal F}_1$ for some ${\cal F}_1 \in \bigcup {\cal C}.$ Hence, if $A\subset B$, then $B \in {\cal F}_1$ implies that $B \in  \bigcup {\cal C}.$ Now let $A,B \in \bigcup {\cal C}.$ Since ${\cal G}$ is a chain with respect to $\subset$, then there is some ${\cal F}_2 \in {\cal G}$ such that $A,\ B \in {\cal F}_2.$ Hence, $A \cap B \in {\cal F}_2$ implies that $\bigcup {\cal C}$ is a filter that contains $\cal F$. By Zorn's Lemma, there is a member of $\cal G$ that is a maximal member $\cal U$ with respect to $\subset.$ If $\cal U$ is not an ultrafilter, then there would be a filter ${\cal U}_1$ not equal to $\cal U$ such that ${\cal U}\subset {\cal U}_1.$ But, then ${\cal U}_1 \in {\cal G};$ a contradiction of the ($\subset$) maximal property for $\cal U.$ This completes this proof. \qed\parm 
I assume that the reader knows what I mean by a first-order language $\cal L$ with equality, where the constants represent objects in $\real,\real^2,\ldots $, and $\power {\real}, \power {\real^2},\ldots.$ Equality is interpreted to be the identity on $\real$ or set-theoretic equality elsewhere. The variables are denoted by Roman font. The first class of atomic formula are $\r x \in Y$ and, for $n>1$, $( x_1,\ldots, x_n)\in Y,$ where $Y$ is a constant, and all possible permutations of members of the $n$-tuples, where $x_k$ is either a constant or variable. I leave to the reader the trivial cases where the various  expressions only contain constants and assume that, for all other formula, at least, one of the symbols that can differ from a constant is a variable. (The $Y$ includes the $+, \cdot,\leq$.  Our result below holds for many other collections of atomic formula that describe members of more comprehensive structures, but I don't use them for this monograph.) Finally, $a = b$, where $a, b $ are both variables or, at most, one is a constant. (Note: the symbols $=$ is interpreted as a special binary relation within our language). Only a special set $K$ of formula built from these atomic formula is used. Further, $P \in K$ if and only if $P$ has only bounded quantifiers. That is each quantifier contained in the $P$ is restricted to subsets of $\real$ or $\real^n.$ Indeed, in this book, all $P$ with bounded quantifiers are equivalent to a form $\forall \r x ((\r x \in X) \to \cdots)$ or the form $\exists \r x((\r x \in X) \land \cdots).$ The reason I'm using these special forms is that the *-transform (Leibniz) property for such formula can be established without the Axiom of Choice. Given any $P,$ then $\Hyper P$ is obtained by placing every constant $A$ in $P$ by $\hyper A.$ (Note: The $+,\cdot,\leq$ are constants that technically should carry the $\hyper {}$ notation, but it's customary to drop this notation when the context is known.) In what follows, $(a)= a$ and $x \in \real$ is considered in two context, either a constant for a member of $\real$ or as varying over a subset of $\real$ as the case may be. 
Let $K$ be our set of formula and $${\cal M} = \langle \real,\ldots, \real^n,\ldots,\power {\real},\ldots, \power {\real^n},\ldots,+, \cdot, \leq \rangle$$ $$\hy {\cal M} = \langle \hyperreal,\ldots, \hyperreal^n,\ldots,\power {\hyperreal},\ldots, \power {\hyperreal^n},\ldots,+, \cdot, \leq \rangle.$$ 

{\bf Theorem A2.} {\it Let $P(x_1,\ldots,x_p) \in K$ contain at least one variable and $X$ is a member of $\cal M$.  Define $ A=\{(x_1,\ldots,x_p)\mid ((x_1,\ldots,x_p) \in X) \land (P {\rm \ holds\ in\ }{\cal M})\}.$ Then 
$$\hyper A = \{(x_1,\ldots,x_p)\mid ((x_1,\ldots,x_p) \in \hyper X) \land  (\hyper P {\rm \ holds\ in\ } \hy {\cal M}) \}.$$ }\pars
Proof. Let $P= (\r x \in Y),$ where $Y \subset \real^n.$ Then $A = \{x\mid (x\in X)\land (x\in Y)\} = \{x\mid x \in X\} \cap \{x\mid x \in Y\}.$ By Theorem 3.2 (vi)(xi), $\Hyper A = \hyper X \cap \Hy Y =\{x\mid (x \in \hyper X)\land (x \in \Hy Y)\}.$ Let $P = (\r x = \r y),$ (or $P = (\r x = \r x)).$  Let $X \subset \real,$ $X\times X = Y$ and $A =\{( x,y)\mid ((x,y) \in Y\land (x=y)\}.$ The result follows from $\hy X \times \hy X = \Hy Y$ and if $a,b \in X$ and $a=b,$ then $[A] = [B].$ For the case that $X \subset \real^n,\ n >1,$ we need either to identify $X \times X$ with the obvious 2n-ary relation or we need to extend the structure to include such objects and extend the results in Theorem 3.2 to cover these objects. This also depends upon your definition of the n-tuple.  The case where $x=a$ or $a=x$ also follows in like manner. Note that since function or term symbols are not used in the language, then the $=$ can be considered as used to generate specific relations that are elements of our structure. \pars

For atomic formula $(\r x_1\ldots \r x_p)\in Y$, the result follows by application of Theorem 3.2 and the definition of $\Hy Y$ noting, of course, that  $(x_1,\ldots,  x_p)\in X.$  \pars
Since every first-order formula is equivalent to a formula that has all of the quantifiers to the left of a formula that contains no quantifiers but only formula built from atomic formula and the connectives  $\land,\lor,\to,\iff,\neg.$ As for the connectives, it is well know that if we assume that the result holds for quantifier free formula $V$ and $W,$ then all we need to do to show, by induction, that the result holds in general for quantifier free formula is to show that it holds for $V\land W$ and for $\neg V.$ For all but the atomic formula $(\r x_1,\ldots, \r x_p)\in Y$ this is immediate by Theorem 3.2. Now let $V$ be the expression $(\r x_1,\ldots, \r x_p)\in Y$ and 
$\hy A = \{( x_1,\ldots,x_p)\mid (x_1,\ldots,x_p)\in \Hy X)\land (( x_1,\ldots,x_p) \in \Hy Y)\}.$ Then consider $\Hy B = \{(x_1, \ldots,x_p)\mid (x_1,\ldots, x_p) \in \Hy X)\land ((x_1,\ldots,x_p) \notin \Hy Y)\}.$ The *-transfer holds since $\Hy B = \hy X - \Hy A$ from Theorem 3.2.\pars

 Now let $V = (\r x_1,\ldots,\r x_p, \r y_1,\ldots, \r y_q)\in Y,\ W = (\r x_1,\ldots,\r x_p, \r z_1,\ldots, \r z_r)\in Z,$ where I assume the possibility that both $V$ and $W$ contain  $\r x_1,\ldots, \r x_p$ and the other constants or variables are distinct from these. Let $(x_1,\ldots,x_p,y_1,\ldots,y_q,z_1,\ldots,z_r) = (x_1,\ldots,z_r).$ Let $B = \{x_1,\ldots,z_r)\mid (x_1,\ldots,z_r) \in X) \land ((x_1,\ldots,x_p, y_1,\ldots, y_q)\in Y)\}$ and $C = \{x_1,\ldots,z_r)\mid (x_1,\ldots,z_r) \in X)\land (( x_1,\ldots,x_p, z_1,\ldots, z_r)\in Z)\}.$ Then 
$A =\{(x_1,\ldots,z_r)\mid (x_1,\ldots,z_r) \in X) \land (V \land W)\} = B \cap C.$  The result holds from the induction hypothesis, in this case, since $\hy A = \Hy B \cap \Hy C.$ \pars

As mentioned, any first-order formula is equivalent to one which can be written as $V = (q\r x_{n+1})\cdots (q\r x_1)W,$ where $\r x_1,\ldots, \r x_{n+1}$ are free variables in $W$ and $W$ is a finite combination
via $\land$ and $\neg$  of all of our quantifier free atomic formula. Hence, represent this formula by $W( \r y_1,\ldots, \r y_p, \r x_1,\ldots,\r x_{n+1}).$ We can always assume that $(q\r x_{n+1}) V= (\exists \r x_{n+1})V$ (for if not, consider $\neg V$) and $V$ is also in this special quantifier form. If $n =0,$ then the result has been established. Assume the result holds for an appropriate member of $K$ with the number of quantifiers $\leq n.$ Under our requirements, $x_{n+1}$ is restricted to a member $Z$ of our structure. Let $D = \{(y_1,\ldots,y_p), x_{n+1})\mid ((y_1,\ldots,y_p),x_{n+1}) \in X \times Z\land (qx_n\ldots qx_1)W\},$ where $X \subset \real^p$ is also in the standard structure. Then, by induction, and Theorem 3.2, 
$\Hy D = \{(y_1,\ldots,y_p), x_{n+1})\mid ((y_1,\ldots,y_p),x_{n+1}) \in \hy X \times \hy Z\land (qx_n\ldots qx_1)\Hy W\}.$ Let 
$A= \{(y_1,\ldots,y_p)\mid ((y_1,\ldots,y_p)\in X)\land(\exists x_{n+1}((x_{n+1} \in Z)\land (qx_n\ldots qx_1)W))\}.$ Using this and a simple modification of the proof of Theorem 3.2 (x), it follows that the domain of $\Hy D= \hy A,$ where 
$\hy A = \{(y_1,\ldots,y_p)\mid ((y_1,\ldots,y_p)\in \hy X)\land(\exists x_{n+1}((x_{n+1} \in \Hy Z)\land (qx_n\ldots qx_1)\Hy W))\}= \{(y_1,\ldots,y_p)\mid ((y_1,\ldots,y_p) \in \Hy X)\land \Hy V\}.$ By induction this completes the proof. \qed\parm
{\bf Theorem A3.} {\it Let $V\in K$ be a sentence with necessary quantifiers or be compose only of connected atomic formula expressed only in constants. 
Then $V$ holds in $\cal M$ iff $\Hy V$ holds in $\hy {\cal M}.$}\pars
Proof. Since $V$ is a sentence it has no free variables. If $V$ contains no quantifiers, then $V$ only contains constants and the result follows from the definition of the hyper-extension and Theorem 3.2. \pars
Now assume that $V$ contains quantifiers and that it is written in the equivalent form (prenex normal form) $V = (q\r x_n\cdots q\r x_1)W$ and $q\r x_n = \exists \r x_n$. (If this is not the case, consider the negation.) To say that $V$ holds in $\cal M$ means that $A =\{x_n\mid (x_n \in X)\land((q\r x_{n-1}\cdots q\r x_1)W\ {\rm holds\ in\ }{\cal M})\} \not= \emptyset,$ where $X$ is the domain for $\exists \r x_n.$ But $A \not=\emptyset$ iff $ {\emptyset} = \Hy {\emptyset} = \hy A = \{x_n\mid (x_n \in \Hy X)\land((q\r x_{n-1}\cdots q\r x_1)\Hy W\ {\rm holds\ in\ }\hy {\cal M})\}.$ This completes the proof. \qed \parm

\centerline{\bf REFERENCES}
\parm
\id{B}ahrens, M. (1972), A local inverse function theorem, in {\it Victoria Symposium on Nonstandard Analysis,} (Lecture Notes in Mathematics \#369, Springer-Verlag, NY, pp. 34-36.\parm
\id{B}urrill, and Kundsen, (1969), {\it Real Analysis}, Holt, Rinehart and Winston, NY.\parm
\id{H}errmann, R. A. (1991), {\it Some Applications of Nonstandard Analysis to Undergraduate Mathematics - Infinitesimal Modeling, Elementary Physics, and Generalized Functions}, Instructional Development Project, Math. Dept., U. S. Naval Academy, Annapolis, MD. http://arxiv.org/abs/math/0312432 \  http://www.serve.com/herrmann/cont2s.htm.\parm
 \id{J}ech, T. J. (1971), {\it Lecture Notes in Set Theory,} Lectures Note In Mathematics \# 217, Springer-Verlag, NY.\parm
\id{L}uxemberg, W. A. J. (1962), {\it Non-Standard Analysis - Lectures on A. Robinson's Theory of Infinitesimals and Infinitely Large Numbers,} Math. Dept., California Institute of Technology Bookstore, Pasadena, CA.\par
\id{R}obinson, A. (1966), {\it Non-standard Analysis,} North-Holland, Amsterdam.\parm
\id{R}obinson, A. (1961), Non-standard analysis, {\it Nederl. Akad. Wetensch. Proc. Ser. A62, and Indag. Math,} {\bf 23}:432-440.\par
\id{R}udin, W. (1953), {\it Principles of Mathematical Analysis,} McGraw Hill, NY. \parm
\id{S}uppes, D. (1960), {\it Axiomatic Set Theory,} D. Van Nostrand, NY.     
\end

\end

\end